\documentclass[a4paper,10pt]{article}
\usepackage[utf8]{inputenc}
\usepackage[T1]{fontenc}
\usepackage{algorithm}
\usepackage{algpseudocode}
\usepackage{amsmath}
\usepackage{enumitem} 
\usepackage{amssymb}
\usepackage{framed}
\usepackage{float}
\usepackage{caption}
\usepackage{geometry}
\usepackage[labelformat=simple]{subcaption}
\usepackage{bm}
\usepackage{multirow}
\usepackage{array}
\usepackage[english]{babel}
\usepackage{blindtext}
\usepackage{booktabs}
\usepackage{pgfplots}
\usepackage{pgfplotstable}
\usepackage{tikz}
\usepackage{array}
\usepackage{authblk}
\usetikzlibrary{matrix}
\usepgfplotslibrary{groupplots, external}
\pgfplotsset{compat=newest}
\pgfplotstableset{
  header = false,
}
\pgfplotsset{
             every legend/.append style={font = \Large}
             }
\newcommand{\bv}[1]
{
  \bm{#1}
}

\makeatletter
\renewenvironment{abstract}{%
      \@beginparpenalty\@lowpenalty
      \small
      \begin{center}%
        \bfseries \abstractname
        \@endparpenalty\@M
      \end{center}}%
     {\par}
\makeatother
\begin{document}
\title{Fast and Accurate Proper Orthogonal Decomposition using Efficient Sampling and Iterative Techniques for Singular Value Decomposition  } 
\author[1]{V. Charumathi}
\author[1]{M. Ramakrishna}
\author[1]{Vinita Vasudevan}
  \affil[1]{Indian Institute of Technology, Madras}
  \date{}
\maketitle
\begin{abstract}
  In this paper, we propose a computationally efficient iterative algorithm for proper orthogonal decomposition (POD) using random sampling based techniques.
  In this algorithm,  additional rows and columns are sampled and a merging technique is used to update the dominant POD modes in each iteration.  
  We derive bounds for the spectral norm of the error introduced by a series of merging operations.
  We use an existing theorem to get an approximate measure of the quality of subspaces obtained on convergence of the iteration. 
  Results on various
 datasets indicate that the POD modes and/or the subspaces are approximated with excellent accuracy with a significant runtime improvement over computing the truncated SVD.
We also propose a method to compute the POD modes of large matrices that do not fit in the RAM using this iterative sampling and merging algorithms.

 \end{abstract}

\maketitle
\section{Introduction}
Proper orthogonal decomposition (POD), also known as Principal component analysis (PCA) or Karhunen-Loeve transform is a dimensionality reduction technique used in a variety of applications including  information compression, pattern mining, clustering and classification,
facial recognition, missing data estimation, data-driven modelling and Galerkin projections.
PCA aims to find an optimal set of basis vectors that captures energetically dominant features of the dataset.
These basis vectors are  also  referred to as eigenfeatures or POD modes or principal directions/components (PCs).
The assumption is that a few of these modes are sufficient
to capture all significant features of the dataset  i.e., the data is highly correlated  and inherently low rank.  In many of these cases, the dataset
is also large, often distributed over several computers and dense.
The two methods commonly used to compute POD modes are via (a) eigendecomposition of $A^TA$  and (b) SVD of $A$. In both cases, a full 
decomposition is expensive and approximate methods are used to obtain the dominant modes. In this paper, we explore the use of random sampling based techniques to obtain
the POD modes accurately and in a computationally efficient manner. The way our data is arranged, these modes are the left singular vectors of the data matrix.

We use the following notation throughout the paper, unless otherwise specified. 
Lower case letters are used to denote scalars.
Lower case bold-face letters are used to denote a vector.
Upper case letters are used to denote a matrix. 
The $i^{th}$ element of vector $\bv{a}$ (bold font) is denoted as $a_i$ (non-bold font). The element in the $i^{th}$ row and $j^{th}$ column
of the matrix $A$ is denoted by $a_{ij}$.
$\bv{a}^i$ denotes the $i^{th}$ row vector  and $\bv{a}_j$ denotes the  $j^{th}$ column vector of matrix $A$. 
 $A_k$ denotes the best rank-$k$ approximation of a matrix $A$  and the approximation to $A_k$ is denoted using tilde, as $\tilde{A}_k$.
Let $A_k = U_k\Sigma_k V_k^T$ and $\tilde{A}_k=\tilde{U}_k\tilde{\Sigma}_k\tilde{V}_k^T$ be the SVDs of 
$A_k$ and $\tilde{A}_k$. If $A$ is obtained by mean centering rows of the data matrix, then $U_k$ and $\tilde{U}_k$ are the accurate and the approximate principal components (PCs or POD modes), respectively.
$||A||_p$ denotes the $p$ norm of matrix $A$. $A$ is assumed to be an $m \times n$ matrix.

Algorithms that approximate $A_k$ can be classified as random projection, random sampling based algorithms or a combination of both.
Random projection based algorithms have a multiplicative error bound \cite{HalMarTro:2011,Sar:2006}.
They have been used for PCA in \cite{RokSzlTyg:2009,NguDoTra:2009,HalMarShkTyg:2011,EriBruKut:2017} 
with algorithmic optimizations to improve accuracy and numerical stability \cite{RokSzlTyg:2009,Gu:2015,EriVorBruKut:2016,LiLinSzlStaKluTyg:2017} and with out-of-core modifications in \cite{HalMarShkTyg:2011,BosKalKonElkPasDri:2019} for large matrices.
\cite{YamTomDon:2017} use random projections to compute the SVD and random sampling (of rows) of additional data to incrementally update the SVD.
\cite{DriMah:2016} also use a combination of projection and sampling techniques, where projection is used to precondition the matrix so that subsequent
uniform sampling works well.
 
 Random sampling methods are of two types: those that sample elements of the matrix and
 those that sample columns/rows of the matrix. 
  The motivation behind sampling elements \cite{AchMcS:2007} from the matrix is to make the matrix sparser, thereby reducing the time taken to compute SVD of the matrix.
 Sampling columns/rows is usually done using (a) uniform probability \cite{DriDriHug:2001,YamTomDon:2017} (b) leverage scores \cite{DriMahMut:2008,LiMilPen:2013,CohLeeMusMusPenSid:2015,YamTomDon:2017} (c) column and/or row norms \cite{DriDriHug:2001,DriFriKanVemVin:2004,FriKanVem:2004,DriKanMah:2006:II}. In all three cases, error bounds are derived assuming sampling is done with replacement.
 Uniform and row/column norm based sampling have additive error bounds, while leverage score based 
 sampling has a multiplicative error bound. However, leverage score based sampling is  computationally more involved since it requires an estimate of the singular vectors.
While row/column norm based sampling require a maximum of two passes over the data to get the sampling probability distribution,
uniform sampling requires none.
But it can be inefficient as it is equally likely to sample unimportant columns or rows, as for example, zero columns/rows.

The results of the experimental evaluation done by \cite{MenElk:2011} indicates the column norm sampling method has low run times, with accuracies well below the error bounds for standard datasets.
Additionally, random sampling of columns/rows has several advantages - it retains the sparsity of the matrix and is easily amenable to incremental improvement.
Multiplicative error bounds using row/column norm based sampling can be obtained using
 adaptive sampling, which has several rounds of sampling based on row/column norms \cite{DesRadVemWan:2006,DesVem:2006}. 

In this paper, our focus is to obtain POD modes of large datasets accurately and in a computationally efficient manner using random sampling based 
techniques. We evaluated the performance of sampling algorithms proposed in \cite{DriKanMah:2006:II}.
Even with tight error parameters, the top $k$ POD modes computed using these sampling algorithms as well as the subspace spanned by these vectors have significant errors. To improve accuracy,
  we propose an iterative technique in which additional rows and columns are sampled and the dominant POD modes are updated using a merge-and truncate 
(MAT) operation \cite{VasRam:2017} in each iteration. The iteration converges when the angle between the current and updated modes is less than 
a threshold. 
  We also derive bounds for the spectral norm of the error introduced by a series of MAT operations. We use an existing theorem to get an approximate measure of the quality of subspaces obtained after convergence.
  Results on various
 datasets indicate that the POD modes and/or the subspaces are approximated with excellent accuracy with a significant runtime improvement over computing the truncated SVD.
Using this iterative algorithm, we also demonstrate a method to compute POD modes of large matrices that do not fit in the RAM of the system.

The paper is organised as follows. In section~\ref{sec:prev_algo}, we discuss a few relevant sampling algorithms found in literature and their drawbacks. 
This motivates the need for an algorithm with more efficiency and better accuracy.
Section~\ref{sec:iterative} contains the proposed iterative algorithm, the error bound for MAT operations and the 
quality measure.
In section~\ref{sec:results} we discuss the results of these algorithms as applied to different datasets.
Section~\ref{sec:conclusion} concludes the paper.
Appendix~\ref{app:eval} contains the experimental results for the datasets  obtained using row/column norm based sampling algorithms (LTSVD, CTSVD) proposed in \cite{DriKanMah:2006:II}.

The following datasets were used for testing all the algorithms.
\begin{enumerate}
	\item Faces data from \cite{ORLfaces:2002}, a collection of 400 images from 40 subjects, each of size $112\times 92$ giving 400 images of 
		dimension 10304 each. The matrix is of size $10304\times400$. \label{data:ORLfaces} This dataset will be referred to as Faces dataset in the paper.
	\item Velocity vector data of a flow generated using a CFD simulation on a grid of size $257\times257$. The snapshot vectors are of dimension 132098 and there are 1024
		of them giving a matrix of size $132098\times 1024$.
    It will be referred as $V_{2D}$
   \item 2414 images of 38 subjects cropped to show only the face of the subject in different illumination conditions 
	   from \cite{LeeHoKri:2005}. The images are of size $192\times168$ giving 2414 images of dimension 32256 each. Therefore, the  matrix size is $32256 \times 2414$.
		This dataset is referred to as cropped Yale faces (CYF).
   \item Faces data from \cite{Yalefaces:2001}, a collection of 16128 images from 28 subjects in different poses
	   and illumination conditions.  18 of the 16128 images in the database were discarded since they were corrupted.
		Each image is of size $480\times 640$. We have 16110 images of dimension 307200 each organised as a $307200\times16110$ matrix.
		This dataset is referred to as Yale faces (YF).
   \end{enumerate}
All the datasets are mean-centered row-wise. The principal components are the left singular vectors.
It is seen that all the datasets have more rows than columns (``tall and thin'').
The first three datasets are small and used for validation of the algorithms. The last dataset is large and the mean-centered data requires 40GB of RAM.

\section{Evaluation of existing sampling algorithms and motivation}\label{sec:prev_algo}
Algorithms to compute low rank approximations of a matrix using sampled columns (LTSVD)/rows and columns(CTSVD) have been proposed in \cite{DriKanMah:2006:II}. 
The sampling is done with replacement and the sampled rows and columns are scaled to preserve the Frobenius norm of the matrix. For both algorithms, the authors show 
that bound for the error in the rank-$k$ approximation 
\begin{equation}
  \label{eqn:LTSVD_err_bound}
  ||A-\tilde{U}_k\tilde{U}_k^TA||_F^2\le ||A-A_k||^2_F + \epsilon||A||_F^2
\end{equation} holds with a probability of at least ${(1-\delta)}$.
The sampling probability of the $i^{th}$ column, $p_i$ and the number of columns sampled, $c$, in LTSVD is
\begin{equation}
  \label{eqn:LTSVD_prob}
  p_i = \frac{|\bv{a}_i|^2}{||A||_F^2}
\end{equation}
\begin{equation}
\label{eqn:c_LTSVD}
    c=4k(1+\sqrt{8\log{(1/\delta)}})^2/\epsilon^2
\end{equation}

In the case of CTSVD, the sampling probabilities for columns and rows ($p_i$ and $q_i$)  are
\begin{equation}
  \label{eqn:CTSVD_prob_col}
  p_i = \frac{|\bv{a}_i|^2}{||A||_F^2}, ~~~ q_j = \frac{|\bv{d}^j|^2}{||D||_F^2},
\end{equation}
Here, $D$ is the matrix containing the sampled and scaled columns.
The number of rows and columns sampled ($c$ and $w$) that need to be sampled are
\begin{equation}
    c=w=k^2(1+\sqrt{\log{(2/\delta)}})^2/\epsilon^4\label{eqn:c_w_CTSVD}
\end{equation}
$\epsilon$ and $\delta$ are the error parameters and the failure probability that need to be set.
An experimental evaluation of the algorithm by \cite{MenElk:2011} concluded that, in practice, the error is much lower than the bound given in equation~(\ref{eqn:LTSVD_err_bound}).

In this paper, our focus is accurate computation of the principal components, which are the left singular vectors. Therefore, we implemented the two algorithms in order to evaluate its performance with respect to the error in the singular vectors. In order to evaluate the accuracy of the modes, we used the following metrics.
\begin{enumerate}
\item The angle between the accurate and approximate left singular vectors ($\theta_i = \cos^{-1}\left( u_i^T\tilde{u}_i \right)$), denoted by mode angle.
\item The principal angles between the rank-$k$ subspaces ($\phi_i$). This is computed as the inverse cosine of the singular values of the matrix $U_k^T\tilde{U}_k$
\cite{BjoGol:1973}.
\end{enumerate}
The pseudo-codes for our implementation and the detailed results are included in Appendix A.
A summary of the results is as follows.
\begin{enumerate}
	\item When $k$ (and hence number of columns sampled) is small, the error in the
singular values can be as high as 20\% for some datasets. For larger values of $k$, this error is of the order of 5\% or less. (see Fig.~\ref{fig:relsigDiff_LTSVD_CTSVD} in Appendix~\ref{app:eval}).
\item There is a large error in the computed POD modes.  As seen in Fig.~\ref{fig:POD_2}, if  the mode angle is greater than about $10^{\circ}$, the POD modes have 
erroneous features. The mode angles obtained using these algorithms are as large as $80^{\circ}$ in some 
cases (see Fig.~\ref{fig:angPOD_LTSVD_CTSVD}). As expected, the accuracy of the dominant subspaces is better, with principal angles 
between subspaces less than $40^{\circ}$ in most cases (see Fig~\ref{fig:pc_LTSVD_CTSVD}). This large error occurs because each column is scaled by a factor that depends on the column norm. While this preserves the Frobenius norm, it tends to distort the subspace spanned by the columns.
There is no significant improvement in accuracy, even when the error parameters are tightened.
\end{enumerate}
  
One major problem with the algorithms is that, it is hard to a priori fix the values of the error parameters $\epsilon$ and $\delta$ (see equations~(\ref{eqn:c_LTSVD}) and 
(\ref{eqn:c_w_CTSVD})). Some of these issues are resolved by the
adaptive sampling algorithm proposed in \cite{DesVem:2006}. This algorithm also has a tighter multiplicative error bound.
It is based on approximate volume sampling which essentially involves multiple rounds of row sampling. In each round, the probability of picking a row is 
proportional to the squared distance from the span of previously sampled rows.

The steps involved in the algorithm are as follows.
\begin{enumerate}
  \item Sample $k$ rows from $A$ based on the approximate volume sampling (see \cite{DesVem:2006} for details) and compute the basis, $Q$, spanning these rows.
  \item Find the space orthogonal to the space spanned by these basis vectors, ${E = A-AQQ^T}$.\label{itm:error_norm}
  \item Sample $s$ rows from $A$ based on row norms of $E$ and compute $Q$, the basis spanning all sampled rows.  
  \item Repeat from (\ref{itm:error_norm}) for $(k+1)\log_2(k+1)$ times.
\end{enumerate}
  In each round $2k$ rows are sampled, except for the last round where $16k/\epsilon$ rows are sampled.
  For $k=10$ and $\epsilon = 0.5$, we would need $38$ rounds of sampling and 1060 rows will be sampled to find the basis.
  Once the iteration is over, $Q$ is an approximation to the basis for the row space. To
  obtain the approximate left-singular vectors, we need to additionally find the SVD of $AQ$. For the $V_{2D}$ dataset, we see that $AQ$ is approximately
  the same size as $A$ and it is larger for the Faces dataset (for which you anyway cannot get more than 400 basis vectors for the column space). 

  The advantage of this 
method is that (a) the error bound is multiplicative and (b) it directly approximates the dominant left and right subspaces. The sampled
  rows are not scaled, thus giving a better approximation to the dominant right singular vectors. 
 
The algorithm can be made more efficient if (a) there is some pruning of basis vectors in $Q$ in each iteration, depending on whether it belongs to the dominant 
subspace or not and (b) computation of $E=A-AQQ^T$ is avoided. Since the memory required by $AQQ^T$ is the same as $A$, the step becomes inefficient when the matrix sizes are large.

\section{Proposed Iterative Sampling Algorithm (ISMA)}\label{sec:iterative}

  In our algorithm, for improved computational efficiency and accuracy, we would like to have several rounds of column or column and row
  sampling rather than just row sampling. More importantly, in each round, we would also like to prune the basis vectors that are not in the dominant
  subspace. We would like to stop the iteration when additional sampling does not improve the quality
  of either the modes or the subspace spanned by the modes. All sampling is done with replacement and duplicates are removed.
 
We use the following steps in our iteration. 
\begin{enumerate}
\item Sample a set of columns from  $A$ based on column norms. Optionally, also sample rows. We use equations~(\ref{eqn:c_LTSVD}) and (\ref{eqn:c_w_CTSVD}) to estimate the
number of columns and rows to be sampled, respectively.
		Compute SVD  of the sampled matrix. Truncate to get the dominant POD modes.
	\item From the remaining columns, sample an additional set of columns or columns and rows. For this, we explore various sampling strategies. 
		These are discussed in section~\ref{sec:samp_strat}. Find distinct columns (optionally, if rows are sampled, remove duplicates and scale 
		remaining rows), compute POD modes of the newly sampled matrix and use it to update the dominant modes.

	\item Find the cosines between the previously computed and newly updated dominant modes or the principal cosines between the subspaces.
  If any one of the cosines is less than the desired value, $\tau$, continue the iteration. 
\end{enumerate}
\begin{algorithm}[!htbp]
  \caption{ISMA; 
    Inputs: matrix $A$, 
  sampling probability $\bv{p}$, 
   rank of approximation $k$,
  number of columns to be sampled $c$ ,
   rank of matrices after BLOCK MERGE $r$, $rows$ (=1  if rows are also to be sampled)
   and
tolerance $\tau$.
   }
   \label{alg:ICS}
  \begin{algorithmic}[1]
    \Procedure{ISMA}
{$A$, $\bv{p}$, $k$, $c$, $w$, $r$, $\tau$, $rows$}
    \State $ S \gets \{1,2, \cdots, n\}$
    \State $\tilde{U}, \tilde{\Sigma}, S \gets \Call{Get Update}{A, S, 0, \bv{p}, c, w, rows}$;
    \If {$rows == 1$} \Comment{Orthonormalise $U$ when there is row sampling}
    \State $Q,R \gets \text{QR}({\tilde{U}}); U,\Sigma,V \gets \text{SVD}(R); \tilde{U} = QU$
    \EndIf
    \State ${\xi_i} \gets 0 \text{ where }i=1,2,\cdots,k$
    \While{$\exists\xi_i : \underset{i=1,2,\cdots,k}{\xi_i}<\tau \text{ and } S$ not empty}
    \State $p_i\gets 1/s$; or $p_i\gets ||\bv{a}_i-\tilde{U}\tilde{U}^T\bv{a}_i||_2^2/\sum||\bv{a}_i-\tilde{U}\tilde{U}^T\bv{a}_i||_F^2\text{, where }i\in S$
    \State $\hat{U}, \hat{\Sigma}, S \gets \Call{Get Update}{A, S, \tilde{U}_k, \bv{p}, c, w, rows}$; 
    \State $\hat{U},\hat{\Sigma}\gets \Call{Block Merge}{\tilde{U},\tilde{\Sigma},\hat{U},\hat{\Sigma},r}$
	  \State $\xi_i \gets {\tilde{\bv{u}}_i}^T\hat{\bv{u}}_i \text{ or } X,\Xi,Y\gets \text{SVD}(\tilde{U}^T_k\hat{U}_k); \xi_i\gets \Xi_{ii}\text{ where }i=1,2,\cdots,k $
    \State $\tilde{U}\gets \hat{U}$; $\tilde{\Sigma}\gets \hat{\Sigma}$
    \EndWhile
	  \State $Q,R\gets\Call{QR}{A^T\tilde{U}}$ \Comment{for computation of error measure}\label{code:extra_iteration_start}
    \State $U_R\Sigma_R V_R^T \gets \Call{SVD}{R}$
    \State $\tilde{U},\tilde{\Sigma} \gets \tilde{U}V_R, \Sigma_R$
    \State $\tilde{V} \gets QU_R$ \Comment{only computed when error measure is needed}\label{code:extra_iteration_end}
    \State \textbf{return } $\tilde{\bv{u}}_i, \tilde{\sigma}_i, \tilde{\bv{v}}_i \text{ where } i=1,2,\cdots,k$
    \EndProcedure
   
    \Procedure{Get Update}{$A$, $S$, $U$, $\bv{p}$, $c$, $w$,$rows$}
    \State $\bv{I}_s\gets\underset{i=1,2,\cdots,c}{s_i}  : \textbf{Pr}[s_i=\alpha]=p_{\alpha} \in \bv{p}$, $s_i \in S$
    \State $\bv{I}_c\gets s_i:\forall \text{ distinct } s_i\in \bv{I}_s $;  $g\gets \text{card}(\bv{I}_c)$; $S \gets S\setminus{\bv{I}_c}$
    \State $\bv{d}_i \gets \bv{a}_{\eta_i},\text{ where }\eta_i\in\bv{I}_c \text{ and } i=1,2,\cdots,g$
    \If {$rows \neq 1$}
    \State $\tilde{V}\tilde{\Sigma}^2\tilde{V}^T \gets \Call{SVD}{D^TD}$; $\tilde{\bv{u}}_i \gets D\tilde{\bv{v}}_i/\tilde{\sigma}_i, i \in 1,2,\dots,r$
    \State \textbf{return } $\tilde{U}$, $\tilde{\Sigma}$, S
    \Else
    \State $q_i\gets ||\bv{d}^i||_2^2$; $q_i\gets q_i/\sum\limits_{i=1}^m q_i$ or if $U\ne 0$ $q_i\gets |\bv{u}^i|_2^2/k$
    \State $W \gets \Call {Sample and Scale Unique Rows}{D,\bv{q},w}$
    \State $\tilde{V}\tilde{\Sigma}^2\tilde{V}^T \gets \Call{SVD}{W^TW}$; $\tilde{\bv{u}}_i \gets D\tilde{\bv{v}}_i/\tilde{\sigma}_i, i \in 1,2,\dots,r$
    \State \textbf{return } $\tilde{U}$, $\tilde{\Sigma}$, S
    \EndIf
    \EndProcedure
    \Procedure {Sample and Scale Unique Rows}{$C,\bv{q},w$}
    \State $\bv{\hat{I}}_s\gets \hat{s}_j \text{ containing $w$ row indices sampled with probability  }
    \bv{q}$
    \State $\bv{I}_r\gets \hat{s}_j; \bv{\hat{t}}\gets \{\hat{t}_j, \text{ number of occurrences of }\hat{s}_j \}:\forall \text{ distinct } \hat{s}_j\in \bv{\hat{I}}_s $
    \State $h\gets \text{card}(\bv{I}_r)$
    \State $\bv{w}^i \gets \bv{c}^{\xi_i}\sqrt{\hat{t}_i/(wq_{\xi_i})},\text{ where }\xi_i\in\bv{I}_r \text{ and }i=1,2,\cdots,h$    
    \State \textbf{return} {$W$}
    \EndProcedure
    
  \end{algorithmic}
\end{algorithm}

  Algorithm~\ref{alg:ICS}  details the steps involved in the iteration.
  The main iteration is the procedure ISMA (Iterative Sampling and Merging Algorithm). 
  It uses three other procedures: GET UPDATE, SAMPLE AND SCALE UNIQUE ROWS and BLOCK MERGE
 The procedure GET UPDATE returns the POD modes of 
the newly sampled columns. In case rows are also sampled (the parameter $rows=1$), GET UPDATE in turn calls the procedure SAMPLE AND SCALE UNIQUE ROWS that performs a row sampling, removes duplicates and scales the remaining rows.
BLOCK MERGE is used to update the dominant subspace. 
It is described further in section~\ref{sec:MAT}.

  The matrix $D$ in GET UPDATE contains the (distinct) sampled columns and the matrix $W$ in SAMPLE AND SCALE UNIQUE ROWS is obtained after additionally sampling rows. If there is no row sampling, GET UPDATE returns the $r$ dominant POD modes, $\tilde{U}_r$ and singular values which are computed 
using a SVD of $D^TD$. Otherwise, the approximate left singular vectors are computed using the right singular vectors of $W$ (line 31).

When using right singular vectors of $W$ as an approximation to those of $D$,  the modes obtained may not be orthonormal. 
Therefore, for the first approximation of the modes, we
orthonormalize using QR and SVD of a smaller matrix (lines 4-6 of  Algorithm~\ref{alg:ICS}). In the subsequent iterations, orthonormalization 
is implicit in the BLOCK MERGE algorithm.

The iteration stops when the POD modes computed in two successive iterations are similar. 
Similarity is  measured using cosines between the modes or principal
cosines between subspaces spanned by the modes. After the iterations converge, the quality of the 
approximated subspace can be measured using Wedin's theorem. It requires estimation of the singular values and 
right singular vectors (lines~\ref{code:extra_iteration_start}-\ref{code:extra_iteration_end} in Algorithm~\ref{alg:ICS}). We describe this measure in section~\ref{sec:meas}

\subsection{Sampling strategies}
\label{sec:samp_strat}
\subsubsection{Column sampling}
In the first round, we do a column norm based sampling, with probabilities as given in equation~(\ref{eqn:LTSVD_prob}),
to get an approximation of the dominant modes.
Additional columns are sampled in subsequent iterations in an attempt to correct the modes. 
The sampling strategy for the additional columns is a compromise between the number of iterations and the computational effort required in each iteration. 
We explored three sampling strategies for adding more columns.
	\begin{enumerate}[label=(\alph*)]
	\item Uniform sampling (UNF) is easy to implement as no extra computation for probabilities is necessary. Larger number of iterations are needed as there is no way to avoid sampling of unimportant or 
highly correlated columns.
\item $L_2$ norm sampling (L2N) has only simple arithmetic manipulation of probabilities from the second iteration but may repeatedly sample columns that are similar or highly correlated.
\item ORT requires us to project the matrix onto the current estimate of the dominant subspace ($\tilde{U}$) and sample based on norm of the component orthogonal to this subspace.
	Probability of sampling columns is computed as,
		\begin{equation}
			p_i=|\bv{a}_i-\tilde{U}\tilde{U}^T\bv{a}_i|_2^2/\sum{|\bv{a}_i-\tilde{U}\tilde{U}^T\bv{a}_i|_2^2}
		\end{equation}
		
Note that, it avoids sampling of highly correlated components and can lead to fewer iterations. 
However, the computational and memory requirements in each iteration is large, making each iteration slow and sometimes impractical for large datasets. 
\end{enumerate}

\subsubsection{Row sampling}
  We have tried two strategies for row sampling - 
	\begin{enumerate}[label=(\alph*)]
	  \item $L_2$ norm of rows and 
	  \item leverage scores (LS).
  \end{enumerate}
  Sampling based on $L_2$ norm of the rows has been used in \cite{DriKanMah:2006:II} as well as in the adaptive sampling algorithm. 
Both leverage score based sampling and uniform sampling has been used by \cite{YamTomDon:2017}, who  propose an incremental algorithm for SVD. 
Their algorithm reduces the size of the matrix by projecting columns and sampling rows. Computing leverage scores requires an estimate of the left singular vectors.
Therefore, it is not done in the first round of sampling.
 
When rows are also sampled, we use the right singular vectors and the singular values of $W$ as an approximation to those  of $D$. These are used to compute the POD modes (see line 31 of Algorithm~\ref{alg:ICS} ). For this approximation to be accurate, the primary requirement is that the singular values computed using $W$ be accurate. In order to maintain the accuracy of the singular values, we scale the sampled rows as specified in the procedure SAMPLE AND SCALE UNIQUE ROWS. 
We remove duplicate rows that were sampled and scale 
the remaining rows with $\sqrt{\hat{t}_i/wq_i}$ where $w$ is the number of rows sampled, $\hat{t}_i$ is the number of occurrences of each sampled row, 
and $q_i$ denotes sampling probability. This scale factor makes the computed SVD equal to the SVD obtained without removing duplicates.
 This is similar to what is done in \cite{SunXieZhaFal:2007} for sampling columns.

\subsection{Merge and truncate (MAT) Operation}
\label{sec:MAT}
Our algorithm requires an efficient technique for estimation of POD vectors in each iteration.  Algorithms to estimate these vectors have been proposed in 
(\cite{QuOstSamGei:2002,BaiChaLuk:2005,LiaBalKanWoo:2014}). 
Also, several techniques have been proposed for incremental SVD \cite{Bra:2002,GuEis:1993,Bra:2006,BakGalDoo:2012}, which can be used to get the principal components. 
We decided to use the algorithm in \cite{VasRam:2017} (which is a generalization of the algorithm in \cite{IweOng:2016}) for the following reasons.
	\begin{enumerate}[label=(\alph*)]
		\item We directly obtain 
the left singular vectors as the output of the operation. 
		\item These left singular vectors are obtained using SVD and QR of smaller matrices. Therefore they are orthogonal by construction.
	\end{enumerate}

Let $A$ be a matrix of size $m\times n$ that consists of sub-matrices $X$ and $Y$ containing $n_1$ and $n_2$ columns respectively, where
$n = n_1+n_2$. Assume that the rank-$r$ approximations of $X$ and $Y$ are $U_{1_r}\Sigma_{1_r}V_{1_r}^T$ and $U_{2_r}\Sigma_{2_r}V_{2_r}^T$. The rank-$r$ approximation of $A$, 
$\tilde{A}_r$,  can be computed from the individual SVDs by first merging the two SVDs and
then truncating it to a rank-$r$ approximation as follows.
 
 The component of $U_{2_r}$ orthogonal to $U_{1_r}$ is  ${U_t =U_{2_r}-U_{1_r}(U_{1_r}^TU_{2_r})}$. If ${U_t = U_oR}$ is the
corresponding QR decomposition, we have
\begin{align}
\begin{bmatrix} X & Y \end{bmatrix} & \approx \begin{bmatrix} U_{1_r} & U_o \end{bmatrix} \begin{bmatrix} \Sigma_{1_r} & (U_{1_r}^T U_{2_r})\Sigma_{2_r} \\ \mathbf{0} & R\Sigma_{2_r} \end{bmatrix} \begin{bmatrix} V_{1_r}^T & \mathbf{0} \\ \mathbf{0} & V_{2_r}^T \end{bmatrix} 
= \begin{bmatrix} U_{1_r} & U_o \end{bmatrix} E \begin{bmatrix} V_{1_r}^T & \mathbf{0} \\ \mathbf{0} & V_{2_r}^T \end{bmatrix}
\label{eqn:orthmerge1}
\end{align}
Now $E$ is a much smaller $(2r) \times (2r)$ matrix.  If $E = U_E \Sigma_E V_E^T$, we get
\begin{align}
\begin{bmatrix} X & Y \end{bmatrix} &\approx \begin{bmatrix} U_{1_r} & U_o \end{bmatrix} U_E \Sigma_E V_E^T \begin{bmatrix} V_{1_r}^T & \mathbf{0} \\ \mathbf{0} & V_{2_r}^T \end{bmatrix} = U \Sigma V^T
\label{eqn:orthmerge}
\end{align}
where $U = \begin{bmatrix} U_{1_r} & U_o \end{bmatrix} U_E$, $\Sigma = \Sigma_E$ and $V =  \begin{bmatrix} V_{1_r} & \mathbf{0} \\ \mathbf{0} & V_{2_r} \end{bmatrix}V_E $. These singular values and vectors are once again truncated to get $\tilde{U}_r$ and $\tilde{\Sigma}_r$. 
   Algorithm~\ref{alg:block_merge} details the steps involved.
   In general, if there are $P$ partitions, a series of MAT operations can be used to obtain the $r$ left singular vectors. Since we are only interested in the left singular vectors, we do not compute $V$.

 The penalty we pay for getting a more accurate estimate of the dominant subspace is finding the SVD of the newly sampled columns to do a MAT operation in each iteration.
However, the number of additional sampled columns in each iteration is much lower 
than the total number of columns in the matrix. 
Therefore, a speedup is obtained since the 
complexity of SVD computation depends quadratically on the number of columns for 
``tall and thin'' matrices.\footnote{
For ``short and fat'' matrices, we need to iteratively sample
rows to get a runtime improvement when run on a single machine.}
Assuming there are $P$ partitions and all $n$ columns are sampled with same number of columns sampled 
in each iteration, the total number of flops 
required for computing a rank-$r$ approximation of $A$ using a series of MAT operations is  $ \approx \frac{14mn^2}{P} + \frac{192 n^3}{P^2}$ where $ r << n$ \cite{VasRam:2017}. 

\begin{algorithm}[htbp!]
	\caption{Pseudocode for block merge algorithm: merges left/right singular vectors of two adjacent blocks.
  Inputs: $U_1 \in \mathbb{R}^{m_1\times n_1}$, $U_2 \in \mathbb{R}^{m_1\times n_2}$, and required rank $r$ 
  }
  \label{alg:block_merge}
  \small
	\begin{algorithmic}[1]
          \Procedure{Block Merge}{$U_1,\Sigma_1,U_2,\Sigma_2,r$}
          \State $U_{1_r},\Sigma_{1_r},U_{2_r},\Sigma_{2_r}$ $\gets$ Truncate( $U_1,\Sigma_1,U_2,\Sigma_2$) retaining only $r$ values
	  \State $U_t\gets U_{2_r}-U_{1_r}(U_{1_r}^TU_{2_r});U_0,R\gets \Call{QR}{U_t}$
	  \State $E\gets \begin{bmatrix} \Sigma_{1_r} & (U_{1_r}^T U_{2_r})\Sigma_{2_r} \\ \mathbf{0} & R\Sigma_{2_r} \end{bmatrix}; U_E\Sigma_EV_E^T\gets \Call{SVD}{E}$
          \State $U\gets \begin{bmatrix}U_{1_r}& U_0\end{bmatrix}U_E; \Sigma\gets \Sigma_E ;r \gets \min(r,\max(i:\sigma_i \neq 0))$
          \State Truncate $U$ and $\Sigma$ retaining only $r$ values
			\State \textbf{return} $U,\Sigma$
		\EndProcedure
	\end{algorithmic}
\end{algorithm}

\subsubsection{Error Analysis for series of MAT operations}
\label{subsec:error_MAT}
We derive an error bound in terms of the spectral norm of the error in the rank-r matrix obtained using a series of MAT operations.
For this analysis, we use the following well known results from literature \cite{Li:2014,Mathias:2014}.
Let $X$ be an $m\times n$ matrix and $q=\min(m,n)$.
\begin{enumerate}
\item If $X = U\Sigma V^T$, and $X_r=U_r\Sigma_rV_r^T$ is a rank-$r$ approximation to $X$ obtained by retaining the top $r$ singular values,
  then, $||X-X_r||_2 = \sigma_{r+1}(X)$ and $||X-X_r||_2\leq||X-B_r||_2$, where $B_r$ is any rank-$r$ matrix.
\item If $Y$ is an $m\times t$ sub-matrix of $X$, then, $\sigma_{i+n-t}(X)\leq\sigma_i(Y)\leq\sigma_i(X)$
\item If $\hat{X}=X+\Delta X$ is a perturbation of $X$, then $|\sigma_i-\hat{\sigma}_i|\leq||\Delta X||_2$, $i = 1, 2, \dots, q$.
\end{enumerate}
We use the following notation. Let the matrix $X = \begin{bmatrix} X_1 & X_2& \cdots &X_P \end{bmatrix}$ be divided into $P$ partitions.
Define $Y_0 = \begin{bmatrix} X_1  \end{bmatrix}$, $Y_1 =  \begin{bmatrix} Y_0 & X_2 \end{bmatrix}$ and so on.
In general $Y_j =  \begin{bmatrix} Y_{j-1} & X_{j+1} \end{bmatrix}$, so that $Y_{P-1} = X$.
Let $X_r$ denote the optimal rank-r approximation of $X$, and ${X_i}_r$ that of $X_i$.
Define $Z_1 = \begin{bmatrix} {X_1}_r & {X_2}_r \end{bmatrix}$ and $Z_j =  \begin{bmatrix} {Z_{j-1}}_r & {X_{j+1}}_r \end{bmatrix}$.
Let $Z = Z_{P-1}$. Therefore the error $E = ||X - Z_r||_2$, where $Z_r$ are the optimal rank-$r$ approximations of $Z$. 

The bound for the error can be derived as follows. 
\begin{align}
 ||X - Z_r||_2 & \leq ||X - Z||_2 + ||Z  - Z_r||_2  \nonumber \\
                       & \leq ||X - Z||_2 + ||Z  - X_r||_2  \nonumber \\
                       & \leq ||X - Z||_2 + ||Z  - X||_2  + ||X_r - X||_2  \nonumber \\
                       & =  ||X_r - X||_2  + 2||X - Z||_2 \nonumber \\
                       & =  \sigma_{r+1}(X) + 2||X - Z||_2
\end{align}

We find $||X - Z||_2$ iteratively. For the MAT operation on the first two partitions, we have
\begin{align}
  ||Z_1 - Y_1||_2 & \leq  ||X_1 - {X_1}_r||_2 + ||X_2 - {X_2}_r||_2 \nonumber \\
  & = \sigma_{r+1}(X_1) + \sigma_{r+1}(X_2) \leq 2\sigma_{r+1}(X)
\end{align}
For subsequent partitions,
\begin{align}
  ||Z_j - Y_j||_2 & \leq  ||{Z_{j-1}}_r - Y_{j-1}||_2 + ||X_{j+1} - {X_{j+1}}_r||_2 \nonumber \\
  & \leq ||{Z_{j-1}}_r - Z_{j-1}||_2 + ||Z_{j-1} - Y_{j-1}||_2 + ||X_{j+1} - {X_{j+1}}_r||_2 
\end{align}

If $Z_i = Y_i + \Delta X$, then $\sigma_{r+1}(Z_i) \leq \sigma_{r+1}(Y_i) + ||Z_i - Y_i||_2$.  Therefore,
\begin{align}
  ||Z_j - Y_j||_2 &   \leq \sigma_{r+1}(Y_{j-1}) + 2 ||Z_{j-1} - Y_{j-1}||_2 + \sigma_{r+1}(X_{j+1}) \nonumber \\
  & \leq \sigma_{r+1}(X) + 2 ||Z_{j-1} - Y_{j-1}||_2 + \sigma_{r+1}(X)\label{eqn:err_MAT}
\end{align}
The last inequality follows since  $Y_{j-1}$  and $X_{j+1}$ are submatrices of $X$. 
 Since $||X - Z||_2 = ||Y_{P-1} - Z_{P-1}||_2$,  it can be computed iteratively using equation~(\ref{eqn:err_MAT}).
Therefore, the bound for error, $||X-Z_r||_2$, is
  $(2^{P+1}-3)\sigma_{r+1}(X)$. \footnote{ The trend for the error is similar in the expression derived in \cite{IweOng:2016} (see Theorem 3 in the reference). They find a bound on the Frobenius norm of the error  which grows as $(1+\sqrt{2})^P$, which is slightly worse than ours. We used the spectral norm as it is more appropriate for noisy data as argued in \cite{LiLinSzlStaKluTyg:2017}. Our derivation is also simpler and more specific to the case we are looking at, namely, an iterative algorithm which has a series of MAT operations. }
We note the following.
\begin{enumerate}
\item In the worst case, $\sigma_{r+1}(X_i)$ is close to $ \sigma_{r+1}(X)$ and the error is close to $({2^{P+1}-3}) \sigma_{r+1}(X)$.
But this
  typically happens when the size of the submatrix is large, which means $P$ is small.
\item If there are many partitions, the number of MAT operations are large, but each  $\sigma_{r+1}(X_i)$ is likely to be small.
If $Y$ is an  $m \times t$ submatrix of $X$, then $\sigma_{r+1+n-t} (X) \leq \sigma_{r+1}(Y) \leq \sigma_{r+1}(X)$.
Since each submatrix is thin and $t$ is small, the lower bound is quite small.
Hence the error in the average case is likely to be much smaller than in the worst case. In practice, we have found this is true.
\end{enumerate}
Overall, if we need $k$ POD vectors to be approximated well, then either the number of partitions should be small or we can set $r$ to be larger than $k$. We set $r$ to $3k$ to get good accuracies. 

\subsection{Quality of converged solution}
\label{sec:meas}
In order to estimate the accuracy of the computed POD modes, we need a measure to quantify the error in the approximation.
After convergence, we can get a bound on the error in the subspaces using 
Wedin's theorem \cite{Wedin:1972}.
Assume that at the end of the iteration, we get $\tilde{A}_r = \tilde{U}_r\tilde{\Sigma}_r \tilde{V}_r^T$, out of which we take the top $k$ singular values and vectors. Define
$R = A\tilde{V}_k - \tilde{U}_k \tilde{\Sigma}_k$ and $S = A^T\tilde{U}_k -  \tilde{V}_k\tilde{\Sigma}_k$. 
Assume $\sin(\Theta)$ and $\sin(\Phi)$ are diagonal matrices containing the sine of principal angles between subspaces spanned by $U_k$ and $\tilde{U}_k$ and
$V_k$ and $\tilde{V}_k$ respectively. The sines of the principal angles are a measure of the distance between the accurate and approximate dominant left and right subspaces.

Wedin's theorem \cite{Wedin:1972,Li:1996,Ste:1990} gives a bound on the norms of the sine of principal angles in terms of the norms of $R$ and $S$. The statement of the theorem is as 
follows.
\begin{equation}
  \label{eqn:ang_meas}
  \sqrt{||\sin(\Theta)||_F^2 + ||\sin(\Phi)||_F^2 } \leq \frac{\sqrt{||R||_F^2 + ||S||_F^2}}{\omega}
\end{equation}
where $\omega \triangleq \min\{ \min\limits_{j=1, \cdots n-k} |\tilde{\sigma}_k - \sigma_{k+j}|, \tilde{\sigma}_k \} > 0$ and $\sigma_i$ are the singular values of $A$. 
In order to compute $\omega$, we need all the singular values of $A$, which are not available. However, at the end of the 
iteration, we have an approximation of $r(=3k)$ singular values. Assuming that the $k^{th}$ and $(k+1)^{th}$
singular values are approximated well (which was the case in all the datasets we have looked at),  we estimate $\omega$ as $ \min\{|\tilde{\sigma}_k - \tilde{\sigma}_{k+1}|, \tilde{\sigma}_k\}$. If our iteration converges to the correct subspaces, the right hand side should be close to zero. If it approaches
$\sqrt{2k}$, either (a) additional columns need to be sampled or (b) the $k^{th}$ and $(k+1)^{th}$ singular values are very close to each other and small
errors in the singular values result in large values for the bound. In the limit when they are identical, any vector in the subspace is a POD mode and the 
bound can become arbitrary. In this case, a better measure of the quality of the solutions can be obtained after increasing $k$. Usually increasing $k$ by one or two 
is sufficient. 

Our main focus is the POD modes. However, in order to compute the bound, the singular values and right singular vectors are also required. Since we do not scale columns and use some subset of the columns for computation of the modes, the error in the singular values (not the vectors) could be significant. This is especially true  if the total number of columns
sampled is much less than the number of columns in the matrix. The left singular vectors will be quite accurate if the columns corresponding to the dominant subspace are captured, but the singular values will not be accurate.
If this is the case, we use the technique followed in the
adaptive sampling method i.e. after the modes have converged, the singular values and the right singular vectors can be obtained by computing the SVD of $\tilde{U}_k^TA$
as detailed in lines~\ref{code:extra_iteration_start}-\ref{code:extra_iteration_end} in Algorithm~\ref{alg:ICS}. 

\section{Results}\label{sec:results}

All algorithms are run in a 4 core Intel\textsuperscript{\textregistered} Core\textsuperscript{TM} i7-6700K processor that runs at a maximum clock frequency of 4GHz with hyper-threading on. The system has a 32GB RAM.  For the Yale Faces dataset (YF), the mean centered data required 40 GB memory. Therefore for this dataset, the runs for the iterative algorithm were done  in a different  system, a 8 core Intel\textsuperscript{\textregistered} Xeon\textsuperscript{\textregistered} CPU E5-2650V2 processor that runs at a maximum clock frequency of 2.6GHz with hyper-threading on, with 64GB of RAM.

The algorithms are written in Python and run in Python version 3.5.3. It uses \emph{numpy} version 1.12.1-3, \emph{scipy} version 0.18.1-2 and \emph{openBLAS} version 0.2.19-3, which is multi-threaded. Double precision was used for all computations. We used the linear algebra routines in Scipy, which are based on LAPACK for SVD computation.
The random number generator in \emph{numpy} is used. Runtime and speedup values are an average of five runs.
We report speed-up obtained when \emph{openBLAS} is  run using a single thread, to give a measure of the operation count in each algorithm 
 and four threads (since our system is a four core system), which gives an indication of how the runtime scales with number of threads.
 We have not explicitly parallelised our algorithm and any improvement with multiple threads is entirely due to \emph{openBLAS}.
 
The accurate POD modes of $A$ can be computed using SVD($A^TA$) as follows, (a) find right singular vectors and singular values of $A^TA$ 
using either SVD or eigendecomposition, (b) compute the $k$ left singular values as $\bv{u}_i=A\bv{v}_i/\sigma_i$. \textit{We denote the algorithm as {SVD-$A^TA$}.}

If $A$ is a tall and skinny matrix, computing SVD($A$) using {SVD-$A^TA$} is faster.
 This is because SVD($A$) takes approximately $\beta mn^2$ floating point operations (FLOPs). SVD-$A^TA$ requires $mn^2$ FLOPs to compute $A^TA$, 
 $(\beta + 16)n^3$ to compute SVD($A^TA$) and another $mn^2$ operations to compute $\bv{u}_i$ 
 from $A\bv{v}_i/\sigma_i$. Assuming $\beta=6$  \cite{GolVan:2007,Bjo:2015}, this gives $m/(m/3 + (1+8/3)n)$ speedup. The results in Table \ref{tab:time_SVD_ATA} are close to what is expected.
Since SVD-$A^TA$ is faster, we benchmark the speedup of the algorithms with respect to time taken for computing $k$ POD modes using SVD-$A^TA$.
 \begin{table}[tph!]
   \centering
   \small
   \caption{Average time taken (in seconds) by SVD($A$) and SVD-$A^TA$ for different datasets. Note that SVD of Yale faces could not be computed in our test system (32GB) since it does not fit in the system's memory.}
 \label{tab:time_SVD_ATA}
   \begin{tabular}{ccc|cc}
     \toprule
     Dataset& \multicolumn{2}{c|}{SVD($A$)} & \multicolumn{2}{c}{SVD-$A^TA$}\\
     &single thread& 4 threads&single thread& 4 threads\\
     \midrule
     $V_{2D}$&$25.28$&$13.94$&$8.39$&$2.75$\\
     Faces&$0.35$&$0.17$&$0.13$&$0.069$\\
     Cropped Yale faces (CYF) &$34.04$&$16.43$&$15.23$&$6.86$\\
     Yale faces (YF) &-&-&-&-\\
     \bottomrule
   \end{tabular}
 \end{table}

\subsection{Accuracy and runtime of the iterative algorithms}
We evaluated the two versions of ISMA algorithm (Algorithm~\ref{alg:ICS}) - an iterative column sampling algorithm (ICS) and an iterative column and row 
sampling algorithm (ICRS).
We tested the performance of the iterative algorithms with two different convergence criteria: cosines between the modes obtained in successive iterations and the principal cosines between the subspaces computed in successive iterations. As indicated previously, we tried different sampling strategies from the second iteration.
Table~\ref{tab:sampl_strategy} lists the different sampling strategies we use in this paper.
\begin{table}[htbp!]
\centering
\caption{Sampling strategies used for sampling rows and columns by iterative sampling algorithms. For all sampling strategies, in the first 
iteration, columns (and rows in the case of column and row sampling) are sampled based on their respective $L_2$ norms.}
\label{tab:sampl_strategy}
\begin{tabular}{cp{5cm}p{5.5cm}}
\toprule
strategy&column sampling&row sampling (when applicable)\\
\midrule
L2N&$L_2$ norm sampling of columns in all iterations&\begin{center}-\end{center}\\
UNF&uniform sampling from second iteration&\begin{center} $L_2$ norm sampling\end{center}\\
ORT&from second iteration, sampling based on orthogonal complement of $A$ projected onto current approximation of modes&\begin{center} $L_2$ norm sampling\end{center}\\
LS&uniform sampling from second iteration&leverage score sampling using current approximation of  modes from second iteration\\
\bottomrule
\end{tabular}
\end{table}

Using $\epsilon$ and $\delta$ around $0.6-0.7$ resulted in a reasonable trade-off between the number of columns (and rows) sampled in each iteration and
the number of iterations. For the Faces dataset, we chose a higher $\epsilon, \delta$ since the number of columns in the matrix is quite small.
ORT was not used in the iterative algorithms for YF, since the overhead for computing the orthogonal complement in each iteration proved to be large.
Table~\ref{tab:param_iter} lists the parameter values used for the iterative algorithms. In all cases, the tolerance parameter for convergence $\tau$, is set to 0.99.

\begin{table}[tbhp!]
  \scriptsize
  \caption{Parameter values ($k$, $\epsilon$, $\delta$) used in the iterative algorithms for various datasets}
  \label{tab:param_iter}
  \centering
  \small
  \begin{tabular}{p{1cm}c|c|c}
    \toprule
    & $V_{2D}$, CYF&Faces&YF\\
    Case&$k$, $\epsilon$, $\delta$&$k$, $\epsilon$, $\delta$&$k$, $\epsilon$, $\delta$\\
    \midrule
    I1&10, 0.7, 0.6&10, 1, 0.75&20, 0.6, 0.6\\
    I2&2, 0.7, 0.6&2, 0.7, 0.6&5, 0.6, 0.6\\
    \bottomrule
  \end{tabular}
\end{table}

\begin{figure}[!hp]
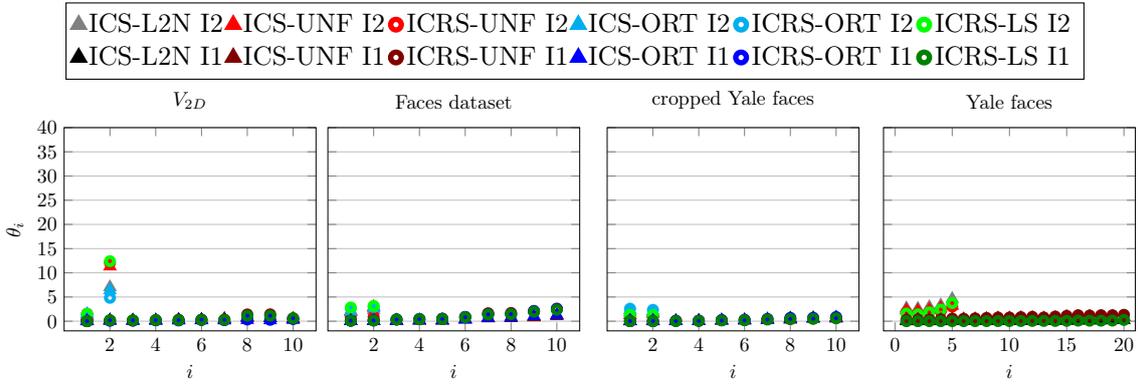

  \begin{minipage}{\textwidth}
    \centering
  \ref{angPOD_iter}
  \end{minipage}
  \begin{minipage}{\textwidth}
   \hspace{-0.37in} \scalebox{0.75}{\input{figures/tikz_plot_angPOD_V_2D_iterate_svdATA.tex}}
   \hspace{-0.31in}\scalebox{0.75}{\input{figures/tikz_plot_angPOD_faces_iterate_svdATA.tex}}
   \hspace{-0.36in} \scalebox{0.75}{ \input{figures/tikz_plot_angPOD_croppedYalefaces_iterate_svdATA.tex}}
   \hspace{-0.28in}\scalebox{0.75}{ \pgfplotstableread{data/angPOD_LTSVD_no_scaling_iterate_Yalefaces.dat}\1
\pgfplotstableread{data/angPOD_LTSVD_US_no_scaling_iterate_Yalefaces.dat}\4
\pgfplotstableread{data/angPOD_CTSVD_US_no_scaling_iterate_Yalefaces.dat}\5
\pgfplotstableread{data/angPOD_CTSVD_LS_no_scaling_iterate_Yalefaces.dat}\7
\begin{tikzpicture}[every mark/.append style={line width=2pt, solid}]
\begin{axis}[
    title = { Yale faces}, 
    xlabel = $i$,
    xmax = 21.000000,
    ytick={0,5,...,40},
    yticklabels={,,,},
    width = 6cm,
    ymajorgrids,
     ymax = 40,
    ymin =-2,
  ]
  \addplot[
		 color=black!50!white,
		 only marks,
		 mark=triangle,
     mark size = 2pt,
		 x filter/.code={\pgfplotstablegetelem{\coordindex}{[index]0}\of{\1}
\pgfmathtruncatemacro{\temp}{abs(\pgfplotsretval-5)==0? 1 : 0}
\ifnum\temp>0
\pgfplotstablegetelem{\coordindex}{[index]1}\of{\1}
\pgfmathtruncatemacro{\temp}{abs(\pgfplotsretval-0.6)==0? 1 : 0}
\ifnum\temp>0
\pgfplotstablegetelem{\coordindex}{[index]2}\of{\1}
\pgfmathtruncatemacro{\temp}{abs(\pgfplotsretval-0.6)==0? 1 : 0}
\ifnum\temp>0
\relax
			\else
  \def\pgfmathresult{}
\fi
\else
  \def\pgfmathresult{}
\fi
\else
  \def\pgfmathresult{}
\fi
},
		]
                table[ 
                x index = 3,
                y index = 4,
                ]
		{\1};
  \addplot[
		 color=red,
		 only marks,
		 mark=triangle,
     mark size = 2pt,
		 x filter/.code={\pgfplotstablegetelem{\coordindex}{[index]0}\of{\4}
\pgfmathtruncatemacro{\temp}{abs(\pgfplotsretval-5)==0? 1 : 0}
\ifnum\temp>0
\pgfplotstablegetelem{\coordindex}{[index]1}\of{\4}
\pgfmathtruncatemacro{\temp}{abs(\pgfplotsretval-0.6)==0? 1 : 0}
\ifnum\temp>0
\pgfplotstablegetelem{\coordindex}{[index]2}\of{\4}
\pgfmathtruncatemacro{\temp}{abs(\pgfplotsretval-0.6)==0? 1 : 0}
\ifnum\temp>0
\relax
			\else
  \def\pgfmathresult{}
\fi
\else
  \def\pgfmathresult{}
\fi
\else
  \def\pgfmathresult{}
\fi
},
		]
                table[ 
                x index = 3,
                y index = 4,
                ]
		{\4};
  \addplot[
		 color=red,
		 only marks,
     mark=o,
     mark size = 2pt,
		 x filter/.code={\pgfplotstablegetelem{\coordindex}{[index]0}\of{\5}
\pgfmathtruncatemacro{\temp}{abs(\pgfplotsretval-5)==0? 1 : 0}
\ifnum\temp>0
\pgfplotstablegetelem{\coordindex}{[index]1}\of{\5}
\pgfmathtruncatemacro{\temp}{abs(\pgfplotsretval-0.6)==0? 1 : 0}
\ifnum\temp>0
\pgfplotstablegetelem{\coordindex}{[index]2}\of{\5}
\pgfmathtruncatemacro{\temp}{abs(\pgfplotsretval-0.6)==0? 1 : 0}
\ifnum\temp>0
\relax
			\else
  \def\pgfmathresult{}
\fi
\else
  \def\pgfmathresult{}
\fi
\else
  \def\pgfmathresult{}
\fi
},
		]
                table[ 
                x index = 3,
                y index = 4,
                ]
		{\5};
		\addplot[
		 color=green,
		 only marks,
     mark=o,
     mark size = 2pt,
		 x filter/.code={\pgfplotstablegetelem{\coordindex}{[index]0}\of{\7}
\pgfmathtruncatemacro{\temp}{abs(\pgfplotsretval-5)==0? 1 : 0}7
\ifnum\temp>0
\pgfplotstablegetelem{\coordindex}{[index]1}\of{\7}
\pgfmathtruncatemacro{\temp}{abs(\pgfplotsretval-0.6)==0? 1 : 0}
\ifnum\temp>0
\pgfplotstablegetelem{\coordindex}{[index]2}\of{\7}
\pgfmathtruncatemacro{\temp}{abs(\pgfplotsretval-0.6)==0? 1 : 0}
\ifnum\temp>0
\relax
			\else
  \def\pgfmathresult{}
\fi
\else
  \def\pgfmathresult{}
\fi
\else
  \def\pgfmathresult{}
\fi
},
		]
                table[ 
                x index = 3,
                y index = 4,
                ]
		{\7};
  \addplot[
		 color=black,
		 only marks,
		 mark=triangle,
     mark size = 2pt,
		 x filter/.code={\pgfplotstablegetelem{\coordindex}{[index]0}\of{\1}
\pgfmathtruncatemacro{\temp}{abs(\pgfplotsretval-20)==0? 1 : 0}
\ifnum\temp>0
\pgfplotstablegetelem{\coordindex}{[index]1}\of{\1}
\pgfmathtruncatemacro{\temp}{abs(\pgfplotsretval-0.6)==0? 1 : 0}
\ifnum\temp>0
\pgfplotstablegetelem{\coordindex}{[index]2}\of{\1}
\pgfmathtruncatemacro{\temp}{abs(\pgfplotsretval-0.6)==0? 1 : 0}
\ifnum\temp>0
\relax
			\else
  \def\pgfmathresult{}
\fi
\else
  \def\pgfmathresult{}
\fi
\else
  \def\pgfmathresult{}
\fi
},
		]
                table[ 
                x index = 3,
                y index = 4,
                ]
		{\1};
  \addplot[
		 color=red!50!black,
		 only marks,
		 mark=triangle,
     mark size = 2pt,
		 x filter/.code={\pgfplotstablegetelem{\coordindex}{[index]0}\of{\4}
\pgfmathtruncatemacro{\temp}{abs(\pgfplotsretval-20)==0? 1 : 0}
\ifnum\temp>0
\pgfplotstablegetelem{\coordindex}{[index]1}\of{\4}
\pgfmathtruncatemacro{\temp}{abs(\pgfplotsretval-0.6)==0? 1 : 0}
\ifnum\temp>0
\pgfplotstablegetelem{\coordindex}{[index]2}\of{\4}
\pgfmathtruncatemacro{\temp}{abs(\pgfplotsretval-0.6)==0? 1 : 0}
\ifnum\temp>0
\relax
			\else
  \def\pgfmathresult{}
\fi
\else
  \def\pgfmathresult{}
\fi
\else
  \def\pgfmathresult{}
\fi
},
		]
                table[ 
                x index = 3,
                y index = 4,
                ]
		{\4};
  \addplot[
		 color=red!50!black,
		 only marks,
     mark=o,
     mark size = 2pt,
		 x filter/.code={\pgfplotstablegetelem{\coordindex}{[index]0}\of{\5}
\pgfmathtruncatemacro{\temp}{abs(\pgfplotsretval-20)==0? 1 : 0}
\ifnum\temp>0
\pgfplotstablegetelem{\coordindex}{[index]1}\of{\5}
\pgfmathtruncatemacro{\temp}{abs(\pgfplotsretval-0.6)==0? 1 : 0}
\ifnum\temp>0
\pgfplotstablegetelem{\coordindex}{[index]2}\of{\5}
\pgfmathtruncatemacro{\temp}{abs(\pgfplotsretval-0.6)==0? 1 : 0}
\ifnum\temp>0
\relax
			\else
  \def\pgfmathresult{}
\fi
\else
  \def\pgfmathresult{}
\fi
\else
  \def\pgfmathresult{}
\fi
},
		]
                table[ 
                x index = 3,
                y index = 4,
                ]
		{\5};
		\addplot[
		 color=green!50!black,
		 only marks,
     mark=o,
     mark size = 2pt,
		 x filter/.code={\pgfplotstablegetelem{\coordindex}{[index]0}\of{\7}
\pgfmathtruncatemacro{\temp}{abs(\pgfplotsretval-20)==0? 1 : 0}
\ifnum\temp>0
\pgfplotstablegetelem{\coordindex}{[index]1}\of{\7}
\pgfmathtruncatemacro{\temp}{abs(\pgfplotsretval-0.6)==0? 1 : 0}
\ifnum\temp>0
\pgfplotstablegetelem{\coordindex}{[index]2}\of{\7}
\pgfmathtruncatemacro{\temp}{abs(\pgfplotsretval-0.6)==0? 1 : 0}
\ifnum\temp>0
\relax
			\else
  \def\pgfmathresult{}
\fi
\else
  \def\pgfmathresult{}
\fi
\else
  \def\pgfmathresult{}
\fi
},
		]
                table[ 
                x index = 3,
                y index = 4,
                ]
		{\7};
\end{axis}
\end{tikzpicture}}
\end{minipage}
  \caption{First $k$ mode angles, $\theta_i$, of different datasets computed by iterative algorithms based on convergence of individual modes.
 ICS:Column sampling, ICRS:column and row sampling
} 
  \label{fig:angPOD_iter}
\end{figure}

Fig.~\ref{fig:angPOD_iter} shows the first $k$ mode angles computed by the iterative algorithms using  L2N, UNF, ORT, and LS.
This is run with convergence of individual modes.
As expected, the accuracies are significantly better than for a single round of sampling.
\begin{figure}[!hp]
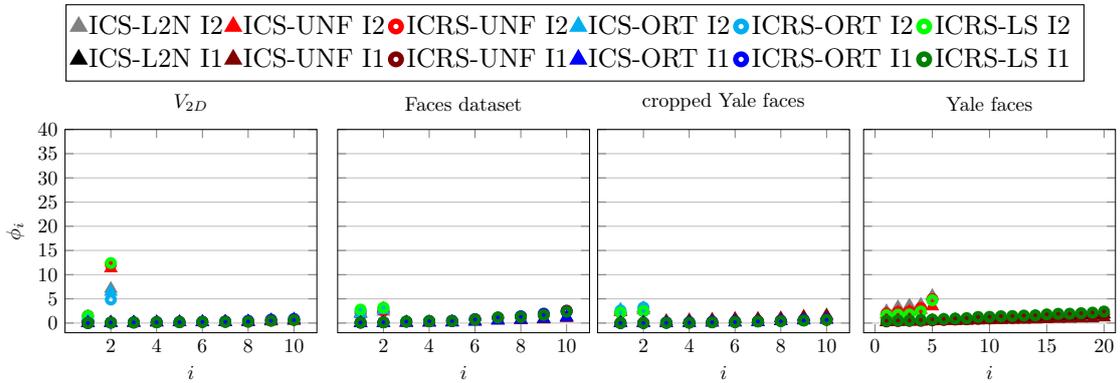

  \begin{minipage}{\textwidth}
    \centering
  \ref{pc_iter}
  \end{minipage}
  \begin{minipage}{\textwidth}
 \hspace{-0.4in}   \scalebox{0.75}{\input{figures/tikz_plot_pc_V_2D_iterate_svdATA.tex}}
 \hspace{-0.37in}    \scalebox{0.75}{\input{figures/tikz_plot_pc_faces_iterate_svdATA.tex}}
 \hspace{-0.4in}    \scalebox{0.75}{ \input{figures/tikz_plot_pc_croppedYalefaces_iterate_svdATA.tex}}
 \hspace{-0.35in}    \scalebox{0.75}{ \pgfplotstableread{data/pc_LTSVD_no_scaling_pc_iterate_Yalefaces.dat}\1
\pgfplotstableread{data/pc_LTSVD_US_no_scaling_pc_iterate_Yalefaces.dat}\4
\pgfplotstableread{data/pc_CTSVD_US_no_scaling_pc_iterate_Yalefaces.dat}\5
\pgfplotstableread{data/pc_CTSVD_LS_no_scaling_pc_iterate_Yalefaces.dat}\7
\begin{tikzpicture}[every mark/.append style={line width=2pt, solid}]
\begin{axis}[
    title = { Yale faces}, 
    xlabel = $i$,
    xmax = 21.000000,
    ytick = {0,5,...,40},
    yticklabels={,,},
    width = 6cm,
    ymajorgrids,
     ymax = 40,
    ymin = -2,
  ]
  \addplot[
		 color=black!50!white,
		 only marks,
		 mark=triangle,
     mark size = 2pt,
		 x filter/.code={\pgfplotstablegetelem{\coordindex}{[index]0}\of{\1}
\pgfmathtruncatemacro{\temp}{abs(\pgfplotsretval-5)==0? 1 : 0}
\ifnum\temp>0
\pgfplotstablegetelem{\coordindex}{[index]1}\of{\1}
\pgfmathtruncatemacro{\temp}{abs(\pgfplotsretval-0.6)==0? 1 : 0}
\ifnum\temp>0
\pgfplotstablegetelem{\coordindex}{[index]2}\of{\1}
\pgfmathtruncatemacro{\temp}{abs(\pgfplotsretval-0.6)==0? 1 : 0}
\ifnum\temp>0
\relax
			\else
  \def\pgfmathresult{}
\fi
\else
  \def\pgfmathresult{}
\fi
\else
  \def\pgfmathresult{}
\fi
},
		]
                table[ 
                x index = 3,
                y index = 4,
                ]
		{\1};
  \addplot[
		 color=red,
		 only marks,
		 mark=triangle,
     mark size = 2pt,
		 x filter/.code={\pgfplotstablegetelem{\coordindex}{[index]0}\of{\4}
\pgfmathtruncatemacro{\temp}{abs(\pgfplotsretval-5)==0? 1 : 0}
\ifnum\temp>0
\pgfplotstablegetelem{\coordindex}{[index]1}\of{\4}
\pgfmathtruncatemacro{\temp}{abs(\pgfplotsretval-0.6)==0? 1 : 0}
\ifnum\temp>0
\pgfplotstablegetelem{\coordindex}{[index]2}\of{\4}
\pgfmathtruncatemacro{\temp}{abs(\pgfplotsretval-0.6)==0? 1 : 0}
\ifnum\temp>0
\relax
			\else
  \def\pgfmathresult{}
\fi
\else
  \def\pgfmathresult{}
\fi
\else
  \def\pgfmathresult{}
\fi
},
		]
                table[ 
                x index = 3,
                y index = 4,
                ]
		{\4};
  \addplot[
		 color=red,
		 only marks,
     mark=o,
     mark size = 2pt,
		 x filter/.code={\pgfplotstablegetelem{\coordindex}{[index]0}\of{\5}
\pgfmathtruncatemacro{\temp}{abs(\pgfplotsretval-5)==0? 1 : 0}
\ifnum\temp>0
\pgfplotstablegetelem{\coordindex}{[index]1}\of{\5}
\pgfmathtruncatemacro{\temp}{abs(\pgfplotsretval-0.6)==0? 1 : 0}
\ifnum\temp>0
\pgfplotstablegetelem{\coordindex}{[index]2}\of{\5}
\pgfmathtruncatemacro{\temp}{abs(\pgfplotsretval-0.6)==0? 1 : 0}
\ifnum\temp>0
\relax
			\else
  \def\pgfmathresult{}
\fi
\else
  \def\pgfmathresult{}
\fi
\else
  \def\pgfmathresult{}
\fi
},
		]
                table[ 
                x index = 3,
                y index = 4,
                ]
		{\5};
		\addplot[
		 color=green,
		 only marks,
     mark=o,
     mark size = 2pt,
		 x filter/.code={\pgfplotstablegetelem{\coordindex}{[index]0}\of{\7}
\pgfmathtruncatemacro{\temp}{abs(\pgfplotsretval-5)==0? 1 : 0}7
\ifnum\temp>0
\pgfplotstablegetelem{\coordindex}{[index]1}\of{\7}
\pgfmathtruncatemacro{\temp}{abs(\pgfplotsretval-0.6)==0? 1 : 0}
\ifnum\temp>0
\pgfplotstablegetelem{\coordindex}{[index]2}\of{\7}
\pgfmathtruncatemacro{\temp}{abs(\pgfplotsretval-0.6)==0? 1 : 0}
\ifnum\temp>0
\relax
			\else
  \def\pgfmathresult{}
\fi
\else
  \def\pgfmathresult{}
\fi
\else
  \def\pgfmathresult{}
\fi
},
		]
                table[ 
                x index = 3,
                y index = 4,
                ]
		{\7};
  \addplot[
		 color=black,
		 only marks,
		 mark=triangle,
     mark size = 2pt,
		 x filter/.code={\pgfplotstablegetelem{\coordindex}{[index]0}\of{\1}
\pgfmathtruncatemacro{\temp}{abs(\pgfplotsretval-20)==0? 1 : 0}
\ifnum\temp>0
\pgfplotstablegetelem{\coordindex}{[index]1}\of{\1}
\pgfmathtruncatemacro{\temp}{abs(\pgfplotsretval-0.6)==0? 1 : 0}
\ifnum\temp>0
\pgfplotstablegetelem{\coordindex}{[index]2}\of{\1}
\pgfmathtruncatemacro{\temp}{abs(\pgfplotsretval-0.6)==0? 1 : 0}
\ifnum\temp>0
\relax
			\else
  \def\pgfmathresult{}
\fi
\else
  \def\pgfmathresult{}
\fi
\else
  \def\pgfmathresult{}
\fi
},
		]
                table[ 
                x index = 3,
                y index = 4,
                ]
		{\1};
  \addplot[
		 color=red!50!black,
		 only marks,
		 mark=triangle,
     mark size = 2pt,
		 x filter/.code={\pgfplotstablegetelem{\coordindex}{[index]0}\of{\4}
\pgfmathtruncatemacro{\temp}{abs(\pgfplotsretval-20)==0? 1 : 0}
\ifnum\temp>0
\pgfplotstablegetelem{\coordindex}{[index]1}\of{\4}
\pgfmathtruncatemacro{\temp}{abs(\pgfplotsretval-0.6)==0? 1 : 0}
\ifnum\temp>0
\pgfplotstablegetelem{\coordindex}{[index]2}\of{\4}
\pgfmathtruncatemacro{\temp}{abs(\pgfplotsretval-0.6)==0? 1 : 0}
\ifnum\temp>0
\relax
			\else
  \def\pgfmathresult{}
\fi
\else
  \def\pgfmathresult{}
\fi
\else
  \def\pgfmathresult{}
\fi
},
		]
                table[ 
                x index = 3,
                y index = 4,
                ]
		{\4};
  \addplot[
		 color=red!50!black,
		 only marks,
     mark=o,
     mark size = 2pt,
		 x filter/.code={\pgfplotstablegetelem{\coordindex}{[index]0}\of{\5}
\pgfmathtruncatemacro{\temp}{abs(\pgfplotsretval-20)==0? 1 : 0}
\ifnum\temp>0
\pgfplotstablegetelem{\coordindex}{[index]1}\of{\5}
\pgfmathtruncatemacro{\temp}{abs(\pgfplotsretval-0.6)==0? 1 : 0}
\ifnum\temp>0
\pgfplotstablegetelem{\coordindex}{[index]2}\of{\5}
\pgfmathtruncatemacro{\temp}{abs(\pgfplotsretval-0.6)==0? 1 : 0}
\ifnum\temp>0
\relax
			\else
  \def\pgfmathresult{}
\fi
\else
  \def\pgfmathresult{}
\fi
\else
  \def\pgfmathresult{}
\fi
},
		]
                table[ 
                x index = 3,
                y index = 4,
                ]
		{\5};
		\addplot[
		 color=green!50!black,
		 only marks,
     mark=o,
     mark size = 2pt,
		 x filter/.code={\pgfplotstablegetelem{\coordindex}{[index]0}\of{\7}
\pgfmathtruncatemacro{\temp}{abs(\pgfplotsretval-20)==0? 1 : 0}
\ifnum\temp>0
\pgfplotstablegetelem{\coordindex}{[index]1}\of{\7}
\pgfmathtruncatemacro{\temp}{abs(\pgfplotsretval-0.6)==0? 1 : 0}
\ifnum\temp>0
\pgfplotstablegetelem{\coordindex}{[index]2}\of{\7}
\pgfmathtruncatemacro{\temp}{abs(\pgfplotsretval-0.6)==0? 1 : 0}
\ifnum\temp>0
\relax
			\else
  \def\pgfmathresult{}
\fi
\else
  \def\pgfmathresult{}
\fi
\else
  \def\pgfmathresult{}
\fi
},
		]
                table[ 
                x index = 3,
                y index = 4,
                ]
		{\7};
\end{axis}
\end{tikzpicture}}
  \end{minipage}
  \caption{First $k$ principal angles, $\phi_i$ between the accurate subspace and the subspace obtained using the iterative algorithm.
 ICS:Column sampling, ICRS:column and row sampling
} 
  \label{fig:pc_iter}
\end{figure}

We also ran the algorithm with convergence of the subspace spanned by modes as the criterion with the same value for $\tau$.
Fig.~\ref{fig:pc_iter} shows the principal angles between  the subspaces spanned by the approximate modes and the modes obtained using the truncated SVD. 
The algorithms give similar results as in Fig.~\ref{fig:angPOD_iter} for all datasets. 
Moreover, it is clear that the accuracies are similar for all the sampling strategies, indicating
it is sufficient to use uniform distribution after the first round of sampling. This is because UNF is the least complex and 
computationally more efficient than the other column sampling strategies.

To evaluate the quality of subspaces obtained, we computed the measure for error in the subspaces 
using equation~(\ref{eqn:ang_meas}).  
Table~\ref{tab:ang_meas} shows the error measure for various cases. In general, the error measure is small indicating good accuracies in the modes. 
In a few cases ($V_{2D}$ ($k=2$) and YF ($k=20$)), the measure is large. 
As discussed, this could happen if $\sigma_k$ and $\sigma_{k+1}$ are clustered together. Small changes in the singular values can result in 
large variations in the measure ( as seen with L2N and UNF). Fig.~\ref{fig:sigma} in Appendix~\ref{app:eval} shows the
dominant singular values of all our test datasets. It can be seen that $\sigma_2$ and $\sigma_3$ of $V_{2D}$ are very close in value, whereas $\sigma_4$ is 
well separated. Therefore when $k$ is increased to $3$,
the error measure reduces significantly. 
This can be seen in Table~\ref{tab:inc_k}. 

The other case for which the measure is large is for YF (UNF). In this case also, several singular values are clustered together as seen in 
Fig.~\ref{fig:sigma}. The error drops when $k$ is increased as seen in Table~\ref{tab:inc_k}. If L2N is used for sampling, the measure is small. It turns that all the 
columns were sampled in this case, leading to a low error.
In general, it can be seen from Figs.~\ref{fig:angPOD_iter} and \ref{fig:pc_iter} that the accuracy of the modes are much better than indicated by the measure.
Also note that, in most cases, our approximation of 
$\omega$ worked well. The average error in $\omega$ was found to be around 6\%.

\begin{table}[tbp!]
  \centering
  \caption{Error measure computed using POD vectors obtained from Iterative algorithms with convergence of individual modes.}
  \label{tab:err_meas}
  \begin{subtable}{\textwidth}
  \centering
  \subcaption{$\sqrt{||R||_F^2+||S||_F^2}/\hat{\omega}$   where $\hat{\omega}\triangleq\min{\{|\tilde{\sigma}_k-\tilde{\sigma}_{k+1}|,\tilde{\sigma}_k\}}$. } 
  \label{tab:ang_meas}
    \begin{tabular}{cccccccccc}
\toprule
    $k$&Dataset&\multicolumn{3}{c}{ICS}&\multicolumn{3}{c}{ICRS}&$\sqrt{2k}$&$\frac{\sigma_k-\sigma_{k+1}}{\sigma_k}$\\
&&L2N&UNF&ORT&UNF&ORT&LS\\
\midrule
    10&faces&0.275&0.252&0.282&0.657&0.677&0.62&4.47&0.107\\
    10&CYF&0.731&0.567&0.665&0.868&0.943&0.748&4.47&0.024\\
    10&$V_{2D}$&0.119&0.123&0.096&0.158&0.132&0.154&4.47&0.11\\
    20&YF&0.428&\textbf{5.053}&-&\textbf{6.48}&-&0.347&\textbf{5.48}&0.026\\
\midrule
    2&faces&0.191&0.23&0.22&0.081&0.25&0.268&2&0.272\\
    2&CYF&0.07&0.055&0.083&0.055&0.097&0.051&2&0.624\\
    2&$V_{2D}$&0.762&\textbf{1.95}&\textbf{2.88}&\textbf{1.278}&\textbf{5.34}&1.115&\textbf{2}&0.022\\
    5&YF&2.08&1.71&-&1.394&-&1.43&3.16&0.073\\
\bottomrule
  \end{tabular} 
  \end{subtable}
  \begin{subtable}{\textwidth}
  \centering
  \subcaption{Error measure with value $k$ increased by one for $V_{2D}$ and $YF$.}
  \label{tab:inc_k}
    \begin{tabular}{cccccccc}
\toprule
    $k$&Dataset&\multicolumn{3}{c}{ICS}&\multicolumn{3}{c}{ICRS}\\
&&L2N&UNF&ORT&UNF&ORT&LS\\
  \midrule
    3&$V_{2D}$&0.237&0.32&0.42&0.356&0.523&0.453\\
    21&YF&0.207&0.88&-&0.88&-&1.46\\
\bottomrule
  \end{tabular}
  \end{subtable}
\end{table}

\begin{figure}[!htbp]
 \begin{minipage}{\textwidth}
   \centering
\ref{spdup_iter}\\
  \end{minipage}
  \begin{minipage}{\textwidth}
    \centering
    Single thread\\
     \scalebox{0.7}{
       \pgfplotstableread{data/speedup_iterate_V_2D_1threads.dat}\1
\pgfplotstableread{data/time_SVDATA_V_2D_1threads.dat}\2
\pgfplotstableset{
     create on use/0/.style={create col/copy column from table={\1}{0}}, 
     create on use/1/.style={create col/copy column from table={\1}{1}}, 
     create on use/2/.style={create col/copy column from table={\1}{2}}, 
     create on use/3/.style={create col/copy column from table={\1}{3}}, 
     create on use/4/.style={create col/copy column from table={\1}{4}}, 
     create on use/5/.style={create col/copy column from table={\1}{5}}, 
     create on use/6/.style={create col/copy column from table={\1}{6}}, 
     create on use/8/.style={create col/copy column from table={\2}{1}}, 
     create on use/9/.style={create col/copy column from table={\2}{2}}, 
}
\pgfplotstablenew[columns={0,1,2,3,4,5,6,8,9}]{\pgfplotstablegetrowsof{\1}}\3
\begin{tikzpicture}
\begin{axis}[
    width = 6cm,
    title = {$V_{2D}$},
    ybar,
    bar width = 6pt,
    enlarge x limits=0.50,
    ylabel={speed up},
    xtick = {1,2},
    xticklabels = {$k=2$,$k=10$},
    nodes near coords,
    nodes near coords align = {west},
    every node near coord/.append style = {rotate=90,  font=\scriptsize, /pgf/number format/.cd, fixed, precision=1},
    ymax = 60,
    ymin = 0,
    yticklabels = {,,},
    legend entries = {$\frac{t_{\text{SVD-}A^TA}}{t_{\text{ICS-L2N}}}$, $\frac{t_{\text{SVD-}A^TA}}{t_{\text{ICS-UNF}}}$,$\frac{t_{\text{SVD-}A^TA}}{t_{\text{ICS-ORT}}}$, $\frac{t_{\text{SVD-}A^TA}}{t_{\text{ICRS-UNF}}}$,$\frac{t_{\text{SVD-}A^TA}}{t_{\text{ICRS-ORT}}}$,$\frac{t_{\text{SVD-}A^TA}}{t_{\text{ICRS-LS}}}$},
    legend style={legend columns = -1},
    legend to name = spdup_iter,
]
  \addplot table[x =0, y expr = \thisrow{6}*\thisrow{8}/\thisrow{9}]{\3};
  \addplot table[x =0, y expr = \thisrow{3}*\thisrow{8}/\thisrow{9}]{\3};
  \addplot table[x =0, y expr = \thisrow{1}*\thisrow{8}/\thisrow{9}]{\3};
  \addplot table[x =0, y expr = \thisrow{4}*\thisrow{8}/\thisrow{9}]{\3};
  \addplot table[x =0, y expr = \thisrow{2}*\thisrow{8}/\thisrow{9}]{\3};
  \addplot table[x =0, y expr = \thisrow{5}*\thisrow{8}/\thisrow{9}]{\3};
\end{axis}
\end{tikzpicture}
}
     \scalebox{0.7}{
       \pgfplotstableread{data/speedup_iterate_faces_1threads.dat}\1
\pgfplotstableread{data/time_SVDATA_faces_1threads.dat}\2
\pgfplotstableset{
     create on use/0/.style={create col/copy column from table={\1}{0}}, 
     create on use/1/.style={create col/copy column from table={\1}{1}}, 
     create on use/2/.style={create col/copy column from table={\1}{2}}, 
     create on use/3/.style={create col/copy column from table={\1}{3}}, 
     create on use/4/.style={create col/copy column from table={\1}{4}}, 
     create on use/5/.style={create col/copy column from table={\1}{5}}, 
     create on use/6/.style={create col/copy column from table={\1}{6}}, 
     create on use/8/.style={create col/copy column from table={\2}{1}}, 
     create on use/9/.style={create col/copy column from table={\2}{2}}, 
}
\pgfplotstablenew[columns={0,1,2,3,4,5,6,8,9}]{\pgfplotstablegetrowsof{\1}}\3
\begin{tikzpicture}
\begin{axis}[
    width = 6cm,
    title = {Faces dataset},
    ybar,
    bar width = 6pt,
    enlarge x limits=0.50,
    xtick = {1,2},
    xticklabels = {$k=2$,$k=10$},
    nodes near coords,
    nodes near coords align = {west},
    every node near coord/.append style = {rotate=90,  font=\scriptsize, /pgf/number format/.cd, fixed, precision=1},
    ymax = 60,
    ymin = 0,
    yticklabels={,,},
]
  \addplot table[x =0, y expr = \thisrow{6}*\thisrow{8}/\thisrow{9}]{\3};
  \addplot table[x =0, y expr = \thisrow{3}*\thisrow{8}/\thisrow{9}]{\3};
  \addplot table[x =0, y expr = \thisrow{1}*\thisrow{8}/\thisrow{9}]{\3};
  \addplot table[x =0, y expr = \thisrow{4}*\thisrow{8}/\thisrow{9}]{\3};
  \addplot table[x =0, y expr = \thisrow{2}*\thisrow{8}/\thisrow{9}]{\3};
  \addplot table[x =0, y expr = \thisrow{5}*\thisrow{8}/\thisrow{9}]{\3};
\end{axis}
\end{tikzpicture}
}
      \scalebox{0.7}{
        \pgfplotstableread{data/speedup_iterate_croppedYalefaces_1threads.dat}\1
\pgfplotstableread{data/time_SVDATA_croppedYalefaces_1threads.dat}\2
\pgfplotstableset{
     create on use/0/.style={create col/copy column from table={\1}{0}}, 
     create on use/1/.style={create col/copy column from table={\1}{1}}, 
     create on use/2/.style={create col/copy column from table={\1}{2}}, 
     create on use/3/.style={create col/copy column from table={\1}{3}}, 
     create on use/4/.style={create col/copy column from table={\1}{4}}, 
     create on use/5/.style={create col/copy column from table={\1}{5}}, 
     create on use/6/.style={create col/copy column from table={\1}{6}}, 
     create on use/8/.style={create col/copy column from table={\2}{1}}, 
     create on use/9/.style={create col/copy column from table={\2}{2}}, 
}
\pgfplotstablenew[columns={0,1,2,3,4,5,6,8,9}]{\pgfplotstablegetrowsof{\1}}\3
\begin{tikzpicture}
\begin{axis}[
    width = 6cm,
    title = {Croppped Yale faces},
    ybar,
    bar width = 6pt,
    enlarge x limits=0.50,
    xtick = {1,2},
    xticklabels = {$k=2$,$k=10$},
    nodes near coords,
    nodes near coords align = {west},
    every node near coord/.append style = {rotate=90,  font=\scriptsize, /pgf/number format/.cd, fixed, precision=1},
    ymax = 60,
    ymin = 0,
    yticklabels = {,,}
]
  \addplot table[x =0, y expr = \thisrow{6}*\thisrow{8}/\thisrow{9}]{\3};
  \addplot table[x =0, y expr = \thisrow{3}*\thisrow{8}/\thisrow{9}]{\3};
  \addplot table[x =0, y expr = \thisrow{1}*\thisrow{8}/\thisrow{9}]{\3};
  \addplot table[x =0, y expr = \thisrow{4}*\thisrow{8}/\thisrow{9}]{\3};
  \addplot table[x =0, y expr = \thisrow{2}*\thisrow{8}/\thisrow{9}]{\3};
  \addplot table[x =0, y expr = \thisrow{5}*\thisrow{8}/\thisrow{9}]{\3};
\end{axis}
\end{tikzpicture}}\\
     \end{minipage}
  \begin{minipage}{\textwidth}
    \centering
    4 threads\\
     \scalebox{0.7}{
       \pgfplotstableread{data/speedup_iterate_V_2D_4threads.dat}\1
\pgfplotstableread{data/time_SVDATA_V_2D_4threads.dat}\2
\pgfplotstableset{
     create on use/0/.style={create col/copy column from table={\1}{0}}, 
     create on use/1/.style={create col/copy column from table={\1}{1}}, 
     create on use/2/.style={create col/copy column from table={\1}{2}}, 
     create on use/3/.style={create col/copy column from table={\1}{3}}, 
     create on use/4/.style={create col/copy column from table={\1}{4}}, 
     create on use/5/.style={create col/copy column from table={\1}{5}}, 
     create on use/6/.style={create col/copy column from table={\1}{6}}, 
     create on use/8/.style={create col/copy column from table={\2}{1}}, 
     create on use/9/.style={create col/copy column from table={\2}{2}}, 
}
\pgfplotstablenew[columns={0,1,2,3,4,5,6,8,9}]{\pgfplotstablegetrowsof{\1}}\3
\begin{tikzpicture}
\begin{axis}[
    width = 6cm,
    title = {$V_{2D}$},
    ybar,
    bar width = 6pt,
    enlarge x limits=0.50,
    ylabel={speed up},
    xtick = {1,2},
    xticklabels = {$k=2$,$k=10$},
    nodes near coords,
    nodes near coords align = {west},
    every node near coord/.append style = {rotate=90,  font=\scriptsize, /pgf/number format/.cd, fixed, precision=1},
    ymax = 60,
    ymin = 0,
    yticklabels = {,,},
]
  \addplot table[x =0, y expr = \thisrow{6}*\thisrow{8}/\thisrow{9}]{\3};
  \addplot table[x =0, y expr = \thisrow{3}*\thisrow{8}/\thisrow{9}]{\3};
  \addplot table[x =0, y expr = \thisrow{1}*\thisrow{8}/\thisrow{9}]{\3};
  \addplot table[x =0, y expr = \thisrow{4}*\thisrow{8}/\thisrow{9}]{\3};
  \addplot table[x =0, y expr = \thisrow{2}*\thisrow{8}/\thisrow{9}]{\3};
  \addplot table[x =0, y expr = \thisrow{5}*\thisrow{8}/\thisrow{9}]{\3};
\end{axis}
\end{tikzpicture}
}
      \scalebox{0.7}{
       \pgfplotstableread{data/speedup_iterate_faces_4threads.dat}\1
\pgfplotstableread{data/time_SVDATA_faces_4threads.dat}\2
\pgfplotstableset{
     create on use/0/.style={create col/copy column from table={\1}{0}}, 
     create on use/1/.style={create col/copy column from table={\1}{1}}, 
     create on use/2/.style={create col/copy column from table={\1}{2}}, 
     create on use/3/.style={create col/copy column from table={\1}{3}}, 
     create on use/4/.style={create col/copy column from table={\1}{4}}, 
     create on use/5/.style={create col/copy column from table={\1}{5}}, 
     create on use/6/.style={create col/copy column from table={\1}{6}}, 
     create on use/8/.style={create col/copy column from table={\2}{1}}, 
     create on use/9/.style={create col/copy column from table={\2}{2}}, 
}
\pgfplotstablenew[columns={0,1,2,3,4,5,6,8,9}]{\pgfplotstablegetrowsof{\1}}\3
\begin{tikzpicture}
\begin{axis}[
    width = 6cm,
    title = {Faces dataset},
    ybar,
    bar width = 6pt,
    enlarge x limits=0.50,
    xtick = {1,2},
    xticklabels = {$k=2$,$k=10$},
    nodes near coords,
    nodes near coords align = {west},
    every node near coord/.append style = {rotate=90,  font=\scriptsize, /pgf/number format/.cd, fixed, precision=1},
    ymax = 60,
    ymin = 0,
    yticklabels = {,,},
]
  \addplot table[x =0, y expr = \thisrow{6}*\thisrow{8}/\thisrow{9}]{\3};
  \addplot table[x =0, y expr = \thisrow{3}*\thisrow{8}/\thisrow{9}]{\3};
  \addplot table[x =0, y expr = \thisrow{1}*\thisrow{8}/\thisrow{9}]{\3};
  \addplot table[x =0, y expr = \thisrow{4}*\thisrow{8}/\thisrow{9}]{\3};
  \addplot table[x =0, y expr = \thisrow{2}*\thisrow{8}/\thisrow{9}]{\3};
  \addplot table[x =0, y expr = \thisrow{5}*\thisrow{8}/\thisrow{9}]{\3};
\end{axis}
\end{tikzpicture}
}
      \scalebox{0.7}{
       \pgfplotstableread{data/speedup_iterate_croppedYalefaces_4threads.dat}\1
\pgfplotstableread{data/time_SVDATA_croppedYalefaces_4threads.dat}\2
\pgfplotstableset{
     create on use/0/.style={create col/copy column from table={\1}{0}}, 
     create on use/1/.style={create col/copy column from table={\1}{1}}, 
     create on use/2/.style={create col/copy column from table={\1}{2}}, 
     create on use/3/.style={create col/copy column from table={\1}{3}}, 
     create on use/4/.style={create col/copy column from table={\1}{4}}, 
     create on use/5/.style={create col/copy column from table={\1}{5}}, 
     create on use/6/.style={create col/copy column from table={\1}{6}}, 
     create on use/8/.style={create col/copy column from table={\2}{1}}, 
     create on use/9/.style={create col/copy column from table={\2}{2}}, 
}
\pgfplotstablenew[columns={0,1,2,3,4,5,6,8,9}]{\pgfplotstablegetrowsof{\1}}\3
\begin{tikzpicture}
\begin{axis}[
    width = 6cm,
    title = {Croppped Yale faces},
    ybar,
    bar width = 6pt,
    enlarge x limits=0.50,
    xtick = {1,2},
    xticklabels = {$k=2$,$k=10$},
    nodes near coords,
    nodes near coords align = {west},
    every node near coord/.append style = {rotate=90,  font=\scriptsize, /pgf/number format/.cd, fixed, precision=1},
    ymax = 60,
    ymin = 0,
    yticklabels = {,,}
]
  \addplot table[x =0, y expr = \thisrow{6}*\thisrow{8}/\thisrow{9}]{\3};
  \addplot table[x =0, y expr = \thisrow{3}*\thisrow{8}/\thisrow{9}]{\3};
  \addplot table[x =0, y expr = \thisrow{1}*\thisrow{8}/\thisrow{9}]{\3};
  \addplot table[x =0, y expr = \thisrow{4}*\thisrow{8}/\thisrow{9}]{\3};
  \addplot table[x =0, y expr = \thisrow{2}*\thisrow{8}/\thisrow{9}]{\3};
  \addplot table[x =0, y expr = \thisrow{5}*\thisrow{8}/\thisrow{9}]{\3};
\end{axis}
\end{tikzpicture}}
  \end{minipage}
  \caption{Speed up for iterative algorithms based on convergence of individual modes. ICS:Column sampling, ICRS:Column and row sampling
  }
  \label{fig:speedup_iter}
\end{figure}
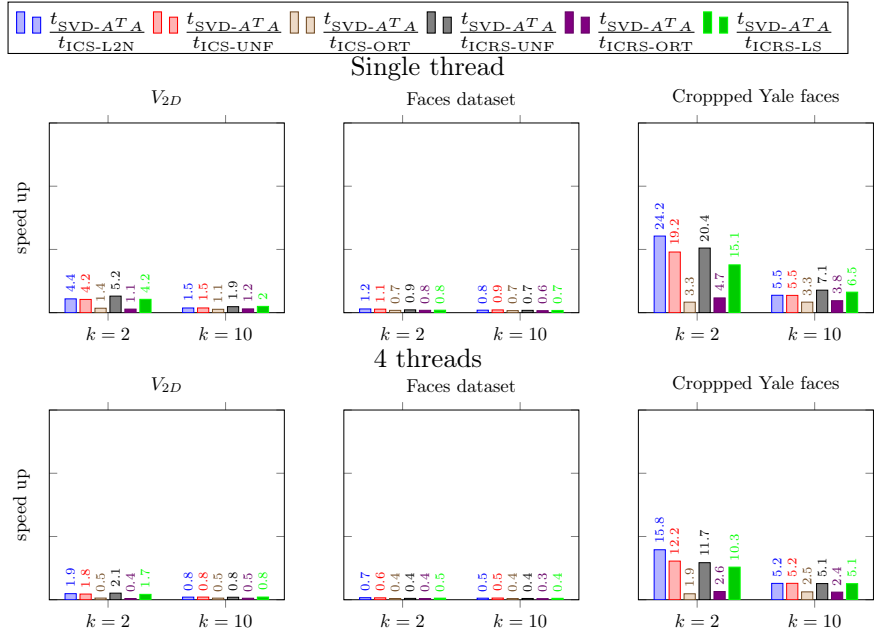
\begin{figure}[!htbp]
 \begin{minipage}{\textwidth}
   \centering
\ref{spdup_iter_pc}\\
  \end{minipage}
  \begin{minipage}{\textwidth}
    \centering
    Single thread\\
     \scalebox{0.7}{
       \pgfplotstableread{data/speedup_iterate_pc_V_2D_1threads.dat}\1
\pgfplotstableread{data/time_SVDATA_V_2D_1threads.dat}\2
\pgfplotstableset{
     create on use/0/.style={create col/copy column from table={\1}{0}}, 
     create on use/1/.style={create col/copy column from table={\1}{1}}, 
     create on use/2/.style={create col/copy column from table={\1}{2}}, 
     create on use/3/.style={create col/copy column from table={\1}{3}}, 
     create on use/4/.style={create col/copy column from table={\1}{4}}, 
     create on use/5/.style={create col/copy column from table={\1}{5}}, 
     create on use/6/.style={create col/copy column from table={\1}{6}}, 
     create on use/8/.style={create col/copy column from table={\2}{1}}, 
     create on use/9/.style={create col/copy column from table={\2}{2}}, 
}
\pgfplotstablenew[columns={0,1,2,3,4,5,6,8,9}]{\pgfplotstablegetrowsof{\1}}\3
\begin{tikzpicture}
\begin{axis}[
    width = 6cm,
    title = {$V_{2D}$},
    ybar,
    bar width = 6pt,
    enlarge x limits=0.50,
    ylabel={speed up},
    xtick = {1,2},
    xticklabels = {$k=2$,$k=10$},
    nodes near coords,
    nodes near coords align = {west},
    every node near coord/.append style = {rotate=90,  font=\scriptsize, /pgf/number format/.cd, fixed, precision=1},
    ymax = 60,
    ymin = 0,
    yticklabels = {,,},
    legend entries = { $\frac{t_{\text{SVD-}A^TA}}{t_{\text{ICS-L2N}}}$,$\frac{t_{\text{SVD-}A^TA}}{t_{\text{ICS-UNF}}}$,$\frac{t_{\text{SVD-}A^TA}}{t_{\text{ICS-ORT}}}$, $\frac{t_{\text{SVD-}A^TA}}{t_{\text{ICRS-UNF}}}$,$\frac{t_{\text{SVD-}A^TA}}{t_{\text{ICRS-ORT}}}$,$\frac{t_{\text{SVD-}A^TA}}{t_{\text{ICRS-LS}}}$},
    legend style={legend columns = -1},
    legend to name = spdup_iter_pc,
]
  \addplot table[x =0, y expr = \thisrow{6}*\thisrow{8}/\thisrow{9}]{\3};
  \addplot table[x =0, y expr = \thisrow{3}*\thisrow{8}/\thisrow{9}]{\3};
  \addplot table[x =0, y expr = \thisrow{1}*\thisrow{8}/\thisrow{9}]{\3};
  \addplot table[x =0, y expr = \thisrow{4}*\thisrow{8}/\thisrow{9}]{\3};
  \addplot table[x =0, y expr = \thisrow{2}*\thisrow{8}/\thisrow{9}]{\3};
  \addplot table[x =0, y expr = \thisrow{5}*\thisrow{8}/\thisrow{9}]{\3};
\end{axis}
\end{tikzpicture}
}
     \scalebox{0.7}{
       \pgfplotstableread{data/speedup_iterate_pc_faces_1threads.dat}\1
\pgfplotstableread{data/time_SVDATA_faces_1threads.dat}\2
\pgfplotstableset{
     create on use/0/.style={create col/copy column from table={\1}{0}}, 
     create on use/1/.style={create col/copy column from table={\1}{1}}, 
     create on use/2/.style={create col/copy column from table={\1}{2}}, 
     create on use/3/.style={create col/copy column from table={\1}{3}}, 
     create on use/4/.style={create col/copy column from table={\1}{4}}, 
     create on use/5/.style={create col/copy column from table={\1}{5}}, 
     create on use/6/.style={create col/copy column from table={\1}{6}}, 
     create on use/8/.style={create col/copy column from table={\2}{1}}, 
     create on use/9/.style={create col/copy column from table={\2}{2}}, 
}
\pgfplotstablenew[columns={0,1,2,3,4,5,6,8,9}]{\pgfplotstablegetrowsof{\1}}\3
\begin{tikzpicture}
\begin{axis}[
    width = 6cm,
    title = {Faces dataset},
    ybar,
    bar width = 6pt,
    enlarge x limits=0.50,
    xtick = {1,2},
    xticklabels = {$k=2$,$k=10$},
    nodes near coords,
    nodes near coords align = {west},
    every node near coord/.append style = {rotate=90,  font=\scriptsize, /pgf/number format/.cd, fixed, precision=1},
    ymax = 60,
    ymin = 0,
    yticklabels={,,},
]
  \addplot table[x =0, y expr = \thisrow{6}*\thisrow{8}/\thisrow{9}]{\3};
  \addplot table[x =0, y expr = \thisrow{3}*\thisrow{8}/\thisrow{9}]{\3};
  \addplot table[x =0, y expr = \thisrow{1}*\thisrow{8}/\thisrow{9}]{\3};
  \addplot table[x =0, y expr = \thisrow{4}*\thisrow{8}/\thisrow{9}]{\3};
  \addplot table[x =0, y expr = \thisrow{2}*\thisrow{8}/\thisrow{9}]{\3};
  \addplot table[x =0, y expr = \thisrow{5}*\thisrow{8}/\thisrow{9}]{\3};
\end{axis}
\end{tikzpicture}
}
      \scalebox{0.7}{
        \pgfplotstableread{data/speedup_iterate_pc_croppedYalefaces_1threads.dat}\1
\pgfplotstableread{data/time_SVDATA_croppedYalefaces_1threads.dat}\2
\pgfplotstableset{
     create on use/0/.style={create col/copy column from table={\1}{0}}, 
     create on use/1/.style={create col/copy column from table={\1}{1}}, 
     create on use/2/.style={create col/copy column from table={\1}{2}}, 
     create on use/3/.style={create col/copy column from table={\1}{3}}, 
     create on use/4/.style={create col/copy column from table={\1}{4}}, 
     create on use/5/.style={create col/copy column from table={\1}{5}}, 
     create on use/6/.style={create col/copy column from table={\1}{6}}, 
     create on use/8/.style={create col/copy column from table={\2}{1}}, 
     create on use/9/.style={create col/copy column from table={\2}{2}}, 
}
\pgfplotstablenew[columns={0,1,2,3,4,5,6,8,9}]{\pgfplotstablegetrowsof{\1}}\3
\begin{tikzpicture}
\begin{axis}[
    width = 6cm,
    title = {Croppped Yale faces},
    ybar,
    bar width = 6pt,
    enlarge x limits=0.50,
    xtick = {1,2},
    xticklabels = {$k=2$,$k=10$},
    nodes near coords,
    nodes near coords align = {west},
    every node near coord/.append style = {rotate=90,  font=\scriptsize, /pgf/number format/.cd, fixed, precision=1},
    ymax = 60,
    ymin = 0,
    yticklabels = {,,}
]
  \addplot table[x =0, y expr = \thisrow{6}*\thisrow{8}/\thisrow{9}]{\3};
  \addplot table[x =0, y expr = \thisrow{3}*\thisrow{8}/\thisrow{9}]{\3};
  \addplot table[x =0, y expr = \thisrow{1}*\thisrow{8}/\thisrow{9}]{\3};
  \addplot table[x =0, y expr = \thisrow{4}*\thisrow{8}/\thisrow{9}]{\3};
  \addplot table[x =0, y expr = \thisrow{2}*\thisrow{8}/\thisrow{9}]{\3};
  \addplot table[x =0, y expr = \thisrow{5}*\thisrow{8}/\thisrow{9}]{\3};
\end{axis}
\end{tikzpicture}
}\\
     \end{minipage}
  \begin{minipage}{\textwidth}
    \centering
    4 threads\\
     \scalebox{0.7}{
       \pgfplotstableread{data/speedup_iterate_pc_V_2D_4threads.dat}\1
\pgfplotstableread{data/time_SVDATA_V_2D_4threads.dat}\2
\pgfplotstableset{
     create on use/0/.style={create col/copy column from table={\1}{0}}, 
     create on use/1/.style={create col/copy column from table={\1}{1}}, 
     create on use/2/.style={create col/copy column from table={\1}{2}}, 
     create on use/3/.style={create col/copy column from table={\1}{3}}, 
     create on use/4/.style={create col/copy column from table={\1}{4}}, 
     create on use/5/.style={create col/copy column from table={\1}{5}}, 
     create on use/6/.style={create col/copy column from table={\1}{6}}, 
     create on use/8/.style={create col/copy column from table={\2}{1}}, 
     create on use/9/.style={create col/copy column from table={\2}{2}}, 
}
\pgfplotstablenew[columns={0,1,2,3,4,5,6,8,9}]{\pgfplotstablegetrowsof{\1}}\3
\begin{tikzpicture}
\begin{axis}[
    width = 6cm,
    title = {$V_{2D}$},
    ybar,
    bar width = 6pt,
    enlarge x limits=0.50,
    ylabel={speed up},
    xtick = {1,2},
    xticklabels = {$k=2$,$k=10$},
    nodes near coords,
    nodes near coords align = {west},
    every node near coord/.append style = {rotate=90,  font=\scriptsize, /pgf/number format/.cd, fixed, precision=1},
    ymax = 60,
    ymin = 0,
    yticklabels = {,,},
]
  \addplot table[x =0, y expr = \thisrow{6}*\thisrow{8}/\thisrow{9}]{\3};
  \addplot table[x =0, y expr = \thisrow{3}*\thisrow{8}/\thisrow{9}]{\3};
  \addplot table[x =0, y expr = \thisrow{1}*\thisrow{8}/\thisrow{9}]{\3};
  \addplot table[x =0, y expr = \thisrow{4}*\thisrow{8}/\thisrow{9}]{\3};
  \addplot table[x =0, y expr = \thisrow{2}*\thisrow{8}/\thisrow{9}]{\3};
  \addplot table[x =0, y expr = \thisrow{5}*\thisrow{8}/\thisrow{9}]{\3};
\end{axis}
\end{tikzpicture}
}
      \scalebox{0.7}{
       \pgfplotstableread{data/speedup_iterate_pc_faces_4threads.dat}\1
\pgfplotstableread{data/time_SVDATA_faces_4threads.dat}\2
\pgfplotstableset{
     create on use/0/.style={create col/copy column from table={\1}{0}}, 
     create on use/1/.style={create col/copy column from table={\1}{1}}, 
     create on use/2/.style={create col/copy column from table={\1}{2}}, 
     create on use/3/.style={create col/copy column from table={\1}{3}}, 
     create on use/4/.style={create col/copy column from table={\1}{4}}, 
     create on use/5/.style={create col/copy column from table={\1}{5}}, 
     create on use/6/.style={create col/copy column from table={\1}{6}}, 
     create on use/8/.style={create col/copy column from table={\2}{1}}, 
     create on use/9/.style={create col/copy column from table={\2}{2}}, 
}
\pgfplotstablenew[columns={0,1,2,3,4,5,6,8,9}]{\pgfplotstablegetrowsof{\1}}\3
\begin{tikzpicture}
\begin{axis}[
    width = 6cm,
    title = {Faces dataset},
    ybar,
    bar width = 6pt,
    enlarge x limits=0.50,
    xtick = {1,2},
    xticklabels = {$k=2$,$k=10$},
    nodes near coords,
    nodes near coords align = {west},
    every node near coord/.append style = {rotate=90,  font=\scriptsize, /pgf/number format/.cd, fixed, precision=1},
    ymax = 60,
    ymin = 0,
    yticklabels = {,,},
]
  \addplot table[x =0, y expr = \thisrow{6}*\thisrow{8}/\thisrow{9}]{\3};
  \addplot table[x =0, y expr = \thisrow{3}*\thisrow{8}/\thisrow{9}]{\3};
  \addplot table[x =0, y expr = \thisrow{1}*\thisrow{8}/\thisrow{9}]{\3};
  \addplot table[x =0, y expr = \thisrow{4}*\thisrow{8}/\thisrow{9}]{\3};
  \addplot table[x =0, y expr = \thisrow{2}*\thisrow{8}/\thisrow{9}]{\3};
  \addplot table[x =0, y expr = \thisrow{5}*\thisrow{8}/\thisrow{9}]{\3};
\end{axis}
\end{tikzpicture}
}
      \scalebox{0.7}{
       \pgfplotstableread{data/speedup_iterate_pc_croppedYalefaces_4threads.dat}\1
\pgfplotstableread{data/time_SVDATA_croppedYalefaces_4threads.dat}\2
\pgfplotstableset{
     create on use/0/.style={create col/copy column from table={\1}{0}}, 
     create on use/1/.style={create col/copy column from table={\1}{1}}, 
     create on use/2/.style={create col/copy column from table={\1}{2}}, 
     create on use/3/.style={create col/copy column from table={\1}{3}}, 
     create on use/4/.style={create col/copy column from table={\1}{4}}, 
     create on use/5/.style={create col/copy column from table={\1}{5}}, 
     create on use/6/.style={create col/copy column from table={\1}{6}}, 
     create on use/8/.style={create col/copy column from table={\2}{1}}, 
     create on use/9/.style={create col/copy column from table={\2}{2}}, 
}
\pgfplotstablenew[columns={0,1,2,3,4,5,6,8,9}]{\pgfplotstablegetrowsof{\1}}\3
\begin{tikzpicture}
\begin{axis}[
    width = 6cm,
    title = {Croppped Yale faces},
    ybar,
    bar width = 6pt,
    enlarge x limits=0.50,
    xtick = {1,2},
    xticklabels = {$k=2$,$k=10$},
    nodes near coords,
    nodes near coords align = {west},
    every node near coord/.append style = {rotate=90,  font=\scriptsize, /pgf/number format/.cd, fixed, precision=1},
    ymax = 60,
    ymin = 0,
    yticklabels = {,,}
]
  \addplot table[x =0, y expr = \thisrow{6}*\thisrow{8}/\thisrow{9}]{\3};
  \addplot table[x =0, y expr = \thisrow{3}*\thisrow{8}/\thisrow{9}]{\3};
  \addplot table[x =0, y expr = \thisrow{1}*\thisrow{8}/\thisrow{9}]{\3};
  \addplot table[x =0, y expr = \thisrow{4}*\thisrow{8}/\thisrow{9}]{\3};
  \addplot table[x =0, y expr = \thisrow{2}*\thisrow{8}/\thisrow{9}]{\3};
  \addplot table[x =0, y expr = \thisrow{5}*\thisrow{8}/\thisrow{9}]{\3};
\end{axis}
\end{tikzpicture}}
  \end{minipage}
  \caption{Speed up for iterative algorithms based on convergence of subspaces formed by modes. ICS:Column sampling, ICRS:Column and row sampling
  }
  \label{fig:speedup_iter_pc}
\end{figure}
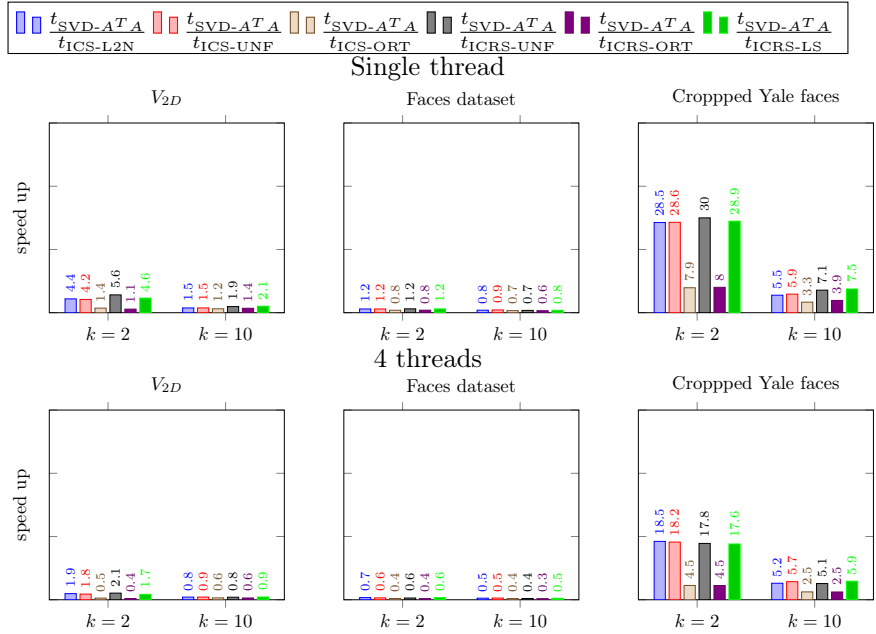

Fig.~\ref{fig:speedup_iter} and Fig.~\ref{fig:speedup_iter_pc} contain the speedup obtained  when the iterative algorithms are run with convergence of individual modes and subspaces formed by the modes as the convergence criterion, respectively. Only the smaller datasets are included. YF is discussed later.
Speedup is evaluated with respect to time taken for computing $k$ POD modes using SVD-$A^TA$.
We note that (a) very good speedup is obtained for CYF for both single and multi-threaded execution. 
The speedup is the least when ORT is used as the sampling strategy, which is as expected. Since the accuracies of all schemes are similar, it makes sense to use the UNF
strategy as it gives maximum speedup. (b) For small datasets such as the FACES dataset, it is better to use SVD-$A^TA$. (c) $V_{2D}$ has about 1000 columns. The results 
show that it may be worthwhile to use the iterative algorithm if $k$ is small.

We found that for $k=10$ (I1), nearly all columns were  sampled for the Faces, $V_{2D}$ and CYF datasets when convergence was achieved, both
in terms of mode and principal angles. 
For I2 ($k=2$), the number of columns sampled was substantially lower than $n$ in all cases and the speedup is much larger.
For the $V_{2D}$ and Faces datasets, the speedup is almost the same for both convergence criteria. 
However, for CYF, speed up obtained is substantially larger with convergence of subspace for $k=2$. This is because the first two singular values are almost identical. Therefore, more columns are sampled to capture individual modes accurately than for subspace spanned by the modes. 
 
As expected, in comparison to other sampling strategies, ORT is slower. UNF, LS, and L2N have very similar speed up for most cases.
This reinforces our conclusion that sampling with uniform sampling from the second iteration is as good as other sampling strategies we tried. 

As indicated earlier, mean centered dataset for Yale Faces occupies 40GB memory. 
We implemented Algorithm~\ref{alg:ICS} in a system with 64GB RAM and compared the results obtained using UNF with  
(a) SVD($A^TA$) and (b) projection algorithms proposed in 
\cite{RokSzlTyg:2009,HalMarTro:2011,EriBruKut:2017}.
We implemented two versions of projection algorithms: one in which the rows are projected onto the random matrix (RSVD) \cite{HalMarTro:2011} and the 
other in which the columns are projected (COLRSVD) \cite{RokSzlTyg:2009,EriBruKut:2017}.
We report results for RSVD and COLRSVD for the YF dataset  with 3 subspace iterations and $s=k+5$ as suggested in \cite{HalMarTro:2011}. 
The wall clock and computational times taken for YF dataset using these algorithms is shown in Table~\ref{tab:spd_compare}.
The wall clock time includes the time required to read the data from the disk. 
Fig.~\ref{fig:angPOD_PRSVD} shows the mode angles obtained.
We can see that for Yale faces dataset with $k=5$, the projection algorithms and sampling algorithms have similar accuracy. Sampling algorithms are $~1.4\times$
faster with respect to wall clock times and $~2\times$ faster with respect to computational times.
 For $k=20$, the wall clock times of projection and sampling algorithms are similar (except for ICS-UNF) but the projection algorithms have a large error in 
comparison to the iterative sampling algorithms.
Also note that, the iterative sampling algorithms are much 
faster ($7-30\times$) than truncated SVD.

\begin{table}[h]
\centering
\caption{Wall clock (WT) and computational times (CT) in seconds (s) for projection and iterative sampling algorithms for Yale faces dataset. Algorithms were run in a 64GB system since YF does not fit in our test system.}
\label{tab:spd_compare}
\begin{tabular}{cccccc}
\toprule
Case&Algorithm&\multicolumn{2}{c}{single thread}&\multicolumn{2}{c}{8 threads}\\
&&WT(s)&CT(s)&WT(s)&CT(s)\\
\midrule
\multirow{5}{*}{$k=5$ (I2)}&{ICS-UNF}&{217.79}&{50.78}&{224.062}&{38.57}\\
&{ICRS-UNF}&{205.71}&{42.11}&{222.226}&{41.73}\\
&RSVD&298.75&108.43&283.5&108.8\\
&COLRSVD&286.44&109.83&281.78&108.97\\
&SVD-$A^TA$&7151.07&6964.7&2149.77&1962.77\\
\midrule
\multirow{5}{*}{$k=20$ (I1)}&{ICS-UNF}&{626.18}&{445.34}&{299.05}&{130.59}\\
&{ICRS-UNF}&{379.16}&{198.19}&{278.74}&{111}\\
&RSVD&335.12&160.77&290.7&159.32\\
&COLRSVD&346.89&159.98&290.18&159.48\\
&SVD-$A^TA$&7403.33&7216.96&2171.17&1984.17\\
\bottomrule
\end{tabular}
\end{table}
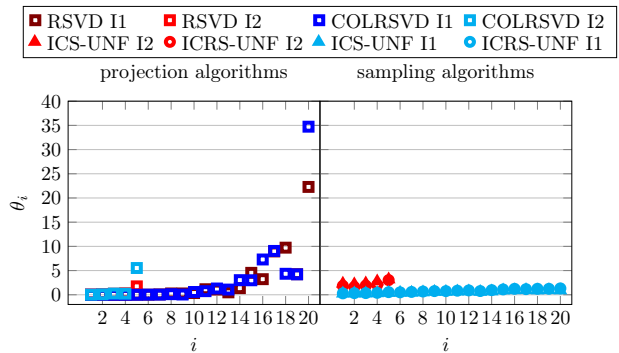
\begin{figure}[h]
    \centering
\scalebox{0.75}{ \pgfplotstableread{data/angPOD_RSVD_Yalefaces.dat}\1
\pgfplotstableread{data/angPOD_COLRSVD_Yalefaces.dat}\2
\pgfplotstableread{data/angPOD_LTSVD_US_no_scaling_iterate_Yalefaces.dat}\4
\pgfplotstableread{data/angPOD_CTSVD_US_no_scaling_iterate_Yalefaces.dat}\5
\begin{tikzpicture}[every mark/.append style={line width=2pt, solid}]
\begin{groupplot}[ 
	group style = {group size = 2 by 1,
  horizontal sep = 0cm,
  },
  xlabel = $i$,
  xmax = 21.000000,
  xtick={2,4,...,20},
  ytick={0,5,...,40},
  width = 6cm,
  ymajorgrids,
   ymax = 40,
  ymin =-2,
]
  \nextgroupplot[
  title = {projection algorithms}, 
  ylabel = $\theta_i$,
]
  \addplot[
		 color=red!50!black,
		 only marks,
		 mark=square,
     mark size = 2pt,
		 x filter/.code={\pgfplotstablegetelem{\coordindex}{[index]0}\of{\1}
\pgfmathtruncatemacro{\temp}{abs(\pgfplotsretval-20)==0? 1 : 0}
\ifnum\temp>0
\pgfplotstablegetelem{\coordindex}{[index]1}\of{\1}
\pgfmathtruncatemacro{\temp}{abs(\pgfplotsretval-3)==0? 1 : 0}
\ifnum\temp>0
\pgfplotstablegetelem{\coordindex}{[index]2}\of{\1}
\pgfmathtruncatemacro{\temp}{abs(\pgfplotsretval-1)==0? 1 : 0}
\ifnum\temp>0
\relax
			\else
  \def\pgfmathresult{}
\fi
\else
  \def\pgfmathresult{}
\fi
\else
  \def\pgfmathresult{}
\fi
},
		]
                table[ 
                x index = 3,
                y index = 4,
                ]
		{\1};
\label{plot:RSVDI1}
  \addplot[
		 color=red,
		 only marks,
		 mark=square,
     mark size = 2pt,
		 x filter/.code={\pgfplotstablegetelem{\coordindex}{[index]0}\of{\1}
\pgfmathtruncatemacro{\temp}{abs(\pgfplotsretval-5)==0? 1 : 0}
\ifnum\temp>0
\pgfplotstablegetelem{\coordindex}{[index]1}\of{\1}
\pgfmathtruncatemacro{\temp}{abs(\pgfplotsretval-3)==0? 1 : 0}
\ifnum\temp>0
\pgfplotstablegetelem{\coordindex}{[index]2}\of{\1}
\pgfmathtruncatemacro{\temp}{abs(\pgfplotsretval-1)==0? 1 : 0}
\ifnum\temp>0
\relax
			\else
  \def\pgfmathresult{}
\fi
\else
  \def\pgfmathresult{}
\fi
\else
  \def\pgfmathresult{}
\fi
},
		]
                table[ 
                x index = 3,
                y index = 4,
                ]
		{\1};
\label{plot:RSVDI2}
  \addplot[
		 color=blue,
		 only marks,
		 mark=square,
     mark size = 2pt,
		 x filter/.code={\pgfplotstablegetelem{\coordindex}{[index]0}\of{\2}
\pgfmathtruncatemacro{\temp}{abs(\pgfplotsretval-20)==0? 1 : 0}
\ifnum\temp>0
\pgfplotstablegetelem{\coordindex}{[index]1}\of{\2}
\pgfmathtruncatemacro{\temp}{abs(\pgfplotsretval-3)==0? 1 : 0}
\ifnum\temp>0
\pgfplotstablegetelem{\coordindex}{[index]2}\of{\2}
\pgfmathtruncatemacro{\temp}{abs(\pgfplotsretval-1)==0? 1 : 0}
\ifnum\temp>0
\relax
			\else
  \def\pgfmathresult{}
\fi
\else
  \def\pgfmathresult{}
\fi
\else
  \def\pgfmathresult{}
\fi
},
		]
                table[ 
                x index = 3,
                y index = 4,
                ]
		{\2};
\label{plot:COLRSVDI1}
  \addplot[
		 color=cyan,
		 only marks,
		 mark=square,
     mark size = 2pt,
		 x filter/.code={\pgfplotstablegetelem{\coordindex}{[index]0}\of{\2}
\pgfmathtruncatemacro{\temp}{abs(\pgfplotsretval-5)==0? 1 : 0}
\ifnum\temp>0
\pgfplotstablegetelem{\coordindex}{[index]1}\of{\2}
\pgfmathtruncatemacro{\temp}{abs(\pgfplotsretval-3)==0? 1 : 0}
\ifnum\temp>0
\pgfplotstablegetelem{\coordindex}{[index]2}\of{\2}
\pgfmathtruncatemacro{\temp}{abs(\pgfplotsretval-1)==0? 1 : 0}
\ifnum\temp>0
\relax
			\else
  \def\pgfmathresult{}
\fi
\else
  \def\pgfmathresult{}
\fi
\else
  \def\pgfmathresult{}
\fi
},
		]
                table[ 
                x index = 3,
                y index = 4,
                ]
		{\2};
\label{plot:COLRSVDI2}
 \coordinate (top) at (rel axis cs:0,1);
  \nextgroupplot[
  title = {sampling algorithms}, 
  yticklabels = {,,,},
  ]
  \addplot[
		 color=red,
		 only marks,
		 mark=triangle,
     mark size = 2pt,
		 x filter/.code={\pgfplotstablegetelem{\coordindex}{[index]0}\of{\4}
\pgfmathtruncatemacro{\temp}{abs(\pgfplotsretval-5)==0? 1 : 0}
\ifnum\temp>0
\pgfplotstablegetelem{\coordindex}{[index]1}\of{\4}
\pgfmathtruncatemacro{\temp}{abs(\pgfplotsretval-0.6)==0? 1 : 0}
\ifnum\temp>0
\pgfplotstablegetelem{\coordindex}{[index]2}\of{\4}
\pgfmathtruncatemacro{\temp}{abs(\pgfplotsretval-0.6)==0? 1 : 0}
\ifnum\temp>0
\relax
			\else
  \def\pgfmathresult{}
\fi
\else
  \def\pgfmathresult{}
\fi
\else
  \def\pgfmathresult{}
\fi
},
		]
                table[ 
                x index = 3,
                y index = 4,
                ]
		{\4};
\label{plot:ICSI2}
  \addplot[
		 color=red,
		 only marks,
     mark=o,
     mark size = 2pt,
		 x filter/.code={\pgfplotstablegetelem{\coordindex}{[index]0}\of{\5}
\pgfmathtruncatemacro{\temp}{abs(\pgfplotsretval-5)==0? 1 : 0}
\ifnum\temp>0
\pgfplotstablegetelem{\coordindex}{[index]1}\of{\5}
\pgfmathtruncatemacro{\temp}{abs(\pgfplotsretval-0.6)==0? 1 : 0}
\ifnum\temp>0
\pgfplotstablegetelem{\coordindex}{[index]2}\of{\5}
\pgfmathtruncatemacro{\temp}{abs(\pgfplotsretval-0.6)==0? 1 : 0}
\ifnum\temp>0
\relax
			\else
  \def\pgfmathresult{}
\fi
\else
  \def\pgfmathresult{}
\fi
\else
  \def\pgfmathresult{}
\fi
},
		]
                table[ 
                x index = 3,
                y index = 4,
                ]
		{\5};
\label{plot:ICRSI2}
  \addplot[
		 color=cyan,
		 only marks,
		 mark=triangle,
     mark size = 2pt,
		 x filter/.code={\pgfplotstablegetelem{\coordindex}{[index]0}\of{\4}
\pgfmathtruncatemacro{\temp}{abs(\pgfplotsretval-20)==0? 1 : 0}
\ifnum\temp>0
\pgfplotstablegetelem{\coordindex}{[index]1}\of{\4}
\pgfmathtruncatemacro{\temp}{abs(\pgfplotsretval-0.6)==0? 1 : 0}
\ifnum\temp>0
\pgfplotstablegetelem{\coordindex}{[index]2}\of{\4}
\pgfmathtruncatemacro{\temp}{abs(\pgfplotsretval-0.6)==0? 1 : 0}
\ifnum\temp>0
\relax
			\else
  \def\pgfmathresult{}
\fi
\else
  \def\pgfmathresult{}
\fi
\else
  \def\pgfmathresult{}
\fi
},
		]
                table[ 
                x index = 3,
                y index = 4,
                ]
		{\4};
\label{plot:ICSI1}
  \addplot[
		 color=cyan,
		 only marks,
     mark=o,
     mark size = 2pt,
		 x filter/.code={\pgfplotstablegetelem{\coordindex}{[index]0}\of{\5}
\pgfmathtruncatemacro{\temp}{abs(\pgfplotsretval-20)==0? 1 : 0}
\ifnum\temp>0
\pgfplotstablegetelem{\coordindex}{[index]1}\of{\5}
\pgfmathtruncatemacro{\temp}{abs(\pgfplotsretval-0.6)==0? 1 : 0}
\ifnum\temp>0
\pgfplotstablegetelem{\coordindex}{[index]2}\of{\5}
\pgfmathtruncatemacro{\temp}{abs(\pgfplotsretval-0.6)==0? 1 : 0}
\ifnum\temp>0
\relax
			\else
  \def\pgfmathresult{}
\fi
\else
  \def\pgfmathresult{}
\fi
\else
  \def\pgfmathresult{}
\fi
},
		]
                table[ 
                x index = 3,
                y index = 4,
                ]
		{\5};
\label{plot:ICRSI1}
 \coordinate (bot) at (rel axis cs:1,0);
\end{groupplot}
\path (top|-current bounding box.north)--
        coordinate(legendpos)
        (bot|-current bounding box.north);
  \matrix[
      matrix of nodes,
      anchor=south,
      draw,
      inner sep=0.2em,
	row 1/.style ={anchor=west},
	row 2/.style ={anchor=west},
    ]at(
legendpos)
    { 
\ref{plot:RSVDI1}& RSVD I1&[2pt]
\ref{plot:RSVDI2}& RSVD I2&[2pt]
\ref{plot:COLRSVDI1}& COLRSVD I1&[2pt]
\ref{plot:COLRSVDI2}& COLRSVD I2&[2pt]\\
\ref{plot:ICSI2}& ICS-UNF I2&[2pt]
\ref{plot:ICRSI2}& ICRS-UNF I2&[2pt]
\ref{plot:ICSI1}& ICS-UNF I1&[2pt]
\ref{plot:ICRSI1}& ICRS-UNF I1&[2pt]
\\};
\end{tikzpicture}}
  \caption{First $k$ mode angles, $\theta_i$, computed by projection algorithms with 3 subspace iterations and
iterative sampling algorithms using sampling strategy UNF.
  }
  \label{fig:angPOD_PRSVD}
\end{figure}

\subsection{Incremental-iterative algorithm for matrices that do not fit in the RAM}
As we already specified before, the mean-centered data matrix for Yale faces dataset requires 40GB memory and does not fit in a system that has 32GB RAM.
There are two possible ways to compute the POD modes of this dataset in this system:
\begin{enumerate}
	\item Run the iterative algorithm (Algorithm~\ref{alg:ICS}). In the first pass, data is read in chunks that fit in the memory to compute the sampling 
		probabilities. In the next pass, the data is read again to form the sampled matrix.
 Assuming we have enough memory for the sampled matrix and its left singular vectors in the RAM,
		$\tilde{U}$ and $\tilde{\Sigma}$ of the matrix are computed.
		Subsequent iterations may require one or both passes on data depending on the sampling strategies used.
		The iteration is continued until Algorithm~\ref{alg:ICS} converges.
              \item Do an incremental computation by running the iterative algorithm on each block of the partitioned data as follows.
                Load the block into memory and use ICS/ICRS to compute $\tilde{U}$, $\tilde{\Sigma}$ of the block. Use the MAT algorithm to get an updated $\tilde{U}$, $\tilde{\Sigma}$ of the dataset. Delete the partition and the sampled columns
                and continue the process until all blocks of data have been processed. Algorithm~\ref{alg:incr} details the steps involved.
              \end{enumerate}
             
\begin{algorithm}[htbp!]
\caption{Incremental POD computed using the iterative sampling algorithms for large matrices that do not fit in system's RAM. Input: parameters necessary for iterative algorithms and number of blocks, $t$.}
\label{alg:incr}
\begin{algorithmic}[1]
\Procedure{Incremental POD}{$A$,$t$,$k$, $c$, $w$, $r$, $\tau$, $rows$}
\State $B\gets \{\bv{a}_i : i\in 1,2,...,\lfloor n/t\rfloor\}$
\State $\tilde{U}, \tilde{\Sigma} \gets \Call{ISMA}{B,k,c,w,r,\tau,rows}$
\For{$i \in 1,...,(t-1)$}
\State $B\gets \{\bv{a}_i : i\in \lceil i n/t\rceil,...,\lfloor (i+1)n/t\rfloor\}$
\State $\hat{U}, \hat{\Sigma} \gets \Call{ISMA}{B,k,c,w,r,\tau,rows}$
\State $\hat{U},\hat{\Sigma}\gets \Call{Block Merge}{\tilde{U},\tilde{\Sigma},\hat{U},\hat{\Sigma},r}$
\State $\tilde{U}\gets \hat{U}$; $\tilde{\Sigma}\gets \hat{\Sigma}$
\EndFor
\State \textbf{return } $\tilde{\bv{u}}_i, \tilde{\sigma}_i \text{ where } i=1,2,\cdots,k$
\EndProcedure
\end{algorithmic}
\end{algorithm}
 
The first alternative requires $i+1$ passes of the matrix for ICS where $i$ is the number of iterations for convergence. 
For ICRS, it takes a maximum of $3i+1$ passes for $i$ iterations. By passes, we mean accessing the entire matrix, though it may be done in blocks.
The second alternative requires only one pass over the entire matrix irrespective of whether ICS or ICRS is used to find $\tilde{U}$, $\tilde{\Sigma}$ 
of the partition. 
Clearly, the first alternative is severely bottle-necked by disk access times. 
For YF, we used the second option and used an incremental computation after partitioning the matrix into 4 blocks.We ran the iterative algorithms on each block using the same parameters as given for YF in Table~\ref{tab:param_iter}.

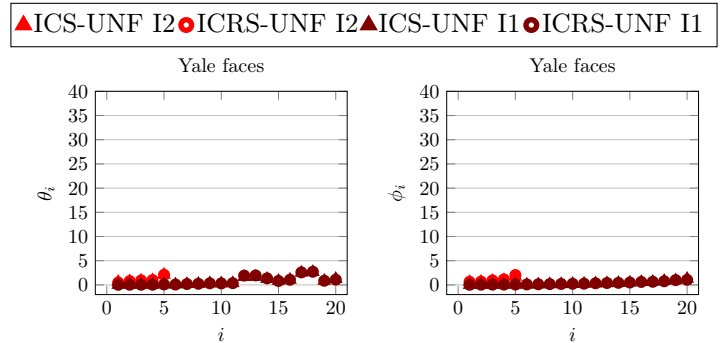
\begin{figure}[!bthp]
\begin{minipage}{\textwidth}
    \centering
\ref{angPOD_split_iter}
\end{minipage}
\begin{minipage}{\textwidth}
    \centering
      \scalebox{0.75}{\pgfplotstableread{data/angPOD_LTSVD_US_no_scaling_Yalefaces_split_into_4.dat}\4
\pgfplotstableread{data/angPOD_CTSVD_US_no_scaling_Yalefaces_split_into_4.dat}\5
\begin{tikzpicture}[every mark/.append style={line width=2pt, solid}]
\begin{axis}[
  title = {Yale faces}, 
  xlabel = $i$,
  xmax = 21.000000,
  xtick={0,5,...,20},
  ytick={0,5,...,40},
  width = 6cm,
  ymajorgrids,
   ymax = 40,
  ymin =-2,
  ylabel = $\theta_i$,
  legend to name = angPOD_split_iter,
  legend style={ legend columns = -1, cells={line width=2pt, solid}},
  legend entries = {ICS-UNF I2, ICRS-UNF I2, ICS-UNF I1, ICRS-UNF I1},
  ]
  \addplot[
		 color=red,
		 only marks,
		 mark=triangle,
     mark size = 2pt,
		 x filter/.code={\pgfplotstablegetelem{\coordindex}{[index]0}\of{\4}
\pgfmathtruncatemacro{\temp}{abs(\pgfplotsretval-5)==0? 1 : 0}
\ifnum\temp>0
\pgfplotstablegetelem{\coordindex}{[index]1}\of{\4}
\pgfmathtruncatemacro{\temp}{abs(\pgfplotsretval-0.6)==0? 1 : 0}
\ifnum\temp>0
\pgfplotstablegetelem{\coordindex}{[index]2}\of{\4}
\pgfmathtruncatemacro{\temp}{abs(\pgfplotsretval-0.6)==0? 1 : 0}
\ifnum\temp>0
\relax
			\else
  \def\pgfmathresult{}
\fi
\else
  \def\pgfmathresult{}
\fi
\else
  \def\pgfmathresult{}
\fi
},
		]
                table[ 
                x index = 3,
                y index = 4,
                ]
		{\4};
  \addplot[
		 color=red,
		 only marks,
     mark=o,
     mark size = 2pt,
		 x filter/.code={\pgfplotstablegetelem{\coordindex}{[index]0}\of{\5}
\pgfmathtruncatemacro{\temp}{abs(\pgfplotsretval-5)==0? 1 : 0}
\ifnum\temp>0
\pgfplotstablegetelem{\coordindex}{[index]1}\of{\5}
\pgfmathtruncatemacro{\temp}{abs(\pgfplotsretval-0.6)==0? 1 : 0}
\ifnum\temp>0
\pgfplotstablegetelem{\coordindex}{[index]2}\of{\5}
\pgfmathtruncatemacro{\temp}{abs(\pgfplotsretval-0.6)==0? 1 : 0}
\ifnum\temp>0
\relax
			\else
  \def\pgfmathresult{}
\fi
\else
  \def\pgfmathresult{}
\fi
\else
  \def\pgfmathresult{}
\fi
},
		]
                table[ 
                x index = 3,
                y index = 4,
                ]
		{\5};
  \addplot[
		 color=red!50!black,
		 only marks,
		 mark=triangle,
     mark size = 2pt,
		 x filter/.code={\pgfplotstablegetelem{\coordindex}{[index]0}\of{\4}
\pgfmathtruncatemacro{\temp}{abs(\pgfplotsretval-20)==0? 1 : 0}
\ifnum\temp>0
\pgfplotstablegetelem{\coordindex}{[index]1}\of{\4}
\pgfmathtruncatemacro{\temp}{abs(\pgfplotsretval-0.6)==0? 1 : 0}
\ifnum\temp>0
\pgfplotstablegetelem{\coordindex}{[index]2}\of{\4}
\pgfmathtruncatemacro{\temp}{abs(\pgfplotsretval-0.6)==0? 1 : 0}
\ifnum\temp>0
\relax
			\else
  \def\pgfmathresult{}
\fi
\else
  \def\pgfmathresult{}
\fi
\else
  \def\pgfmathresult{}
\fi
},
		]
                table[ 
                x index = 3,
                y index = 4,
                ]
		{\4};
  \addplot[
		 color=red!50!black,
		 only marks,
     mark=o,
     mark size = 2pt,
		 x filter/.code={\pgfplotstablegetelem{\coordindex}{[index]0}\of{\5}
\pgfmathtruncatemacro{\temp}{abs(\pgfplotsretval-20)==0? 1 : 0}
\ifnum\temp>0
\pgfplotstablegetelem{\coordindex}{[index]1}\of{\5}
\pgfmathtruncatemacro{\temp}{abs(\pgfplotsretval-0.6)==0? 1 : 0}
\ifnum\temp>0
\pgfplotstablegetelem{\coordindex}{[index]2}\of{\5}
\pgfmathtruncatemacro{\temp}{abs(\pgfplotsretval-0.6)==0? 1 : 0}
\ifnum\temp>0
\relax
			\else
  \def\pgfmathresult{}
\fi
\else
  \def\pgfmathresult{}
\fi
\else
  \def\pgfmathresult{}
\fi
},
		]
                table[ 
                x index = 3,
                y index = 4,
                ]
		{\5};
\end{axis}
\end{tikzpicture}}
      \scalebox{0.75}{\pgfplotstableread{data/pc_LTSVD_US_no_scaling_pc_Yalefaces_split_into_4.dat}\4
\pgfplotstableread{data/pc_CTSVD_US_no_scaling_pc_Yalefaces_split_into_4.dat}\5
\begin{tikzpicture}[every mark/.append style={line width=2pt, solid}]
\begin{axis}[
    title = {Yale faces}, 
    xlabel = $i$,
    ylabel = $\phi_i$,
    xmax = 21.000000,
  xtick={0,5,...,20},
    ytick={0,5,...,40},
    width = 6cm,
    ymajorgrids,
     ymax = 40,
    ymin =-2,
  ]
  \addplot[
		 color=red,
		 only marks,
		 mark=triangle,
     mark size = 2pt,
		 x filter/.code={\pgfplotstablegetelem{\coordindex}{[index]0}\of{\4}
\pgfmathtruncatemacro{\temp}{abs(\pgfplotsretval-5)==0? 1 : 0}
\ifnum\temp>0
\pgfplotstablegetelem{\coordindex}{[index]1}\of{\4}
\pgfmathtruncatemacro{\temp}{abs(\pgfplotsretval-0.6)==0? 1 : 0}
\ifnum\temp>0
\pgfplotstablegetelem{\coordindex}{[index]2}\of{\4}
\pgfmathtruncatemacro{\temp}{abs(\pgfplotsretval-0.6)==0? 1 : 0}
\ifnum\temp>0
\relax
			\else
  \def\pgfmathresult{}
\fi
\else
  \def\pgfmathresult{}
\fi
\else
  \def\pgfmathresult{}
\fi
},
		]
                table[ 
                x index = 3,
                y index = 4,
                ]
		{\4};
  \addplot[
		 color=red,
		 only marks,
     mark=o,
     mark size = 2pt,
		 x filter/.code={\pgfplotstablegetelem{\coordindex}{[index]0}\of{\5}
\pgfmathtruncatemacro{\temp}{abs(\pgfplotsretval-5)==0? 1 : 0}
\ifnum\temp>0
\pgfplotstablegetelem{\coordindex}{[index]1}\of{\5}
\pgfmathtruncatemacro{\temp}{abs(\pgfplotsretval-0.6)==0? 1 : 0}
\ifnum\temp>0
\pgfplotstablegetelem{\coordindex}{[index]2}\of{\5}
\pgfmathtruncatemacro{\temp}{abs(\pgfplotsretval-0.6)==0? 1 : 0}
\ifnum\temp>0
\relax
			\else
  \def\pgfmathresult{}
\fi
\else
  \def\pgfmathresult{}
\fi
\else
  \def\pgfmathresult{}
\fi
},
		]
                table[ 
                x index = 3,
                y index = 4,
                ]
		{\5};
  \addplot[
		 color=red!50!black,
		 only marks,
		 mark=triangle,
     mark size = 2pt,
		 x filter/.code={\pgfplotstablegetelem{\coordindex}{[index]0}\of{\4}
\pgfmathtruncatemacro{\temp}{abs(\pgfplotsretval-20)==0? 1 : 0}
\ifnum\temp>0
\pgfplotstablegetelem{\coordindex}{[index]1}\of{\4}
\pgfmathtruncatemacro{\temp}{abs(\pgfplotsretval-0.6)==0? 1 : 0}
\ifnum\temp>0
\pgfplotstablegetelem{\coordindex}{[index]2}\of{\4}
\pgfmathtruncatemacro{\temp}{abs(\pgfplotsretval-0.6)==0? 1 : 0}
\ifnum\temp>0
\relax
			\else
  \def\pgfmathresult{}
\fi
\else
  \def\pgfmathresult{}
\fi
\else
  \def\pgfmathresult{}
\fi
},
		]
                table[ 
                x index = 3,
                y index = 4,
                ]
		{\4};
  \addplot[
		 color=red!50!black,
		 only marks,
     mark=o,
     mark size = 2pt,
		 x filter/.code={\pgfplotstablegetelem{\coordindex}{[index]0}\of{\5}
\pgfmathtruncatemacro{\temp}{abs(\pgfplotsretval-20)==0? 1 : 0}
\ifnum\temp>0
\pgfplotstablegetelem{\coordindex}{[index]1}\of{\5}
\pgfmathtruncatemacro{\temp}{abs(\pgfplotsretval-0.6)==0? 1 : 0}
\ifnum\temp>0
\pgfplotstablegetelem{\coordindex}{[index]2}\of{\5}
\pgfmathtruncatemacro{\temp}{abs(\pgfplotsretval-0.6)==0? 1 : 0}
\ifnum\temp>0
\relax
			\else
  \def\pgfmathresult{}
\fi
\else
  \def\pgfmathresult{}
\fi
\else
  \def\pgfmathresult{}
\fi
},
		]
                table[ 
                x index = 3,
                y index = 4,
                ]
		{\5};
\end{axis}
\end{tikzpicture}}
\end{minipage}
  \caption{Left: mode angles of Yale faces computed using convergence of individual modes. Right: principal angles computed using convergence of subspace formed by modes. Results are of the incremental sampling algorithms applied on YF dataset partitioned into 4. Parameters are set as given in Table~\ref{tab:param_iter}.
  }
  \label{fig:angPOD_inc}
\end{figure}

Fig.~\ref{fig:angPOD_inc} shows that both mode angles and the principal angles are less than $5^{\circ}$.
Table~\ref{tab:time_incr} shows the wall clock time and computational time for ICS-UNF, ICRS-UNF with convergence of (a) modes and (b) subspace spanned by modes. It can be seen that for both convergence criteria,
time taken is very similar. As expected, ICRS is faster than ICS. Therefore, we favour ICRS over ICS.

\begin{table}[htbp!]
\centering
\caption{Wall clock (WT) and computational times (CT) in seconds (s) for incremental computation of POD modes for Yale faces dataset run with a single thread in a 32GB system.}
\label{tab:time_incr}
Convergence of individual modes\\
\begin{tabular}{cp{1.7cm}cccc}
\toprule
Case&sampling strategy&\multicolumn{2}{c}{ICS}&\multicolumn{2}{c}{ICRS}\\
\midrule
&&WT(s)&CT(s)&WT(s)&CT(s)\\
\cmidrule(lr){3-4}\cmidrule(lr){5-6}
$k=5$ (I2)&UNF&146.73&51.14&134.09&35.92\\
\midrule
$k=20$ (I1)&UNF&302.06&205.27&217.22&121.41\\
\bottomrule
\end{tabular}\\
Convergence of subspace spanned by modes\\
\begin{tabular}{cp{1.7cm}cccc}
\toprule
Case&sampling strategy&\multicolumn{2}{c}{ICS}&\multicolumn{2}{c}{ICRS}\\
\midrule
&&WT(s)&CT(s)&WT(s)&CT(s)\\
\cmidrule(lr){3-4}\cmidrule(lr){5-6}
$k=5$ (I2)&UNF&146.59&50.87&125.76&30.69\\
\midrule
$k=20$ (I1)&UNF&299.73&204.11&211.72&116.45\\
\bottomrule
\end{tabular}
\end{table}

\section{Conclusion}\label{sec:conclusion}
In this paper, 
we proposed an iterative algorithm to improve  the modes/subspaces spanned by the modes using a previously proposed MAT algorithm.
 Unlike the earlier methods used for multiple rounds of sampling, we estimate the POD modes in each iteration and stop sampling when no further
 improvement is obtained in either the modes or subspaces spanned by the modes. The algorithm resulted in good accuracies with
 significant improvement in runtime even if all columns are sampled when convergence is achieved. We found that using column-norm sampling
 in the first round and uniform sampling in subsequent rounds resulted in good speedups, with accuracy comparable to using norm of
 the orthogonal component of the sampled columns.

We proposed a measure using Wedin's theorem to quantify the accuracy in the computed subspaces. In most cases, the measure was much less than its upper limit,
corroborating the fact that the approximated modes were very accurate. In few cases, the measure was large due to clustering of singular values.
When the rank of approximation, $k$, was increased so that this clustering is avoided, the error measure reduced to a large extent.
In general, the accuracy of the modes are much better than indicated by the measure.
 
  For large matrices that do not fit in the RAM, we used the iterative algorithm on each partition of the data (that fit in the RAM) to approximate the 
 dominant left singular vectors and singular values of the partition 
and then used MAT operation to get the modes of the entire data. As mentioned, computing the modes using incremental sampling and MAT is advantageous 
as it requires only one pass over the data. We obtained very good accuracy using 
this method when applied on our large dataset.


\appendix
\section{Experimental evaluation of LTSVD and CTSVD}\label{app:eval}
We evaluated the column and column and row norm based sampling algorithms (referred to as LTSVD and CTSVD).
Details are contained in Algorithms~\ref{alg:LTSVD} and \ref{alg:CTSVD}.
 The parameters
to be set in the algorithms are (a) $\epsilon, \delta$: error and failure probability parameters and (b) $k$: desired number of POD modes. Based on these parameters, number of columns/rows to be sampled, $c/w$ are set according to equations~(\ref{eqn:c_LTSVD}) and (\ref{eqn:c_w_CTSVD}).
We evaluated the algorithms based on (a) runtime and (b) accuracy for various sets of parameters. 
\begin{algorithm}[htbp!]
  \small
  \caption{LTSVD algorithm. 
  Inputs: matrix $A \in \mathbb{R}^{m\times n}$, 
  vector containing sampling probabilities $\bv{p}$,
  number of columns to be sampled $c$, and 
   rank of approximation, $k$.}
  \label{alg:LTSVD}
  \begin{algorithmic}[1]
      \Procedure{LTSVD}{$A$, $\bv{p}$, $k$, $c$}
    \For {$i \gets 1,2,\dots,c$}
      \State Pick $s_i \in \{1,2,\cdots,n\} \text{ with probability  } \textbf{Pr}[s_i=\alpha]=p_{\alpha}, \alpha = 1,2,\cdots,n$
      \State $\bv{c}_i \gets \bv{a}_{s_i}/\sqrt{cp_{s_i}}$
      \EndFor
    \State $\tilde{V} \tilde{\Sigma}^2\tilde{V}^T \gets \text{SVD}(C^TC)$ \Comment \small{$\tilde{V}, \tilde{\Sigma}$ are the right singular vectors and values of $C$}
    \State $\tilde{\bv{u}}_i \gets C\tilde{\bv{v}}_i/\tilde{\sigma}_i \text{ where } i=1,2,\cdots,l \text{ and }l = \min(k, r)$; $r$: rank of $C$
    \State \textbf{return } $\tilde{\bv{u}}_i, \tilde{\sigma}_i \text{ where } i=1,2,\cdots,l$
    \EndProcedure
  \end{algorithmic}
\end{algorithm}
\begin{algorithm}[htbp!]
  \small
  \caption{CTSVD Algorithm.
  Inputs: matrix $A \in \mathbb{R}^{m\times n}$, vector containing  
	sampling probabilities for column, $\bv{p}$, 
  number of columns $c$
  and rows $w$ to be sampled,  
  and  rank of approximation, $k$}
 \label{alg:CTSVD}
  \begin{algorithmic}[1]
    \Procedure{CTSVD}{$A$, $\bv{p}$, $k$, $c$, $w$}
    \For {$i \gets 1,2,\dots,c$}
      \State Pick $s_i \in \{1,2,\cdots,n\} \text{ with probability  } \textbf{Pr}[s_i=\alpha]=p_{\alpha}, \alpha = 1,2,\cdots,n$
      \State $\bv{c}_i \gets \bv{a}_{s_i}/\sqrt{cp_{s_i}}$
      \EndFor
    \State $q_i\gets |\bv{c}^i|_2^2;q_i\gets q_i/\sum\limits_{i=1}^m q_i$
    \For {$j \gets 1,2,\dots,w$}
      \State Pick $\hat{s}_j \in \{1,2,\cdots,m\} \text{ with probability  } \textbf{Pr}[\hat{s}_j=\alpha]=q_{\alpha}, \alpha = 1,2,\cdots,m$
      \State $\bv{w}^j \gets \bv{c}^{\hat{s}_j}/\sqrt{wq_{\hat{s}_j}}$
      \EndFor
    \State $\tilde{V} \tilde{\Sigma}^2\tilde{V}^T \gets \text{SVD}(W^TW)$
    \State $\gamma \gets \epsilon/(100k)$
    \State $l\gets \min(k, \max(r:\tilde{\sigma}_r^2\ge \gamma||W||_F^2))$ \label{line:CTSVD_filter} 
    \State $\tilde{\bv{u}}_i \gets C\tilde{\bv{v}}_i/\tilde{\sigma}_i \text{ where } i=1,2,\cdots,l$ \label{line:CTSVD_find_u}
    \State \textbf{return } $\tilde{\bv{u}}_i, \tilde{\sigma}_i \text{ where } i=1,2,\cdots,l$
    \EndProcedure
  \end{algorithmic}
\end{algorithm}
\subsection{Results}
We would like the error parameters, $\epsilon$ and $\delta$, to be small ($<$ 1). 
However, since the number of columns sampled depends inversely on powers of $\epsilon$, it increases rapidly and often
exceeds the number of columns in the matrix. This limits the parameter values that can be used.
To get some representative results, we looked at results for three cases namely, 
C1: $c \approx n$ except for the YF dataset, where it is lower,  
C2: $c \approx n/2$ (lower for CYF), and C3: $ c << n$ (higher for YF).
For $V_{2D}$, Faces and CYF, we used $k = 10,5,2$ for the three cases C1, C2, and C3 respectively.
Since YF is a larger dataset, with a more gradual decay of singular values, we used  $k=20,10,5$ for 
the three cases.
For these cases, $k,\epsilon,\delta,c,w$ are shown in Tables~\ref{tab:parval}.
\begin{table}[!htbp]
  \small
  \caption{Parameter values ($k$, $\epsilon$, $\delta$) used in LTSVD and CTSVD algorithms for various datasets}
  \label{tab:parval}
  \centering
  \begin{tabular}{p{1cm}c|c|c|c}
    \toprule
    & $V_{2D}$&Faces&CYF&YF\\
    Case&$k$, $\epsilon$, $\delta$&$k$, $\epsilon$, $\delta$&$k$, $\epsilon$, $\delta$&$k$, $\epsilon$, $\delta$\\
    \midrule
    &\multicolumn{4}{c}{LTSVD}\\
    \midrule
    C1&10, 0.7, 0.45&10, 0.75, 0.8&10, 0.46, 0.44&20, 0.35, 0.35\\
    C2&5, 1, 0.1&5, 0.75, 0.75&5, 1, 0.1&10, 0.3, 0.25\\
    C3&2, 1, 0.8&2, 1, 0.8&2, 1, 0.8&5, 0.3, 0.25\\
    \midrule
    &\multicolumn{4}{c}{CTSVD}\\
    \midrule
    C1&10, 1.04, 0.94&10, 1.3, 1&10, 0.9, 0.7&20, 0.83, 0.9\\
    C2&5, 1, 0.1&5, 1, 0.1&5, 1, 1&10, 0.62, 0.9\\
    C3&2, 1, 0.8&2, 1, 0.8&2, 1, 0.8&5, 0.5, 0.9\\
        \bottomrule
  \end{tabular}
\end{table}

%

\begin{figure}[!htbp]
\hspace{-0.1in}    \scalebox{0.6}{\pgfplotstableread{data/sigma_V_2D.dat1}\1
\begin{tikzpicture}[every mark/.append style={solid}]
\begin{semilogyaxis}[ 
    scaled y ticks = false,
    y tick label style = {/pgf/number format/sci},
  title = {$V_{2D}$}, 
    xlabel = $i$,
    ylabel = $\sigma_i$,
    legend pos =  north east,
    width = 6cm,
    xmax = 20,
    ]
    \addplot[
            only marks, 
		        mark=square,
            black, 
            ]table[ x expr = \thisrowno{0}+1, y index = 1] {\1};
\end{semilogyaxis}
\end{tikzpicture}}
\hspace{-0.1in}
  \scalebox{0.6}{\pgfplotstableread{data/sigma_faces.dat1}\1
\begin{tikzpicture}[every mark/.append style={solid}]
\begin{axis}[ 
    scaled y ticks = false,
  title = {Faces dataset}, 
    xlabel = $i$,
    legend pos =  north east,
    width = 6cm,
    xmax = 20,
    scaled y ticks={base 10:-5},
    ]
    \addplot[
            only marks, 
		        mark=square,
            black, 
            ]table[ x expr = \thisrowno{0}+1, y index = 1] {\1};
\end{axis}
\end{tikzpicture}}
\hspace{-0.1in}
  \scalebox{0.6}{ \pgfplotstableread{data/sigma_croppedYalefaces.dat1}\1
\begin{tikzpicture}[every mark/.append style={solid}]
\begin{axis}[ 
    scaled y ticks = false,
  title = {cropped Yale faces}, 
    xlabel = $i$,
    legend pos =  north east,
    width = 6cm,
    xmax = 20,
    ymin=5000,
    ymax=400000,
    scaled y ticks={base 10:-5},
    ]
    \addplot[
            only marks, 
		        mark=square,
            black, 
            ]table[ x expr = \thisrowno{0}+1, y index = 1] {\1};
\end{axis}
\end{tikzpicture}}
\hspace{-0.1in}
  \scalebox{0.6}{ \pgfplotstableread{data/sigma_Yalefaces.dat1}\1
\begin{tikzpicture}[every mark/.append style={ solid}]
\begin{axis}[ 
    scaled y ticks = false,
  title = {Yale faces}, 
    xlabel = $i$,
    legend pos =  north east,
    width = 6cm,
    xmax = 20,
    ymin = 80000,
    ymax = 2000000,
    scaled y ticks={base 10:-6},
    ]
    \addplot[
            only marks, 
		        mark=square,
            black, 
            ]table[ x expr = \thisrowno{0}+1, y index = 1] {\1};
\end{axis}
\end{tikzpicture}}
  \caption{First 20 singular values of datasets considered.}
  \label{fig:sigma}
\end{figure}
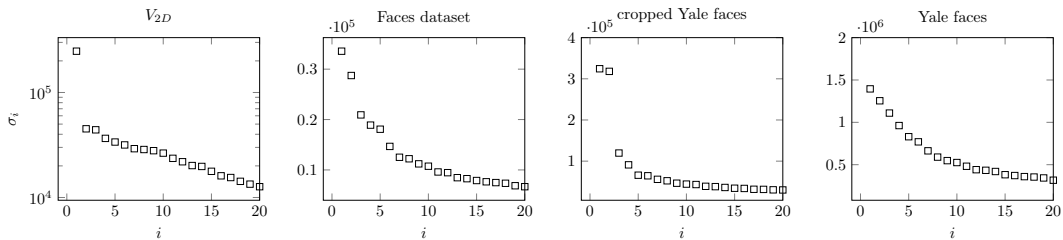

Fig.~\ref{fig:sigma} shows the first 20 singular values of the datasets. $V_{2D}$ has a dominant singular value, followed by a sharp decay.
The most gradual decay of singular values is for the YF dataset, where the first two modes capture only 30\% of the energy.
For CYF the first two singular values are very close. For all datasets, but especially for CYF,
there is a clustering of the singular values beyond the first ten values.
Note that the singular values of $V_{2D}$ is plotted on a log scale, while others are on a linear scale.

\begin{figure}[!htbp]
  \begin{minipage}{\textwidth}
\centering
  \ref{relsigDiff_dri}\\
  \end{minipage}
  \begin{minipage}{\textwidth}
\hspace{-0.4in}    \scalebox{0.75}{\pgfplotstableread{data/relsigDiff_LTSVD_drineas_V_2D.dat}\4
\pgfplotstableread{data/relsigDiff_CTSVD_drineas_V_2D.dat}\1

\begin{tikzpicture}[every mark/.append style={line width=2pt, solid}]
\begin{axis}[ 
  title = {$V_{2D}$}, 
  scaled y ticks = false,
    y tick label style = {/pgf/number format/.cd ={fixed, precision = 2}},
    xlabel = $i$,
    ylabel = $\frac{|\sigma_i-\tilde{\sigma}_i|}{\sigma_i}*100$,
    legend to name = relsigDiff_dri,
    legend style={ at={(0.1,1.05)},anchor=south,legend columns = 6, cells={line width=2pt, solid}},
    legend entries = {LTSVD C1, CTSVD C1, LTSVD C2, CTSVD C2, LTSVD C3, CTSVD C3},
    xmax = 11.000000,
    ymax = 32,
    ytick={0,5,10,20},
    xtick={2,4,...,10},
    ymajorgrids,
    width = 6cm,
    ]
   \addplot[
		 color=cyan,
		 only marks,
		 mark=triangle,
     mark size = 2pt,
		 x filter/.code={\pgfplotstablegetelem{\coordindex}{[index]0}\of{\4}
\pgfmathtruncatemacro{\temp}{abs(\pgfplotsretval-10)==0? 1 : 0}
\ifnum\temp>0
\pgfplotstablegetelem{\coordindex}{[index]1}\of{\4}
\pgfmathtruncatemacro{\temp}{abs(\pgfplotsretval-0.7)==0? 1 : 0}
\ifnum\temp>0
\pgfplotstablegetelem{\coordindex}{[index]2}\of{\4}
\pgfmathtruncatemacro{\temp}{abs(\pgfplotsretval-0.45)==0? 1 : 0}
\ifnum\temp>0
\relax
			\else
  \def\pgfmathresult{}
\fi
\else
  \def\pgfmathresult{}
\fi
\else
  \def\pgfmathresult{}
\fi
},
		]
                table[ 
                x index = 3,
                 y expr = \thisrowno{4}*100,
                ]
		{\4};
  \addplot[
		 color=cyan,
		 only marks,
     mark=o,
     mark size = 2pt,
		 x filter/.code={\pgfplotstablegetelem{\coordindex}{[index]0}\of{\1}
\pgfmathtruncatemacro{\temp}{abs(\pgfplotsretval-10)==0? 1 : 0}
\ifnum\temp>0
\pgfplotstablegetelem{\coordindex}{[index]1}\of{\1}
\pgfmathtruncatemacro{\temp}{abs(\pgfplotsretval-1.04)==0? 1 : 0}
\ifnum\temp>0
\pgfplotstablegetelem{\coordindex}{[index]2}\of{\1}
\pgfmathtruncatemacro{\temp}{abs(\pgfplotsretval-0.94)==0? 1 : 0}
\ifnum\temp>0
\relax
			\else
  \def\pgfmathresult{}
\fi
\else
  \def\pgfmathresult{}
\fi
\else
  \def\pgfmathresult{}
\fi
},
		]
                table[ 
                x index = 3,
                 y expr = \thisrowno{4}*100,
                ]
		{\1};
                \addplot[
		 color=green!50!black,
		 only marks,
		 mark=triangle,
     mark size = 1.5pt,
		 x filter/.code={\pgfplotstablegetelem{\coordindex}{[index]0}\of{\4}
\pgfmathtruncatemacro{\temp}{abs(\pgfplotsretval-5)==0? 1 : 0}
\ifnum\temp>0
\pgfplotstablegetelem{\coordindex}{[index]1}\of{\4}
\pgfmathtruncatemacro{\temp}{abs(\pgfplotsretval-1.0)==0? 1 : 0}
\ifnum\temp>0
\pgfplotstablegetelem{\coordindex}{[index]2}\of{\4}
\pgfmathtruncatemacro{\temp}{abs(\pgfplotsretval-0.1)==0? 1 : 0}
\ifnum\temp>0
\relax
			\else
  \def\pgfmathresult{}
\fi
\else
  \def\pgfmathresult{}
\fi
\else
  \def\pgfmathresult{}
\fi
},
		]
                table[ 
                x index = 3,
                 y expr = \thisrowno{4}*100,
                ]
		{\4};
                \addplot[
		 color=green!50!black,
		 only marks,
		 mark=o,
		 x filter/.code={\pgfplotstablegetelem{\coordindex}{[index]0}\of{\1}
\pgfmathtruncatemacro{\temp}{abs(\pgfplotsretval-5)==0? 1 : 0}
\ifnum\temp>0
\pgfplotstablegetelem{\coordindex}{[index]1}\of{\1}
\pgfmathtruncatemacro{\temp}{abs(\pgfplotsretval-1.0)==0? 1 : 0}
\ifnum\temp>0
\pgfplotstablegetelem{\coordindex}{[index]2}\of{\1}
\pgfmathtruncatemacro{\temp}{abs(\pgfplotsretval-0.1)==0? 1 : 0}
\ifnum\temp>0
\relax
			\else
  \def\pgfmathresult{}
\fi
\else
  \def\pgfmathresult{}
\fi
\else
  \def\pgfmathresult{}
\fi
},
		]
                table[ 
                x index = 3,
                 y expr = \thisrowno{4}*100,
                ]
		{\1};
   \addplot[
		 color=red,
		 only marks,
		 mark=triangle,
     mark size = 2pt,
		 x filter/.code={\pgfplotstablegetelem{\coordindex}{[index]0}\of{\4}
\pgfmathtruncatemacro{\temp}{abs(\pgfplotsretval-2)==0? 1 : 0}
\ifnum\temp>0
\pgfplotstablegetelem{\coordindex}{[index]1}\of{\4}
\pgfmathtruncatemacro{\temp}{abs(\pgfplotsretval-1.0)==0? 1 : 0}
\ifnum\temp>0
\pgfplotstablegetelem{\coordindex}{[index]2}\of{\4}
\pgfmathtruncatemacro{\temp}{abs(\pgfplotsretval-0.8)==0? 1 : 0}
\ifnum\temp>0
\relax
			\else
  \def\pgfmathresult{}
\fi
\else
  \def\pgfmathresult{}
\fi
\else
  \def\pgfmathresult{}
\fi
},
		]
                table[ 
                x index = 3,
                 y expr = \thisrowno{4}*100,
                ]
		{\4};
  \addplot[
		 color=red,
		 only marks,
     mark=o,
     mark size = 2pt,
		 x filter/.code={\pgfplotstablegetelem{\coordindex}{[index]0}\of{\1}
\pgfmathtruncatemacro{\temp}{abs(\pgfplotsretval-2)==0? 1 : 0}
\ifnum\temp>0
\pgfplotstablegetelem{\coordindex}{[index]1}\of{\1}
\pgfmathtruncatemacro{\temp}{abs(\pgfplotsretval-1.0)==0? 1 : 0}
\ifnum\temp>0
\pgfplotstablegetelem{\coordindex}{[index]2}\of{\1}
\pgfmathtruncatemacro{\temp}{abs(\pgfplotsretval-0.8)==0? 1 : 0}
\ifnum\temp>0
\relax
			\else
  \def\pgfmathresult{}
\fi
\else
  \def\pgfmathresult{}
\fi
\else
  \def\pgfmathresult{}
\fi
},
		]
                table[ 
                x index = 3,
                 y expr = \thisrowno{4}*100,
                ]
		{\1};
    \end{axis}
    \end{tikzpicture}}
\hspace{-0.2in}    \scalebox{0.75}{\pgfplotstableread{data/relsigDiff_LTSVD_drineas_faces.dat}\4
\pgfplotstableread{data/relsigDiff_CTSVD_drineas_faces.dat}\1

\begin{tikzpicture}[every mark/.append style={line width=2pt, solid}]
\begin{axis}[ 
  title = {Faces dataset}, 
  scaled y ticks = false,
    y tick label style = {/pgf/number format/.cd = {fixed, precision = 2}},
    xlabel = $i$,
     xmax = 11.000000,
     ymax=32,
    xtick={2,4,...,10},
    ytick={0,5, 10,20},
    yticklabels={,,},
    ymajorgrids,
    width = 6cm,
    ]
   \addplot[
		 color=cyan,
		 only marks,
		 mark=triangle,
     mark size = 2pt,
		 x filter/.code={\pgfplotstablegetelem{\coordindex}{[index]0}\of{\4}
\pgfmathtruncatemacro{\temp}{abs(\pgfplotsretval-10)==0? 1 : 0}
\ifnum\temp>0
\pgfplotstablegetelem{\coordindex}{[index]1}\of{\4}
\pgfmathtruncatemacro{\temp}{abs(\pgfplotsretval-0.75)==0? 1 : 0}
\ifnum\temp>0
\pgfplotstablegetelem{\coordindex}{[index]2}\of{\4}
\pgfmathtruncatemacro{\temp}{abs(\pgfplotsretval-0.8)==0? 1 : 0}
\ifnum\temp>0
\relax
			\else
  \def\pgfmathresult{}
\fi
\else
  \def\pgfmathresult{}
\fi
\else
  \def\pgfmathresult{}
\fi
},
		]
                table[ 
                x index = 3,
                 y expr = \thisrowno{4}*100,
                ]
		{\4};
  \addplot[
		 color=cyan,
		 only marks,
     mark=o,
     mark size = 2pt,
		 x filter/.code={\pgfplotstablegetelem{\coordindex}{[index]0}\of{\1}
\pgfmathtruncatemacro{\temp}{abs(\pgfplotsretval-10)==0? 1 : 0}
\ifnum\temp>0
\pgfplotstablegetelem{\coordindex}{[index]1}\of{\1}
\pgfmathtruncatemacro{\temp}{abs(\pgfplotsretval-1.3)==0? 1 : 0}
\ifnum\temp>0
\pgfplotstablegetelem{\coordindex}{[index]2}\of{\1}
\pgfmathtruncatemacro{\temp}{abs(\pgfplotsretval-1)==0? 1 : 0}
\ifnum\temp>0
\relax
			\else
  \def\pgfmathresult{}
\fi
\else
  \def\pgfmathresult{}
\fi
\else
  \def\pgfmathresult{}
\fi
},
		]
                table[ 
                x index = 3,
                 y expr = \thisrowno{4}*100,
                ]
		{\1};
                \addplot[
		 color=green!50!black,
		 only marks,
		 mark=triangle,
     mark size = 1.5pt,
		 x filter/.code={\pgfplotstablegetelem{\coordindex}{[index]0}\of{\4}
\pgfmathtruncatemacro{\temp}{abs(\pgfplotsretval-5)==0? 1 : 0}
\ifnum\temp>0
\pgfplotstablegetelem{\coordindex}{[index]1}\of{\4}
\pgfmathtruncatemacro{\temp}{abs(\pgfplotsretval-0.75)==0? 1 : 0}
\ifnum\temp>0
\pgfplotstablegetelem{\coordindex}{[index]2}\of{\4}
\pgfmathtruncatemacro{\temp}{abs(\pgfplotsretval-0.75)==0? 1 : 0}
\ifnum\temp>0
\relax
			\else
  \def\pgfmathresult{}
\fi
\else
  \def\pgfmathresult{}
\fi
\else
  \def\pgfmathresult{}
\fi
},
		]
                table[ 
                x index = 3,
                 y expr = \thisrowno{4}*100,
                ]
		{\4};
                \addplot[
                color=green!50!black,
		 only marks,
		 mark=o,
		 x filter/.code={\pgfplotstablegetelem{\coordindex}{[index]0}\of{\1}
\pgfmathtruncatemacro{\temp}{abs(\pgfplotsretval-5)==0? 1 : 0}
\ifnum\temp>0
\pgfplotstablegetelem{\coordindex}{[index]1}\of{\1}
\pgfmathtruncatemacro{\temp}{abs(\pgfplotsretval-1.0)==0? 1 : 0}
\ifnum\temp>0
\pgfplotstablegetelem{\coordindex}{[index]2}\of{\1}
\pgfmathtruncatemacro{\temp}{abs(\pgfplotsretval-1)==0? 1 : 0}
\ifnum\temp>0
\relax
			\else
  \def\pgfmathresult{}
\fi
\else
  \def\pgfmathresult{}
\fi
\else
  \def\pgfmathresult{}
\fi
},
		]
                table[ 
                x index = 3,
                 y expr = \thisrowno{4}*100,
                ]
		{\1};
   \addplot[
		 color=red,
		 only marks,
		 mark=triangle,
     mark size = 2pt,
		 x filter/.code={\pgfplotstablegetelem{\coordindex}{[index]0}\of{\4}
\pgfmathtruncatemacro{\temp}{abs(\pgfplotsretval-2)==0? 1 : 0}
\ifnum\temp>0
\pgfplotstablegetelem{\coordindex}{[index]1}\of{\4}
\pgfmathtruncatemacro{\temp}{abs(\pgfplotsretval-1.0)==0? 1 : 0}
\ifnum\temp>0
\pgfplotstablegetelem{\coordindex}{[index]2}\of{\4}
\pgfmathtruncatemacro{\temp}{abs(\pgfplotsretval-0.8)==0? 1 : 0}
\ifnum\temp>0
\relax
			\else
  \def\pgfmathresult{}
\fi
\else
  \def\pgfmathresult{}
\fi
\else
  \def\pgfmathresult{}
\fi
},
		]
                table[ 
                x index = 3,
                 y expr = \thisrowno{4}*100,
                ]
		{\4};
  \addplot[
		 color=red,
		 only marks,
     mark=o,
     mark size = 2pt,
		 x filter/.code={\pgfplotstablegetelem{\coordindex}{[index]0}\of{\1}
\pgfmathtruncatemacro{\temp}{abs(\pgfplotsretval-2)==0? 1 : 0}
\ifnum\temp>0
\pgfplotstablegetelem{\coordindex}{[index]1}\of{\1}
\pgfmathtruncatemacro{\temp}{abs(\pgfplotsretval-1.0)==0? 1 : 0}
\ifnum\temp>0
\pgfplotstablegetelem{\coordindex}{[index]2}\of{\1}
\pgfmathtruncatemacro{\temp}{abs(\pgfplotsretval-0.8)==0? 1 : 0}
\ifnum\temp>0
\relax
			\else
  \def\pgfmathresult{}
\fi
\else
  \def\pgfmathresult{}
\fi
\else
  \def\pgfmathresult{}
\fi
},
		]
                table[ 
                x index = 3,
                 y expr = \thisrowno{4}*100,
                ]
		{\1};

    \end{axis}
    \end{tikzpicture}}
\hspace{-0.2in}     \scalebox{0.75}{ \pgfplotstableread{data/relsigDiff_LTSVD_drineas_croppedYalefaces.dat}\4
\pgfplotstableread{data/relsigDiff_CTSVD_drineas_croppedYalefaces.dat}\1

\begin{tikzpicture}[every mark/.append style={line width=2pt, solid}]
\begin{axis}[ 
    scaled y ticks = false,
    y tick label style = {/pgf/number format/.cd = {fixed, precision = 2}},
  title = {cropped Yale faces}, 
    xlabel = $i$,
    xmax = 11.000000,
    ymax=32,
    xtick={2,4,...,10},
    ytick={0,5,10,20},
    yticklabels={,,},
    ymajorgrids,
    width = 6cm,
    ]
   \addplot[
		 color=cyan,
		 only marks,
		 mark=triangle,
     mark size = 2pt,
		 x filter/.code={\pgfplotstablegetelem{\coordindex}{[index]0}\of{\4}
\pgfmathtruncatemacro{\temp}{abs(\pgfplotsretval-10)==0? 1 : 0}
\ifnum\temp>0
\pgfplotstablegetelem{\coordindex}{[index]1}\of{\4}
\pgfmathtruncatemacro{\temp}{abs(\pgfplotsretval-0.46)==0? 1 : 0}
\ifnum\temp>0
\pgfplotstablegetelem{\coordindex}{[index]2}\of{\4}
\pgfmathtruncatemacro{\temp}{abs(\pgfplotsretval-0.44)==0? 1 : 0}
\ifnum\temp>0
\relax
			\else
  \def\pgfmathresult{}
\fi
\else
  \def\pgfmathresult{}
\fi
\else
  \def\pgfmathresult{}
\fi
},
		]
                table[ 
                x index = 3,
                 y expr = \thisrowno{4}*100,
                ]
		{\4};
  \addplot[
		 color=cyan,
		 only marks,
     mark=o,
     mark size = 2pt,
		 x filter/.code={\pgfplotstablegetelem{\coordindex}{[index]0}\of{\1}
\pgfmathtruncatemacro{\temp}{abs(\pgfplotsretval-10)==0? 1 : 0}
\ifnum\temp>0
\pgfplotstablegetelem{\coordindex}{[index]1}\of{\1}
\pgfmathtruncatemacro{\temp}{abs(\pgfplotsretval-0.9)==0? 1 : 0}
\ifnum\temp>0
\pgfplotstablegetelem{\coordindex}{[index]2}\of{\1}
\pgfmathtruncatemacro{\temp}{abs(\pgfplotsretval-0.7)==0? 1 : 0}
\ifnum\temp>0
\relax
			\else
  \def\pgfmathresult{}
\fi
\else
  \def\pgfmathresult{}
\fi
\else
  \def\pgfmathresult{}
\fi
},
		]
                table[ 
                x index = 3,
                 y expr = \thisrowno{4}*100,
                ]
		{\1};
                \addplot[
		 color=green!50!black,
		 only marks,
		 mark=triangle,
     mark size = 1.5pt,
		 x filter/.code={\pgfplotstablegetelem{\coordindex}{[index]0}\of{\4}
\pgfmathtruncatemacro{\temp}{abs(\pgfplotsretval-5)==0? 1 : 0}
\ifnum\temp>0
\pgfplotstablegetelem{\coordindex}{[index]1}\of{\4}
\pgfmathtruncatemacro{\temp}{abs(\pgfplotsretval-1.0)==0? 1 : 0}
\ifnum\temp>0
\pgfplotstablegetelem{\coordindex}{[index]2}\of{\4}
\pgfmathtruncatemacro{\temp}{abs(\pgfplotsretval-0.1)==0? 1 : 0}
\ifnum\temp>0
\relax
			\else
  \def\pgfmathresult{}
\fi
\else
  \def\pgfmathresult{}
\fi
\else
  \def\pgfmathresult{}
\fi
},
		]
                table[ 
                x index = 3,
                 y expr = \thisrowno{4}*100,
                ]
		{\4};
                \addplot[
		 color=green!50!black,
		 only marks,
		 mark=o,
		 x filter/.code={\pgfplotstablegetelem{\coordindex}{[index]0}\of{\1}
\pgfmathtruncatemacro{\temp}{abs(\pgfplotsretval-5)==0? 1 : 0}
\ifnum\temp>0
\pgfplotstablegetelem{\coordindex}{[index]1}\of{\1}
\pgfmathtruncatemacro{\temp}{abs(\pgfplotsretval-1.0)==0? 1 : 0}
\ifnum\temp>0
\pgfplotstablegetelem{\coordindex}{[index]2}\of{\1}
\pgfmathtruncatemacro{\temp}{abs(\pgfplotsretval-0.1)==0? 1 : 0}
\ifnum\temp>0
\relax
			\else
  \def\pgfmathresult{}
\fi
\else
  \def\pgfmathresult{}
\fi
\else
  \def\pgfmathresult{}
\fi
},
		]
                table[ 
                x index = 3,
                 y expr = \thisrowno{4}*100,
                ]
		{\1};
   \addplot[
		 color=red,
		 only marks,
		 mark=triangle,
     mark size = 2pt,
		 x filter/.code={\pgfplotstablegetelem{\coordindex}{[index]0}\of{\4}
\pgfmathtruncatemacro{\temp}{abs(\pgfplotsretval-2)==0? 1 : 0}
\ifnum\temp>0
\pgfplotstablegetelem{\coordindex}{[index]1}\of{\4}
\pgfmathtruncatemacro{\temp}{abs(\pgfplotsretval-1.0)==0? 1 : 0}
\ifnum\temp>0
\pgfplotstablegetelem{\coordindex}{[index]2}\of{\4}
\pgfmathtruncatemacro{\temp}{abs(\pgfplotsretval-0.8)==0? 1 : 0}
\ifnum\temp>0
\relax
			\else
  \def\pgfmathresult{}
\fi
\else
  \def\pgfmathresult{}
\fi
\else
  \def\pgfmathresult{}
\fi
},
		]
                table[ 
                x index = 3,
                 y expr = \thisrowno{4}*100,
                ]
		{\4};
  \addplot[
		 color=red,
		 only marks,
     mark=o,
     mark size = 2pt,
		 x filter/.code={\pgfplotstablegetelem{\coordindex}{[index]0}\of{\1}
\pgfmathtruncatemacro{\temp}{abs(\pgfplotsretval-2)==0? 1 : 0}
\ifnum\temp>0
\pgfplotstablegetelem{\coordindex}{[index]1}\of{\1}
\pgfmathtruncatemacro{\temp}{abs(\pgfplotsretval-1.0)==0? 1 : 0}
\ifnum\temp>0
\pgfplotstablegetelem{\coordindex}{[index]2}\of{\1}
\pgfmathtruncatemacro{\temp}{abs(\pgfplotsretval-0.8)==0? 1 : 0}
\ifnum\temp>0
\relax
			\else
  \def\pgfmathresult{}
\fi
\else
  \def\pgfmathresult{}
\fi
\else
  \def\pgfmathresult{}
\fi
},
		]
                table[ 
                x index = 3,
                 y expr = \thisrowno{4}*100,
                ]
		{\1};

\end{axis}
\end{tikzpicture}}
\hspace{-0.22in}     \scalebox{0.75}{ \pgfplotstableread{data/relsigDiff_LTSVD_drineas_ext_Yalefaces_wo_Ambient.dat}\4
\pgfplotstableread{data/relsigDiff_CTSVD_drineas_ext_Yalefaces_wo_Ambient.dat}\1

\begin{tikzpicture}[every mark/.append style={line width=2pt, solid}]
\begin{axis}[
    title = {Yale faces}, 
    xlabel = $i$,
    xmax = 21.000000,
    ymax=32,
    ytick={0, 5, 10, 20},
    yticklabels={,,},
    width = 6cm,
    ymajorgrids,
  ]
  \addplot[
		 color=cyan,
		 only marks,
		 mark=triangle,
     mark size = 2pt,
		 x filter/.code={\pgfplotstablegetelem{\coordindex}{[index]0}\of{\4}
\pgfmathtruncatemacro{\temp}{abs(\pgfplotsretval-20)==0? 1 : 0}
\ifnum\temp>0
\pgfplotstablegetelem{\coordindex}{[index]1}\of{\4}
\pgfmathtruncatemacro{\temp}{abs(\pgfplotsretval-0.35)==0? 1 : 0}
\ifnum\temp>0
\pgfplotstablegetelem{\coordindex}{[index]2}\of{\4}
\pgfmathtruncatemacro{\temp}{abs(\pgfplotsretval-0.35)==0? 1 : 0}
\ifnum\temp>0
\relax
			\else
  \def\pgfmathresult{}
\fi
\else
  \def\pgfmathresult{}
\fi
\else
  \def\pgfmathresult{}
\fi
},
		]
                table[ 
                x index = 3,
                 y expr = \thisrowno{4}*100,
                ]
		{\4};
  \addplot[
		 color=cyan,
		 only marks,
     mark=o,
     mark size = 2pt,
		 x filter/.code={\pgfplotstablegetelem{\coordindex}{[index]0}\of{\1}
\pgfmathtruncatemacro{\temp}{abs(\pgfplotsretval-20)==0? 1 : 0}
\ifnum\temp>0
\pgfplotstablegetelem{\coordindex}{[index]1}\of{\1}
\pgfmathtruncatemacro{\temp}{abs(\pgfplotsretval-0.83)==0? 1 : 0}
\ifnum\temp>0
\pgfplotstablegetelem{\coordindex}{[index]2}\of{\1}
\pgfmathtruncatemacro{\temp}{abs(\pgfplotsretval-0.9)==0? 1 : 0}
\ifnum\temp>0
\relax
			\else
  \def\pgfmathresult{}
\fi
\else
  \def\pgfmathresult{}
\fi
\else
  \def\pgfmathresult{}
\fi
},
		]
                table[ 
                x index = 3,
                 y expr = \thisrowno{4}*100,
                ]
		{\1};
                \addplot[
		 color=green!50!black,
		 only marks,
		 mark=triangle,
     mark size = 1.5pt,
		 x filter/.code={\pgfplotstablegetelem{\coordindex}{[index]0}\of{\4}
\pgfmathtruncatemacro{\temp}{abs(\pgfplotsretval-10)==0? 1 : 0}
\ifnum\temp>0
\pgfplotstablegetelem{\coordindex}{[index]1}\of{\4}
\pgfmathtruncatemacro{\temp}{abs(\pgfplotsretval-0.3)==0? 1 : 0}
\ifnum\temp>0
\pgfplotstablegetelem{\coordindex}{[index]2}\of{\4}
\pgfmathtruncatemacro{\temp}{abs(\pgfplotsretval-0.25)==0? 1 : 0}
\ifnum\temp>0
\relax
			\else
  \def\pgfmathresult{}
\fi
\else
  \def\pgfmathresult{}
\fi
\else
  \def\pgfmathresult{}
\fi
},
		]
                table[ 
                x index = 3,
                 y expr = \thisrowno{4}*100,
                ]
		{\4};
                \addplot[
		 color=green!50!black,
		 only marks,
		 mark=o,
		 x filter/.code={\pgfplotstablegetelem{\coordindex}{[index]0}\of{\1}
\pgfmathtruncatemacro{\temp}{abs(\pgfplotsretval-10)==0? 1 : 0}
\ifnum\temp>0
\pgfplotstablegetelem{\coordindex}{[index]1}\of{\1}
\pgfmathtruncatemacro{\temp}{abs(\pgfplotsretval-0.62)==0? 1 : 0}
\ifnum\temp>0
\pgfplotstablegetelem{\coordindex}{[index]2}\of{\1}
\pgfmathtruncatemacro{\temp}{abs(\pgfplotsretval-0.9)==0? 1 : 0}
\ifnum\temp>0
\relax
			\else
  \def\pgfmathresult{}
\fi
\else
  \def\pgfmathresult{}
\fi
\else
  \def\pgfmathresult{}
\fi
},
		]
                table[ 
                x index = 3,
                 y expr = \thisrowno{4}*100,
                ]
		{\1};
   \addplot[
		 color=red,
		 only marks,
		 mark=triangle,
     mark size = 2pt,
		 x filter/.code={\pgfplotstablegetelem{\coordindex}{[index]0}\of{\4}
\pgfmathtruncatemacro{\temp}{abs(\pgfplotsretval-5)==0? 1 : 0}
\ifnum\temp>0
\pgfplotstablegetelem{\coordindex}{[index]1}\of{\4}
\pgfmathtruncatemacro{\temp}{abs(\pgfplotsretval-0.3)==0? 1 : 0}
\ifnum\temp>0
\pgfplotstablegetelem{\coordindex}{[index]2}\of{\4}
\pgfmathtruncatemacro{\temp}{abs(\pgfplotsretval-0.25)==0? 1 : 0}
\ifnum\temp>0
\relax
			\else
  \def\pgfmathresult{}
\fi
\else
  \def\pgfmathresult{}
\fi
\else
  \def\pgfmathresult{}
\fi
},
		]
                table[ 
                x index = 3,
                 y expr = \thisrowno{4}*100,
                ]
		{\4};
  \addplot[
		 color=red,
		 only marks,
     mark=o,
     mark size = 2pt,
		 x filter/.code={\pgfplotstablegetelem{\coordindex}{[index]0}\of{\1}
\pgfmathtruncatemacro{\temp}{abs(\pgfplotsretval-5)==0? 1 : 0}
\ifnum\temp>0
\pgfplotstablegetelem{\coordindex}{[index]1}\of{\1}
\pgfmathtruncatemacro{\temp}{abs(\pgfplotsretval-0.5)==0? 1 : 0}
\ifnum\temp>0
\pgfplotstablegetelem{\coordindex}{[index]2}\of{\1}
\pgfmathtruncatemacro{\temp}{abs(\pgfplotsretval-0.9)==0? 1 : 0}
\ifnum\temp>0
\relax
			\else
  \def\pgfmathresult{}
\fi
\else
  \def\pgfmathresult{}
\fi
\else
  \def\pgfmathresult{}
\fi
},
		]
                table[ 
                x index = 3,
                 y expr = \thisrowno{4}*100,
                ]
		{\1};
                
    \end{axis}
    \end{tikzpicture}}
  \end{minipage}
  \caption{Percentage error in the  singular values of different datasets for different parameter values, as given in Table~\ref{tab:parval}.}
  \label{fig:relsigDiff_LTSVD_CTSVD}
\end{figure}
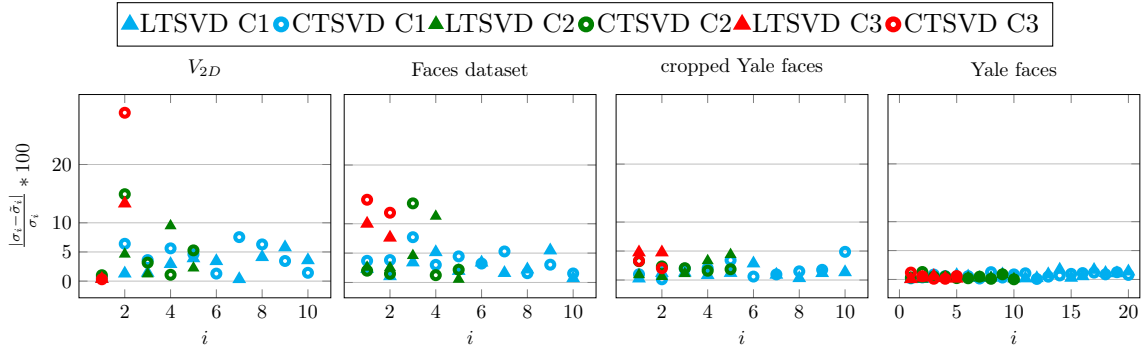

We found that the  error in approximating $A_k$ 
is well below the theoretical error bound for the respective algorithms in all cases.  
Fig.~\ref{fig:relsigDiff_LTSVD_CTSVD} shows the percentage
error in the singular values. \footnote{In the case of CTSVD C1, we could not obtain 10 values with the filter (see line \ref{line:CTSVD_filter} of Algorithm~\ref{alg:LTSVD}) and the values plotted are obtained without the filter.} It is seen that the error is less than 5\% for CYF and YF. The two dominant singular values are captured
accurately in CYF even when the error parameter values are close to one.
The error is slightly larger for $V_{2D}$ and Faces dataset, but is still less
than 20\%, except for CTSVD (C3) for the $V_{2D}$ dataset.

\begin{figure}[!htbp]
\centering
  \begin{tabular}{cccccc}
    \addlinespace[0.5pt]
     \multicolumn{6}{c}{$V_{2D}$}\\
    LTSVD C1&
     \hspace{-0.15in}   SVD&
     \hspace{-0.15in}   CTSVD C1&
     \hspace{-0.15in}   LTSVD C2&
     \hspace{-0.15in}   SVD&
     \hspace{-0.15in}   CTSVD C2\\
    \includegraphics[keepaspectratio, scale= 0.4]{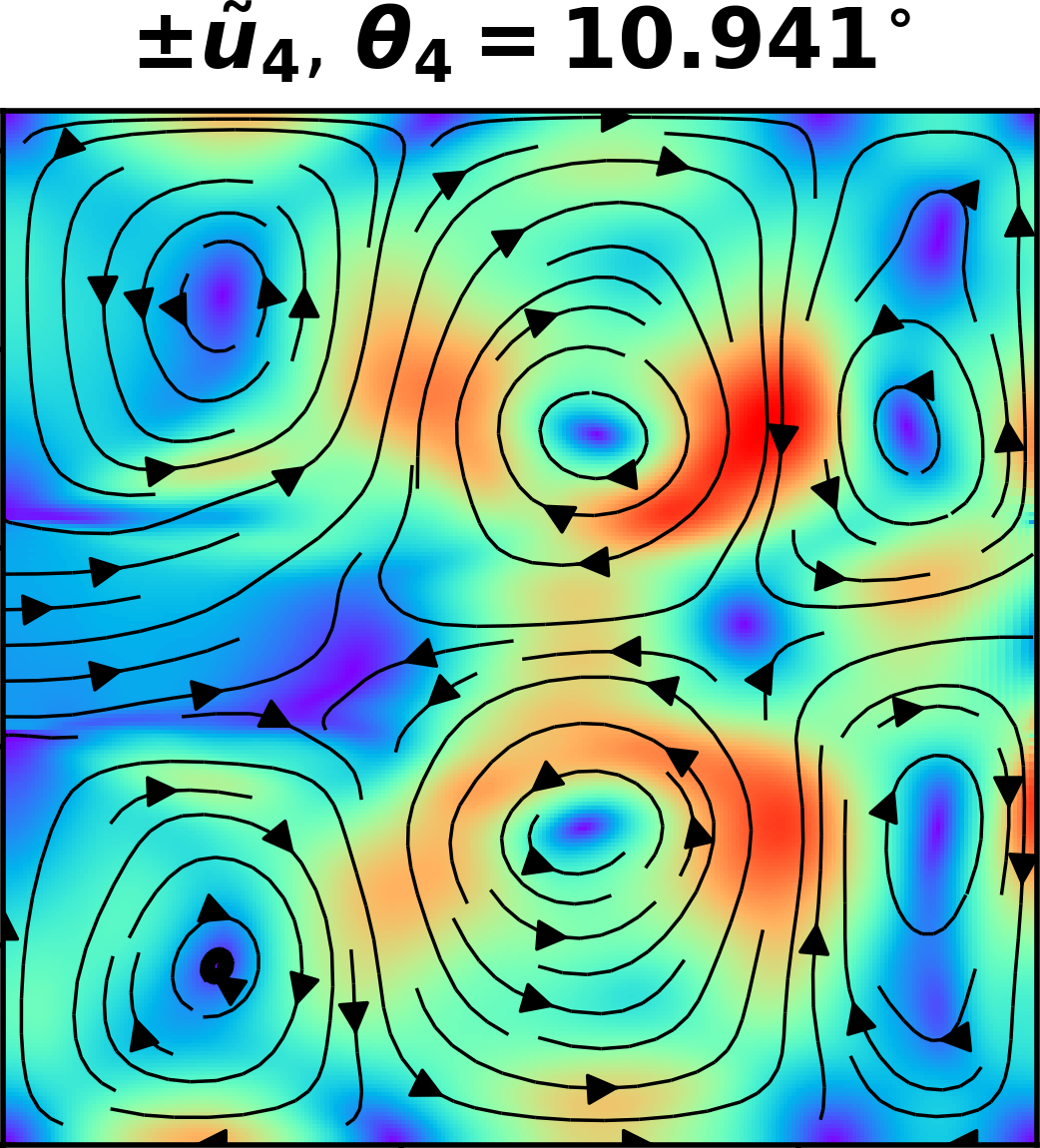}&
     \hspace{-0.15in}   \includegraphics[keepaspectratio, scale= 0.4]{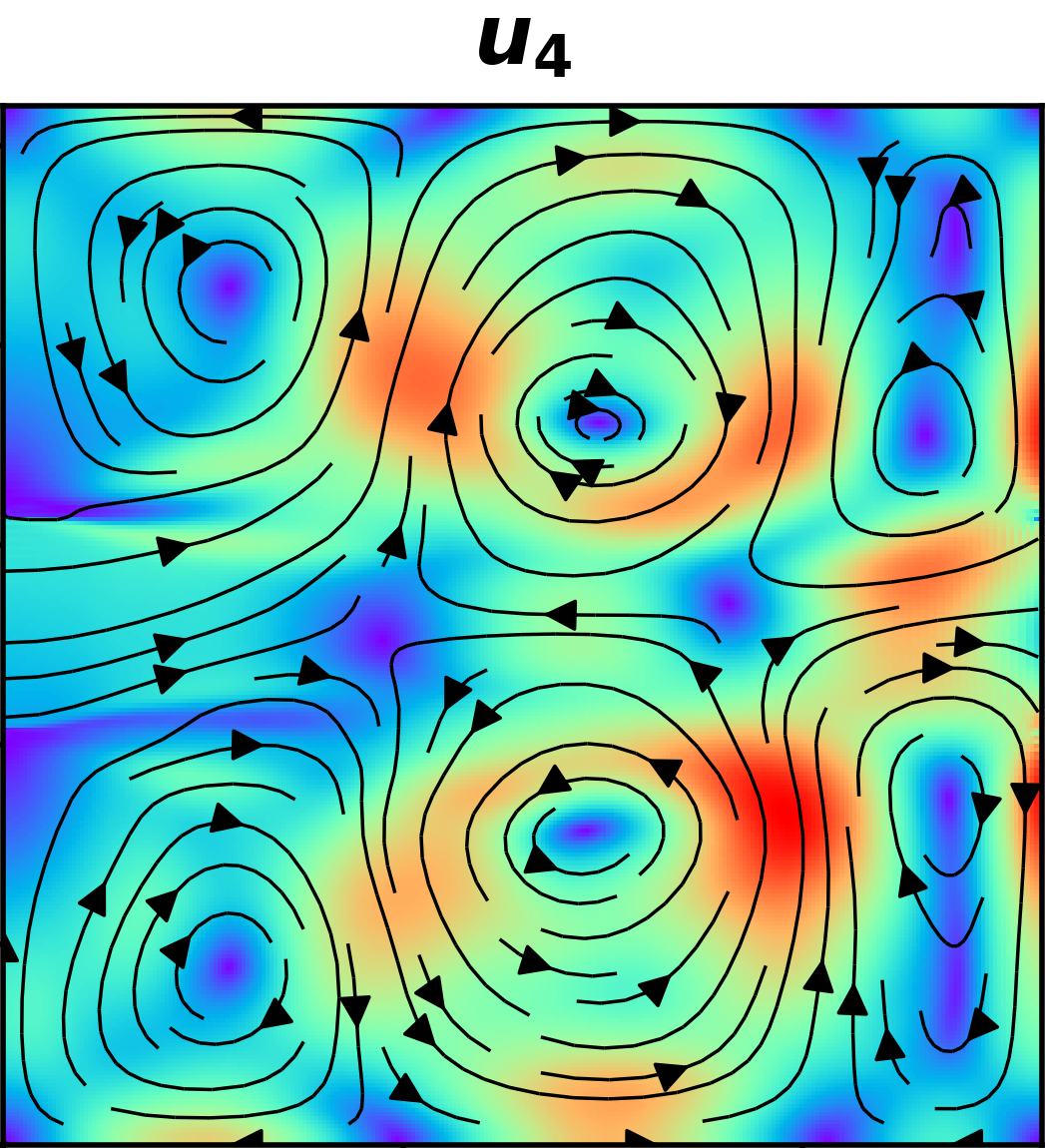}&
    \hspace{-0.15in}    \includegraphics[keepaspectratio, scale= 0.4]{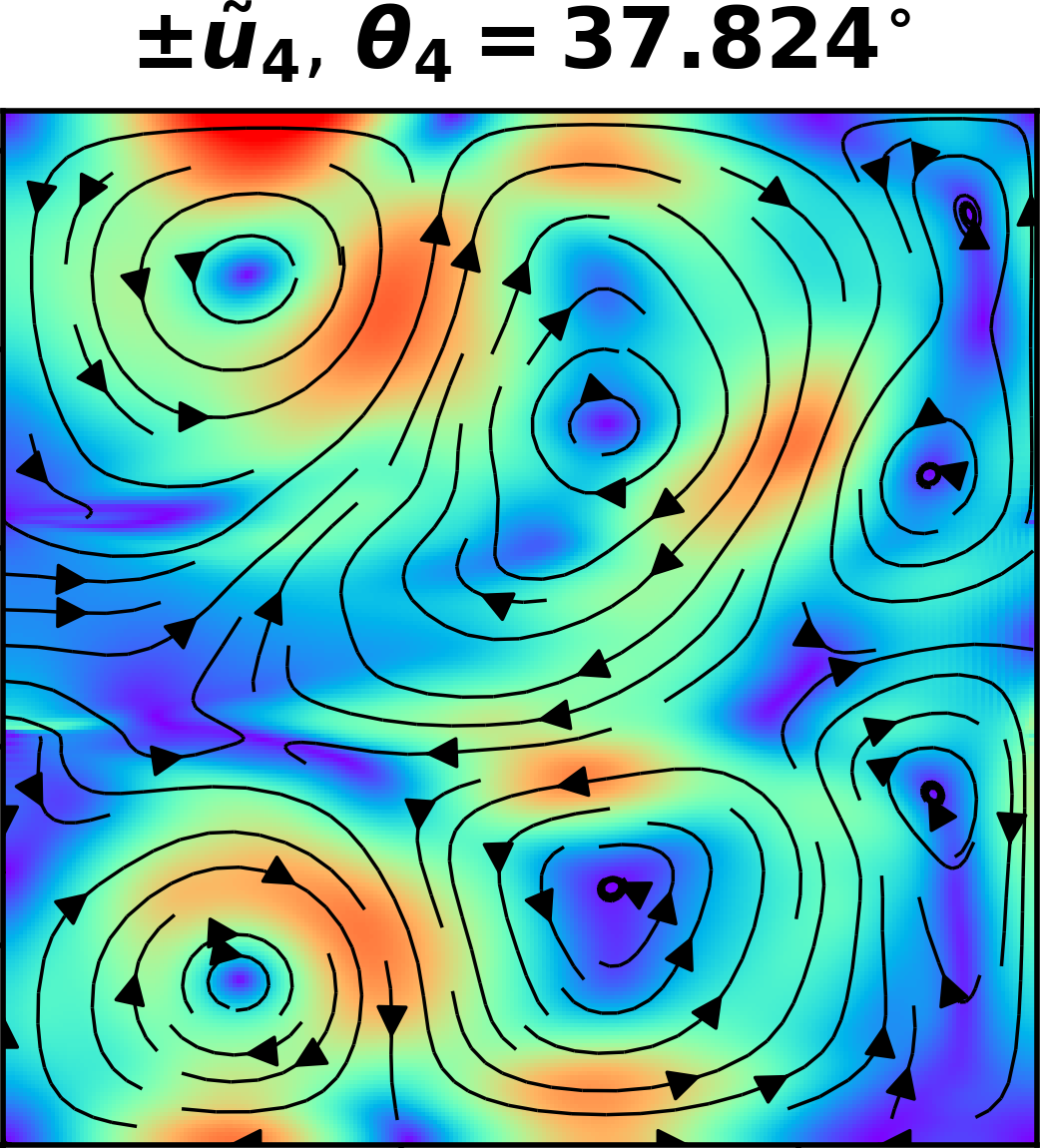}&
     \hspace{-0.15in}    \includegraphics[keepaspectratio, scale= 0.4]{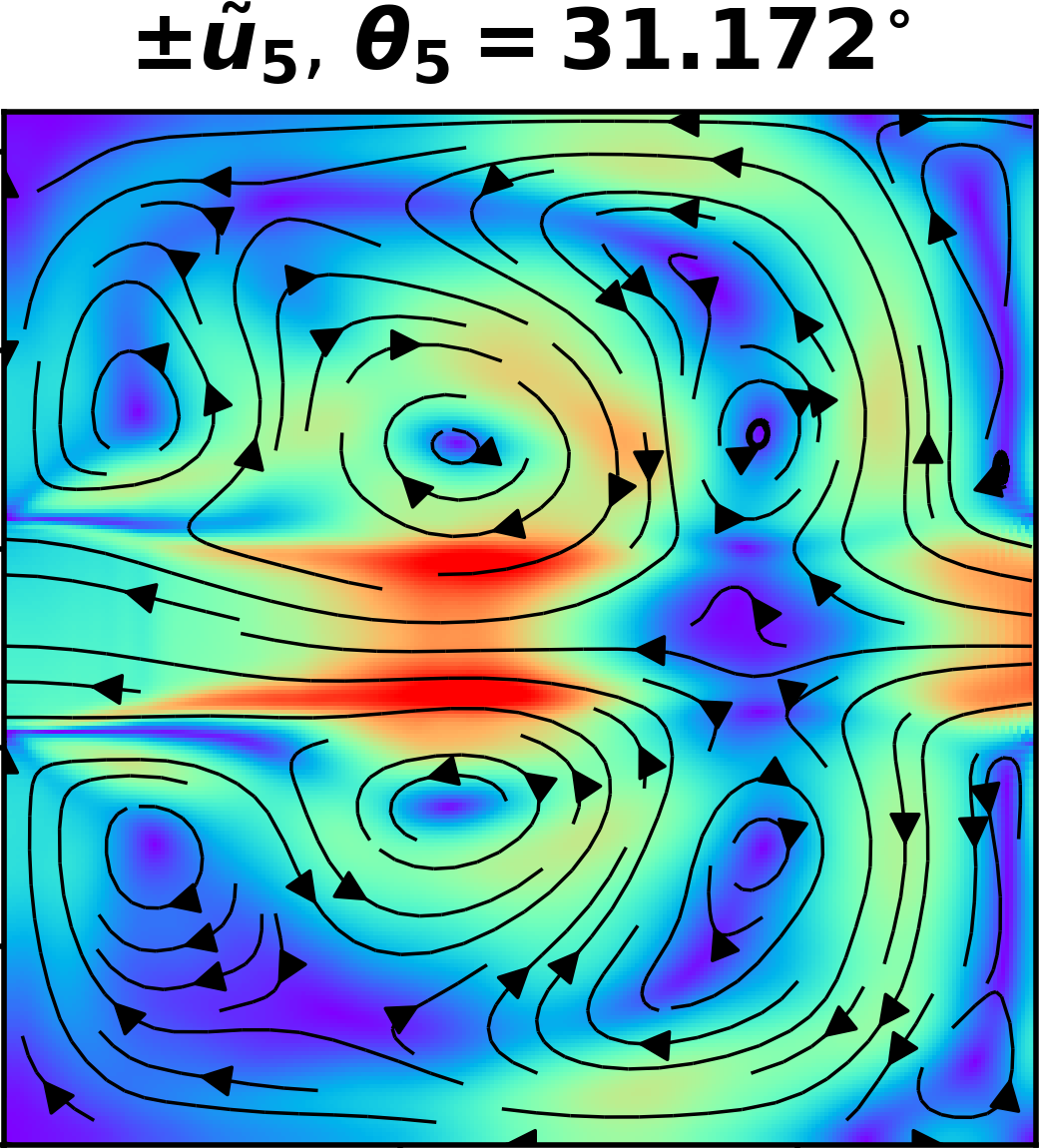}&
     \hspace{-0.15in}   \includegraphics[keepaspectratio, scale= 0.4]{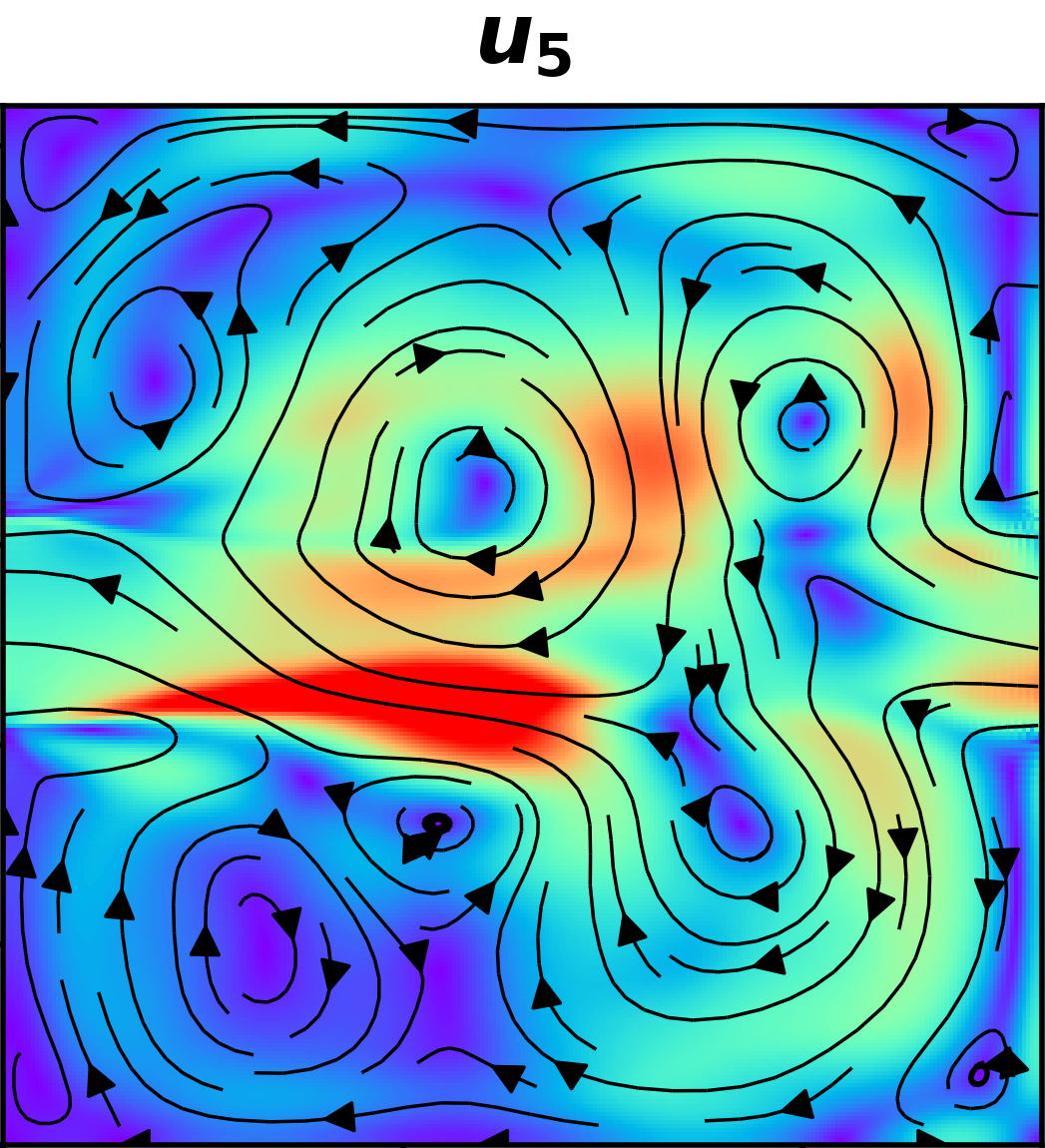}&
     \hspace{-0.15in}    \includegraphics[keepaspectratio, scale= 0.4]{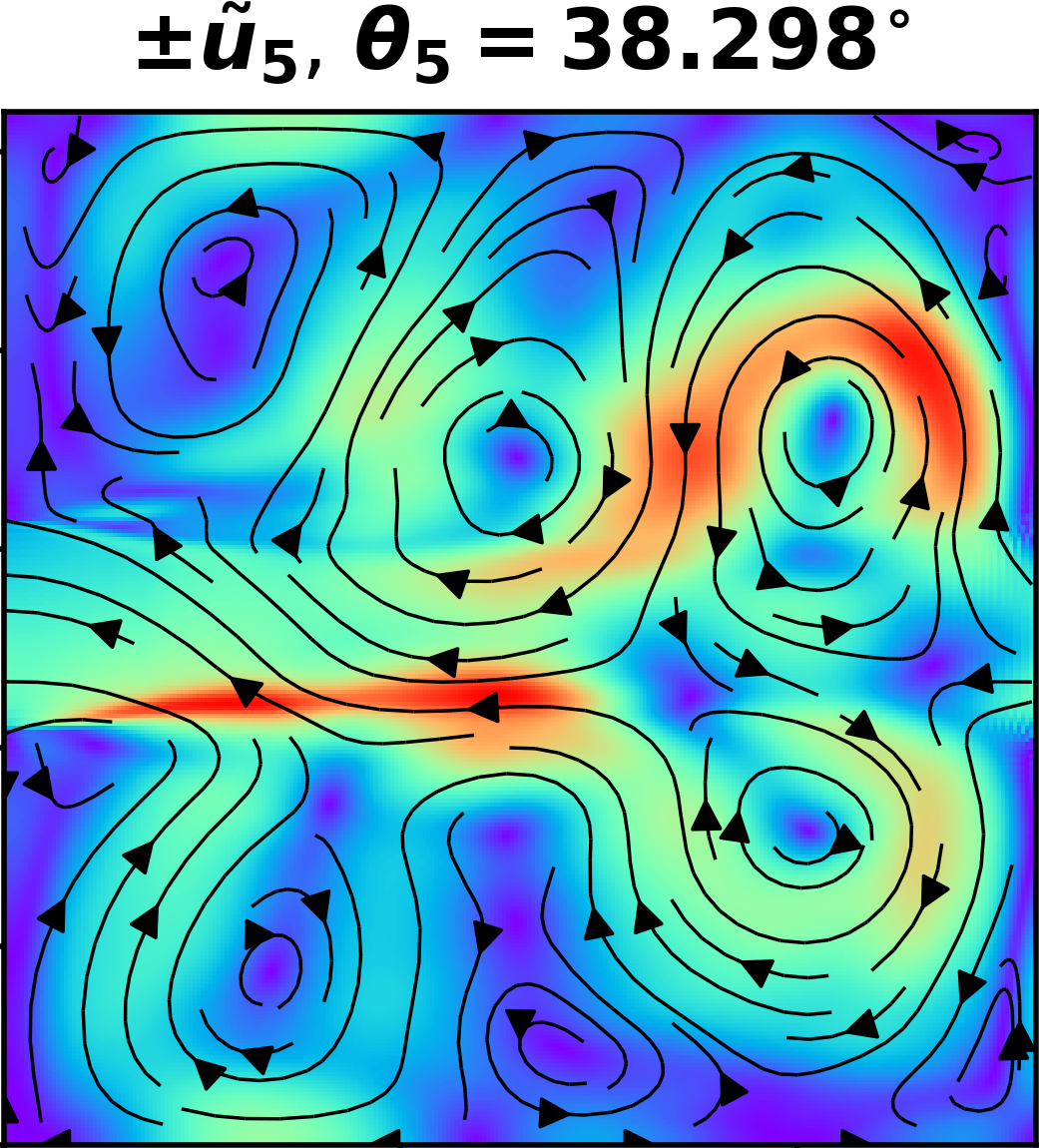}\\
     \multicolumn{6}{c}{Faces dataset}\\
    LTSVD C1&
     \hspace{-0.15in}   SVD&
     \hspace{-0.15in}   CTSVD C1&
     \hspace{-0.15in}   LTSVD C2&
     \hspace{-0.15in}   SVD&
     \hspace{-0.15in}   CTSVD C2\\
    \includegraphics[keepaspectratio, scale= 0.4]{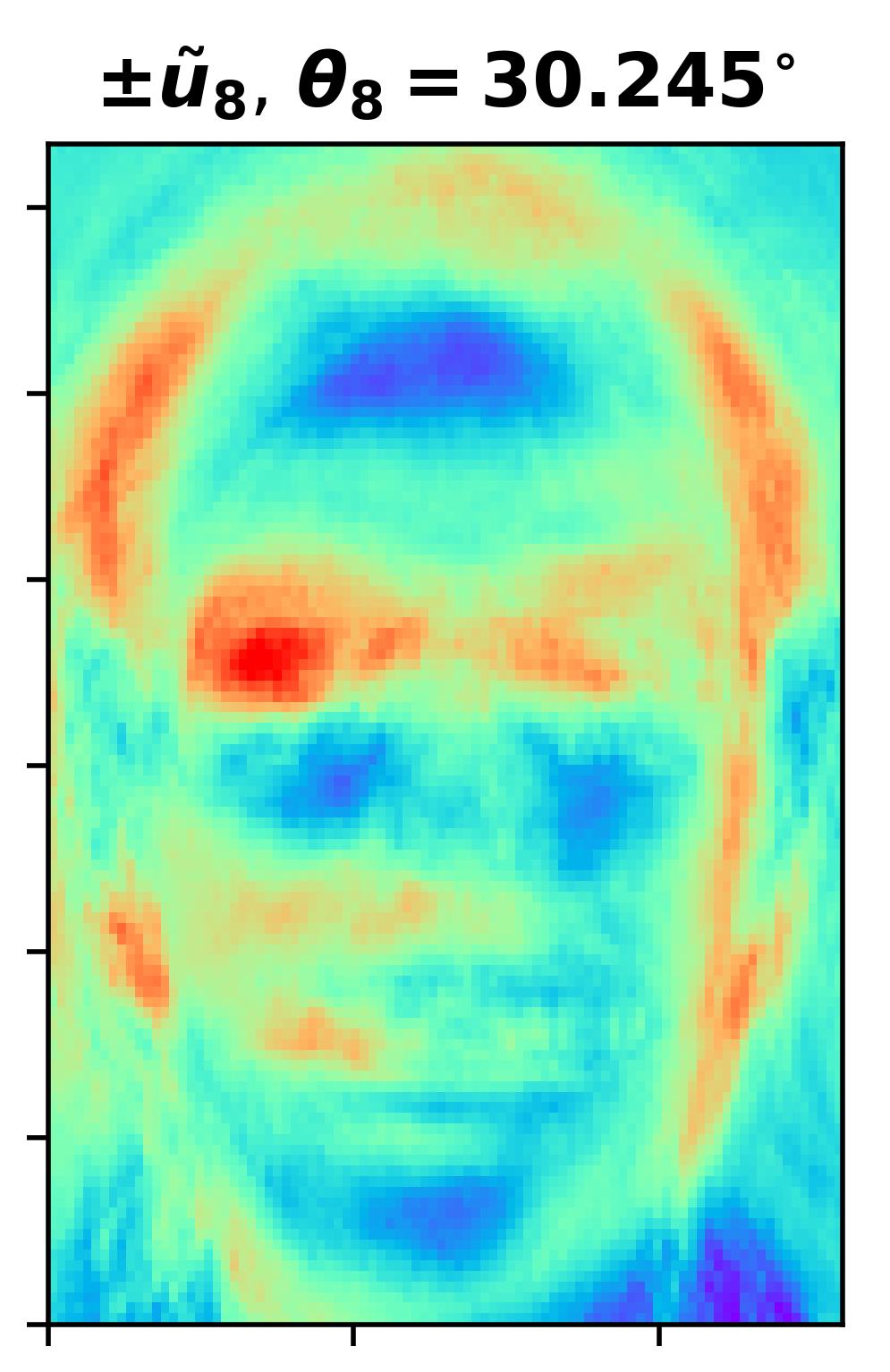}&
     \hspace{-0.15in}   \includegraphics[keepaspectratio, scale= 0.4]{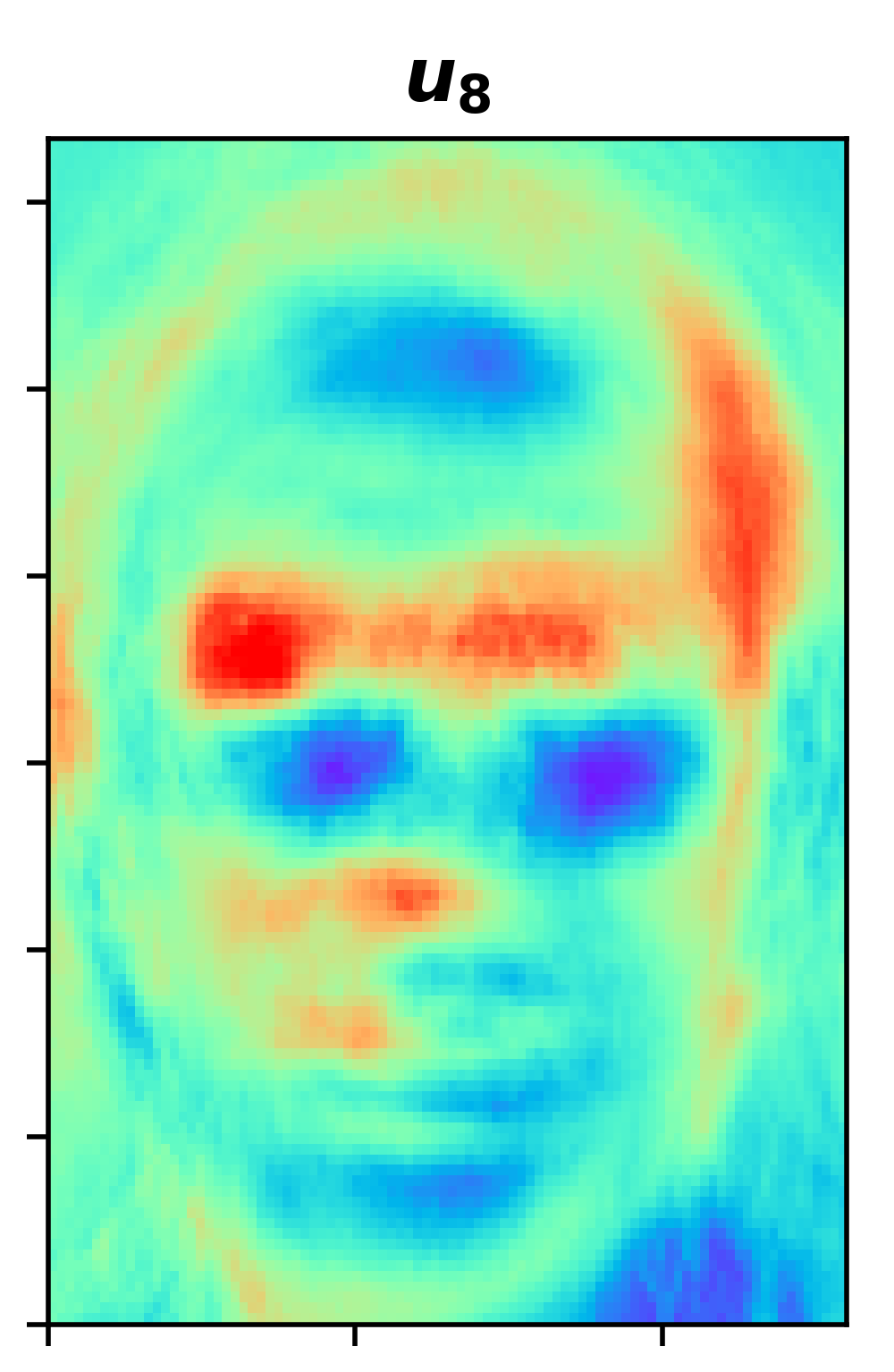}&
     \hspace{-0.15in}   \includegraphics[keepaspectratio, scale= 0.4]{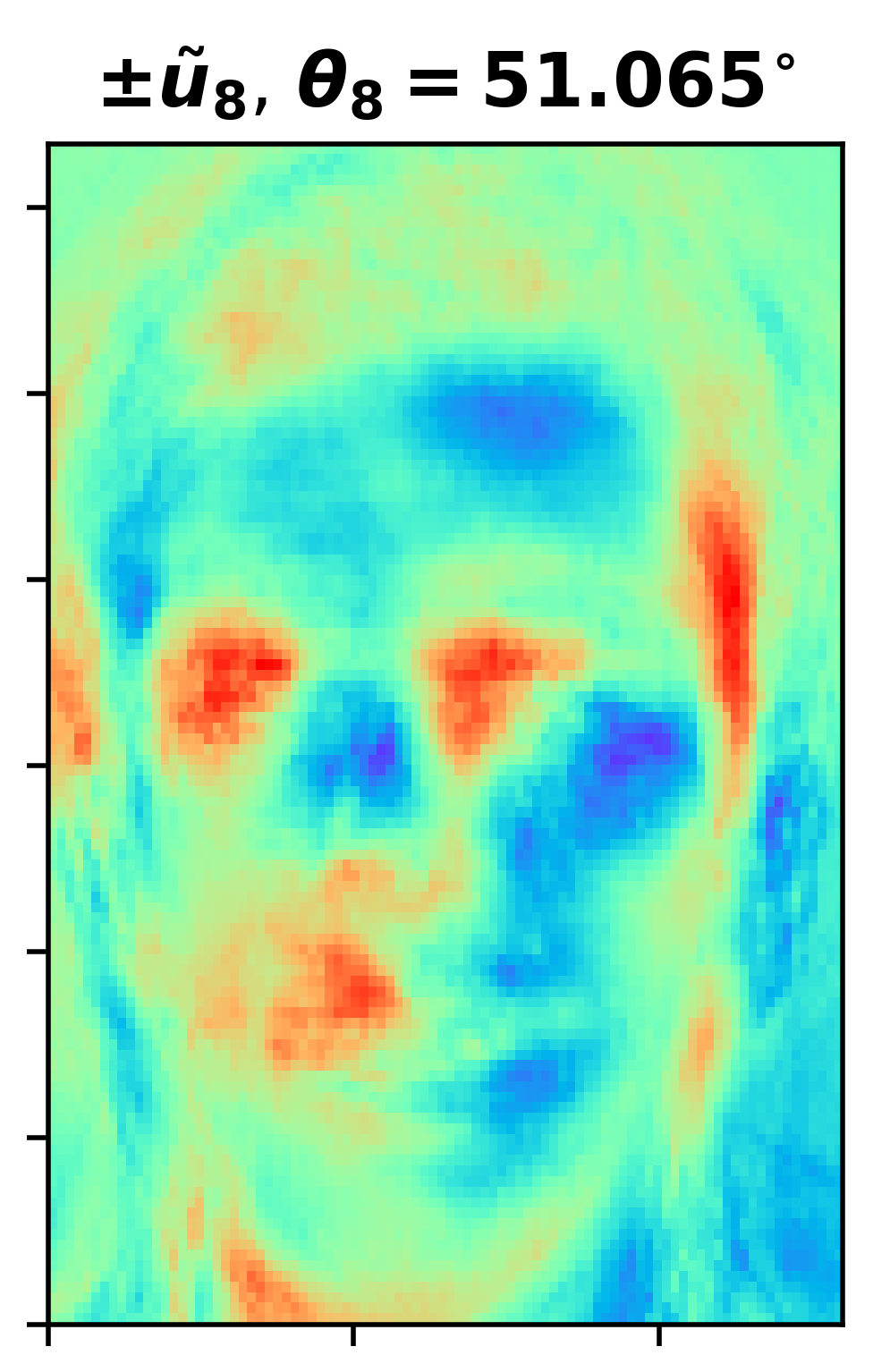}&
     \hspace{-0.15in}    \includegraphics[keepaspectratio, scale= 0.4]{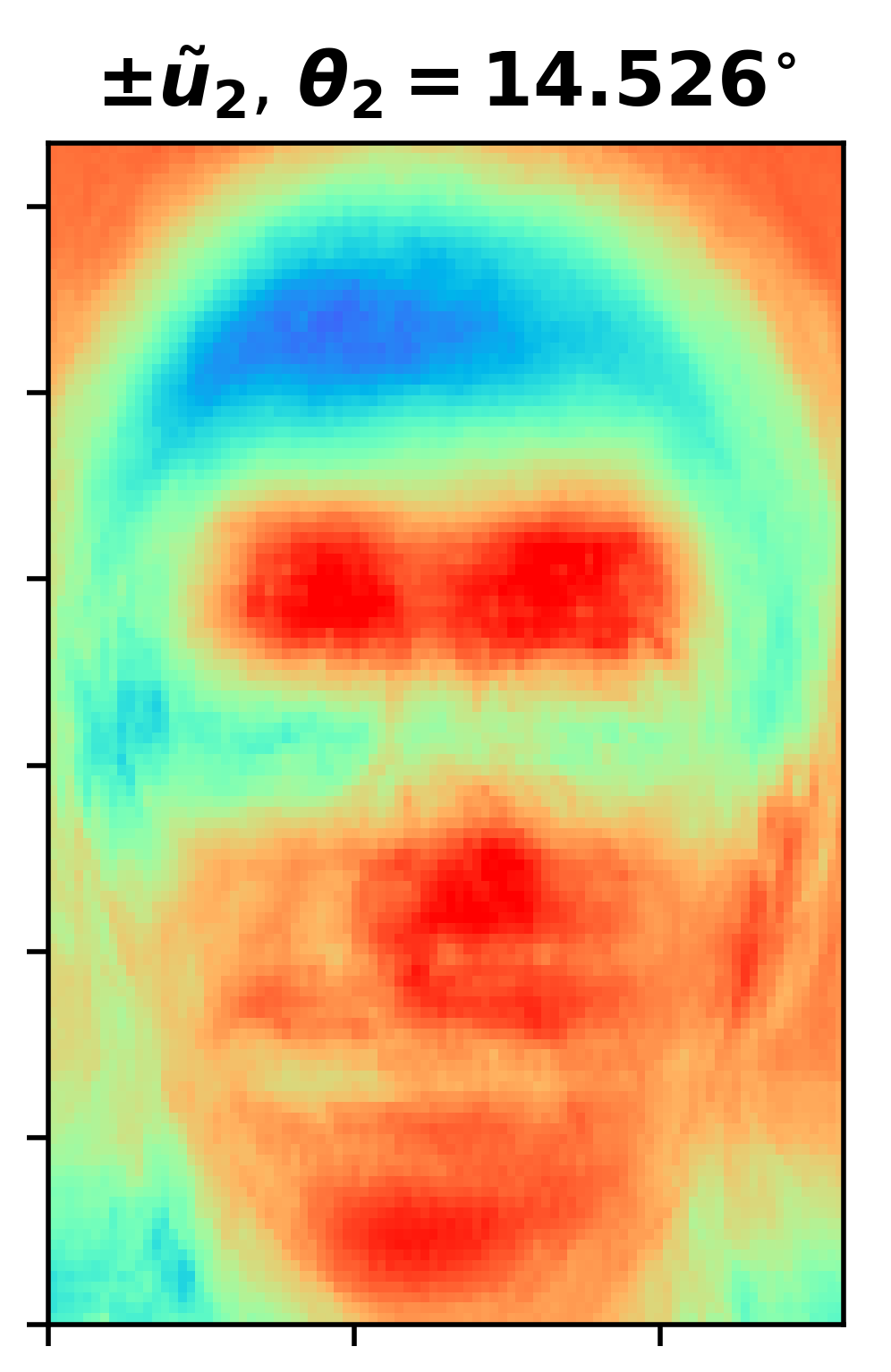}&
     \hspace{-0.15in}   \includegraphics[keepaspectratio, scale= 0.4]{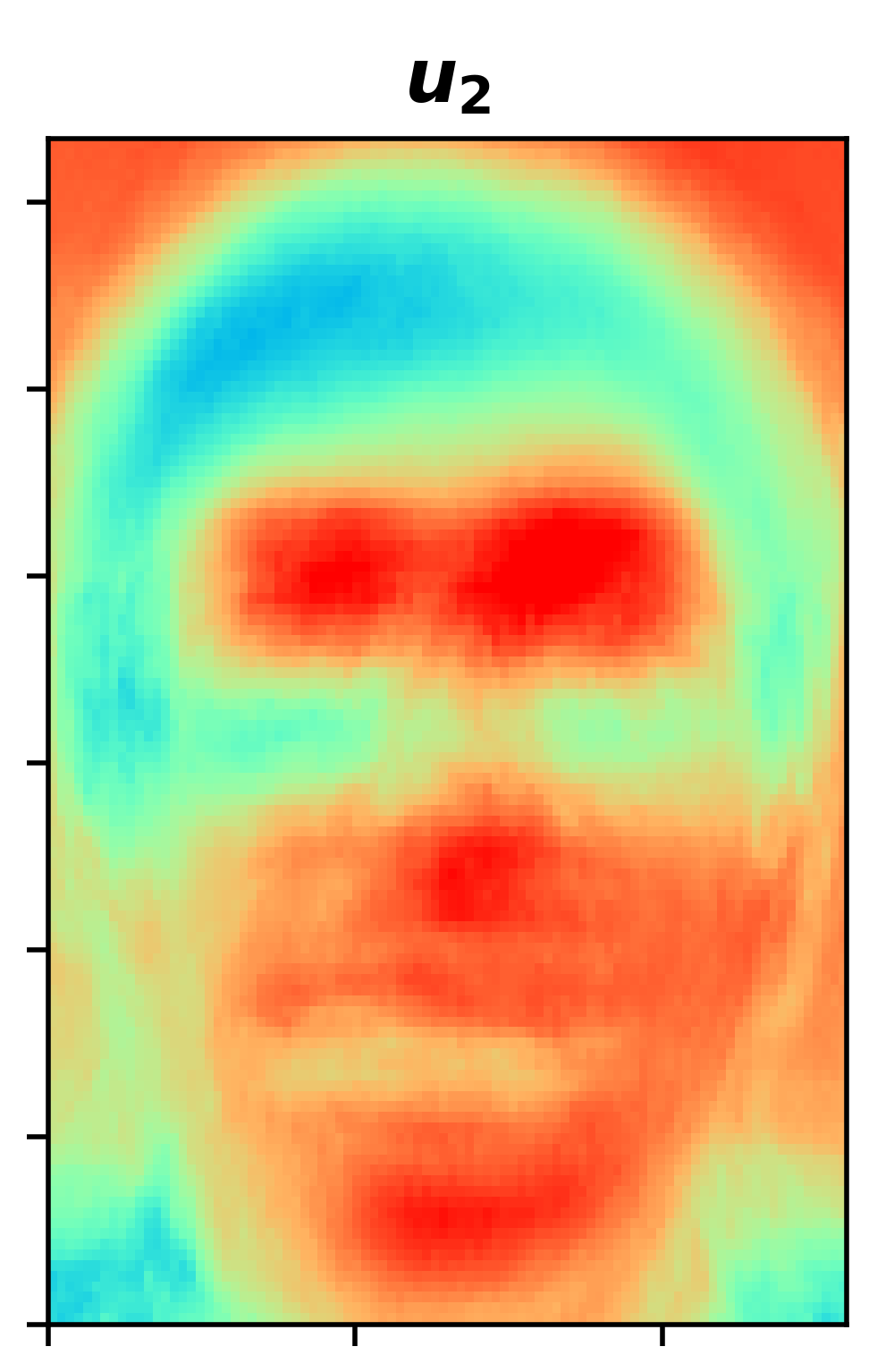}&
     \hspace{-0.15in}    \includegraphics[keepaspectratio, scale= 0.4]{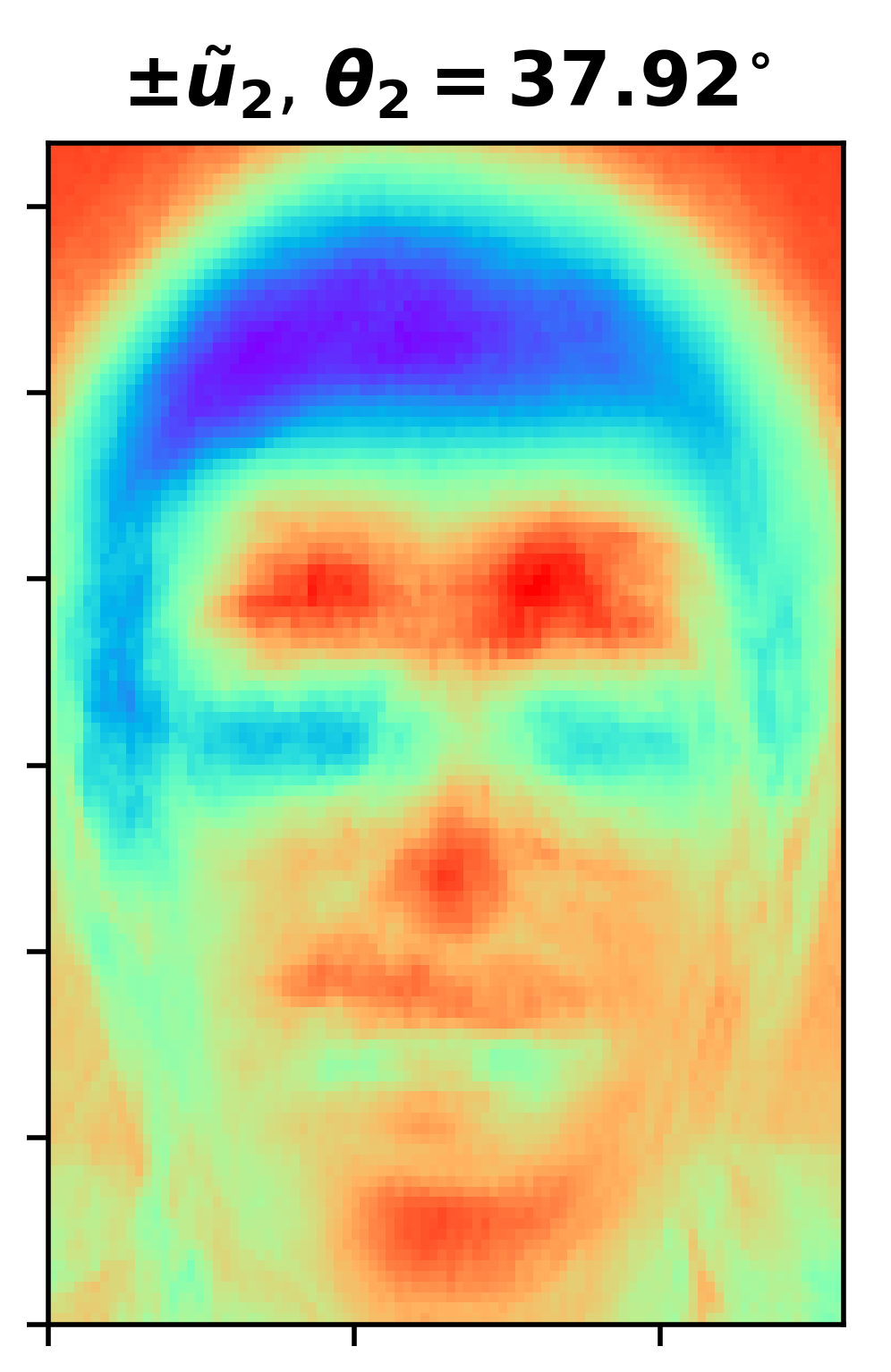}\\
   
  \end{tabular}
\caption{POD modes of various datasets computed by SVD, LTSVD (Algorithm~\ref{alg:LTSVD}) 
and CTSVD (Algorithm~\ref{alg:CTSVD}).
POD mode of $V_{2D}$ consists of both the $x$ and $y$ components of velocity vector at each point in space. 
The color at each point denotes the magnitude of this vector at that point.
Streamlines shown, denote the direction of the flow captured by the POD mode. 
}
\label{fig:POD_2}
\end{figure}

\begin{figure}[h]
  \begin{minipage}{\textwidth}
\centering
  \ref{angPOD_dri}\\
  \end{minipage}
  \begin{minipage}{\textwidth}
\hspace{-0.4in}    \scalebox{0.75}{  \pgfplotstableread{data/angPOD_LTSVD_drineas_V_2D.dat}\4
\pgfplotstableread{data/angPOD_CTSVD_drineas_V_2D.dat}\1

\begin{tikzpicture}[every mark/.append style={line width=2pt, solid}]
\begin{axis}[ 
  title = {$V_{2D}$}, 
  xlabel = $i$,
  ylabel = $\theta_i$,
  xmax = 11.000000,
  ymax=100,
  ymin = -10,
  xtick={2,4,...,10},
  ytick = {0,20,...,80},
  width = 6cm,
  ymajorgrids,
  legend to name = angPOD_dri,
  legend style={at={(0.1,1.05)},anchor=south, legend columns = 6, cells={line width=2pt, solid}},
  legend entries = {LTSVD C1, CTSVD C1, LTSVD C2, CTSVD C2, LTSVD C3, CTSVD C3},
  ]
  \addplot[
		 color=cyan,
		 only marks,
		 mark=triangle,
     mark size = 2pt,
		 x filter/.code={\pgfplotstablegetelem{\coordindex}{[index]0}\of{\4}
\pgfmathtruncatemacro{\temp}{abs(\pgfplotsretval-10)==0? 1 : 0}
\ifnum\temp>0
\pgfplotstablegetelem{\coordindex}{[index]1}\of{\4}
\pgfmathtruncatemacro{\temp}{abs(\pgfplotsretval-0.7)==0? 1 : 0}
\ifnum\temp>0
\pgfplotstablegetelem{\coordindex}{[index]2}\of{\4}
\pgfmathtruncatemacro{\temp}{abs(\pgfplotsretval-0.45)==0? 1 : 0}
\ifnum\temp>0
\relax
			\else
  \def\pgfmathresult{}
\fi
\else
  \def\pgfmathresult{}
\fi
\else
  \def\pgfmathresult{}
\fi
},
		]
                table[ 
                x index = 3,
                y index = 4,
                ]
		{\4};
  \addplot[
		 color=cyan,
		 only marks,
     mark=o,
     mark size = 2pt,
		 x filter/.code={\pgfplotstablegetelem{\coordindex}{[index]0}\of{\1}
\pgfmathtruncatemacro{\temp}{abs(\pgfplotsretval-10)==0? 1 : 0}
\ifnum\temp>0
\pgfplotstablegetelem{\coordindex}{[index]1}\of{\1}
\pgfmathtruncatemacro{\temp}{abs(\pgfplotsretval-1.04)==0? 1 : 0}
\ifnum\temp>0
\pgfplotstablegetelem{\coordindex}{[index]2}\of{\1}
\pgfmathtruncatemacro{\temp}{abs(\pgfplotsretval-0.94)==0? 1 : 0}
\ifnum\temp>0
\relax
			\else
  \def\pgfmathresult{}
\fi
\else
  \def\pgfmathresult{}
\fi
\else
  \def\pgfmathresult{}
\fi
},
		]
                table[ 
                x index = 3,
                y index = 4,
                ]
		{\1};
                \addplot[
		 color=green!50!black,
		 only marks,
		 mark=triangle,
     mark size = 1.5pt,
		 x filter/.code={\pgfplotstablegetelem{\coordindex}{[index]0}\of{\4}
\pgfmathtruncatemacro{\temp}{abs(\pgfplotsretval-5)==0? 1 : 0}
\ifnum\temp>0
\pgfplotstablegetelem{\coordindex}{[index]1}\of{\4}
\pgfmathtruncatemacro{\temp}{abs(\pgfplotsretval-1.0)==0? 1 : 0}
\ifnum\temp>0
\pgfplotstablegetelem{\coordindex}{[index]2}\of{\4}
\pgfmathtruncatemacro{\temp}{abs(\pgfplotsretval-0.1)==0? 1 : 0}
\ifnum\temp>0
\relax
			\else
  \def\pgfmathresult{}
\fi
\else
  \def\pgfmathresult{}
\fi
\else
  \def\pgfmathresult{}
\fi
},
		]
                table[ 
                x index = 3,
                y index = 4,
                ]
		{\4};
                \addplot[
		 color=green!50!black,
		 only marks,
		 mark=o,
		 x filter/.code={\pgfplotstablegetelem{\coordindex}{[index]0}\of{\1}
\pgfmathtruncatemacro{\temp}{abs(\pgfplotsretval-5)==0? 1 : 0}
\ifnum\temp>0
\pgfplotstablegetelem{\coordindex}{[index]1}\of{\1}
\pgfmathtruncatemacro{\temp}{abs(\pgfplotsretval-1.0)==0? 1 : 0}
\ifnum\temp>0
\pgfplotstablegetelem{\coordindex}{[index]2}\of{\1}
\pgfmathtruncatemacro{\temp}{abs(\pgfplotsretval-0.1)==0? 1 : 0}
\ifnum\temp>0
\relax
			\else
  \def\pgfmathresult{}
\fi
\else
  \def\pgfmathresult{}
\fi
\else
  \def\pgfmathresult{}
\fi
},
		]
                table[ 
                x index = 3,
                y index = 4,
                ]
		{\1};
   \addplot[
		 color=red,
		 only marks,
		 mark=triangle,
     mark size = 2pt,
		 x filter/.code={\pgfplotstablegetelem{\coordindex}{[index]0}\of{\4}
\pgfmathtruncatemacro{\temp}{abs(\pgfplotsretval-2)==0? 1 : 0}
\ifnum\temp>0
\pgfplotstablegetelem{\coordindex}{[index]1}\of{\4}
\pgfmathtruncatemacro{\temp}{abs(\pgfplotsretval-1.0)==0? 1 : 0}
\ifnum\temp>0
\pgfplotstablegetelem{\coordindex}{[index]2}\of{\4}
\pgfmathtruncatemacro{\temp}{abs(\pgfplotsretval-0.8)==0? 1 : 0}
\ifnum\temp>0
\relax
			\else
  \def\pgfmathresult{}
\fi
\else
  \def\pgfmathresult{}
\fi
\else
  \def\pgfmathresult{}
\fi
},
		]
                table[ 
                x index = 3,
                y index = 4,
                ]
		{\4};
  \addplot[
		 color=red,
		 only marks,
     mark=o,
     mark size = 2pt,
		 x filter/.code={\pgfplotstablegetelem{\coordindex}{[index]0}\of{\1}
\pgfmathtruncatemacro{\temp}{abs(\pgfplotsretval-2)==0? 1 : 0}
\ifnum\temp>0
\pgfplotstablegetelem{\coordindex}{[index]1}\of{\1}
\pgfmathtruncatemacro{\temp}{abs(\pgfplotsretval-1.0)==0? 1 : 0}
\ifnum\temp>0
\pgfplotstablegetelem{\coordindex}{[index]2}\of{\1}
\pgfmathtruncatemacro{\temp}{abs(\pgfplotsretval-0.8)==0? 1 : 0}
\ifnum\temp>0
\relax
			\else
  \def\pgfmathresult{}
\fi
\else
  \def\pgfmathresult{}
\fi
\else
  \def\pgfmathresult{}
\fi
},
		]
                table[ 
                x index = 3,
                y index = 4,
                ]
		{\1};
    \end{axis}
    \end{tikzpicture}}
 \hspace{-0.22in}    \scalebox{0.75}{ \pgfplotstableread{data/angPOD_LTSVD_drineas_faces.dat}\4
\pgfplotstableread{data/angPOD_CTSVD_drineas_faces.dat}\1

\begin{tikzpicture}[every mark/.append style={line width=2pt, solid}]
\begin{axis}[ 
  title = {Faces dataset}, 
  xlabel = $i$,
  xmax = 11.000000,
  ymax=100,
  ymin = -10,
  xtick={2,4,...,10},
  ytick = {0,20,...,80},
  yticklabels={,,},
  width = 6cm,
  ymajorgrids,
  ]
  \addplot[
		 color=cyan,
		 only marks,
		 mark=triangle,
     mark size = 2pt,
		 x filter/.code={\pgfplotstablegetelem{\coordindex}{[index]0}\of{\4}
\pgfmathtruncatemacro{\temp}{abs(\pgfplotsretval-10)==0? 1 : 0}
\ifnum\temp>0
\pgfplotstablegetelem{\coordindex}{[index]1}\of{\4}
\pgfmathtruncatemacro{\temp}{abs(\pgfplotsretval-0.75)==0? 1 : 0}
\ifnum\temp>0
\pgfplotstablegetelem{\coordindex}{[index]2}\of{\4}
\pgfmathtruncatemacro{\temp}{abs(\pgfplotsretval-0.8)==0? 1 : 0}
\ifnum\temp>0
\relax
			\else
  \def\pgfmathresult{}
\fi
\else
  \def\pgfmathresult{}
\fi
\else
  \def\pgfmathresult{}
\fi
},
		]
                table[ 
                x index = 3,
                y index = 4,
                ]
		{\4};
  \addplot[
		 color=cyan,
		 only marks,
     mark=o,
     mark size = 2pt,
		 x filter/.code={\pgfplotstablegetelem{\coordindex}{[index]0}\of{\1}
\pgfmathtruncatemacro{\temp}{abs(\pgfplotsretval-10)==0? 1 : 0}
\ifnum\temp>0
\pgfplotstablegetelem{\coordindex}{[index]1}\of{\1}
\pgfmathtruncatemacro{\temp}{abs(\pgfplotsretval-1.3)==0? 1 : 0}
\ifnum\temp>0
\pgfplotstablegetelem{\coordindex}{[index]2}\of{\1}
\pgfmathtruncatemacro{\temp}{abs(\pgfplotsretval-1)==0? 1 : 0}
\ifnum\temp>0
\relax
			\else
  \def\pgfmathresult{}
\fi
\else
  \def\pgfmathresult{}
\fi
\else
  \def\pgfmathresult{}
\fi
},
		]
                table[ 
                x index = 3,
                y index = 4,
                ]
		{\1};
                \addplot[
		 color=green!50!black,
		 only marks,
		 mark=triangle,
     mark size = 1.5pt,
		 x filter/.code={\pgfplotstablegetelem{\coordindex}{[index]0}\of{\4}
\pgfmathtruncatemacro{\temp}{abs(\pgfplotsretval-5)==0? 1 : 0}
\ifnum\temp>0
\pgfplotstablegetelem{\coordindex}{[index]1}\of{\4}
\pgfmathtruncatemacro{\temp}{abs(\pgfplotsretval-0.75)==0? 1 : 0}
\ifnum\temp>0
\pgfplotstablegetelem{\coordindex}{[index]2}\of{\4}
\pgfmathtruncatemacro{\temp}{abs(\pgfplotsretval-0.75)==0? 1 : 0}
\ifnum\temp>0
\relax
			\else
  \def\pgfmathresult{}
\fi
\else
  \def\pgfmathresult{}
\fi
\else
  \def\pgfmathresult{}
\fi
},
		]
                table[ 
                x index = 3,
                y index = 4,
                ]
		{\4};
                \addplot[
		 color=green!50!black,
		 only marks,
		 mark=o,
		 x filter/.code={\pgfplotstablegetelem{\coordindex}{[index]0}\of{\1}
\pgfmathtruncatemacro{\temp}{abs(\pgfplotsretval-5)==0? 1 : 0}
\ifnum\temp>0
\pgfplotstablegetelem{\coordindex}{[index]1}\of{\1}
\pgfmathtruncatemacro{\temp}{abs(\pgfplotsretval-1.0)==0? 1 : 0}
\ifnum\temp>0
\pgfplotstablegetelem{\coordindex}{[index]2}\of{\1}
\pgfmathtruncatemacro{\temp}{abs(\pgfplotsretval-1)==0? 1 : 0}
\ifnum\temp>0
\relax
			\else
  \def\pgfmathresult{}
\fi
\else
  \def\pgfmathresult{}
\fi
\else
  \def\pgfmathresult{}
\fi
},
		]
                table[ 
                x index = 3,
                y index = 4,
                ]
		{\1};
   \addplot[
		 color=red,
		 only marks,
		 mark=triangle,
     mark size = 2pt,
		 x filter/.code={\pgfplotstablegetelem{\coordindex}{[index]0}\of{\4}
\pgfmathtruncatemacro{\temp}{abs(\pgfplotsretval-2)==0? 1 : 0}
\ifnum\temp>0
\pgfplotstablegetelem{\coordindex}{[index]1}\of{\4}
\pgfmathtruncatemacro{\temp}{abs(\pgfplotsretval-1.0)==0? 1 : 0}
\ifnum\temp>0
\pgfplotstablegetelem{\coordindex}{[index]2}\of{\4}
\pgfmathtruncatemacro{\temp}{abs(\pgfplotsretval-0.8)==0? 1 : 0}
\ifnum\temp>0
\relax
			\else
  \def\pgfmathresult{}
\fi
\else
  \def\pgfmathresult{}
\fi
\else
  \def\pgfmathresult{}
\fi
},
		]
                table[ 
                x index = 3,
                y index = 4,
                ]
		{\4};
  \addplot[
		 color=red,
		 only marks,
     mark=o,
     mark size = 2pt,
		 x filter/.code={\pgfplotstablegetelem{\coordindex}{[index]0}\of{\1}
\pgfmathtruncatemacro{\temp}{abs(\pgfplotsretval-2)==0? 1 : 0}
\ifnum\temp>0
\pgfplotstablegetelem{\coordindex}{[index]1}\of{\1}
\pgfmathtruncatemacro{\temp}{abs(\pgfplotsretval-1.0)==0? 1 : 0}
\ifnum\temp>0
\pgfplotstablegetelem{\coordindex}{[index]2}\of{\1}
\pgfmathtruncatemacro{\temp}{abs(\pgfplotsretval-0.8)==0? 1 : 0}
\ifnum\temp>0
\relax
			\else
  \def\pgfmathresult{}
\fi
\else
  \def\pgfmathresult{}
\fi
\else
  \def\pgfmathresult{}
\fi
},
		]
                table[ 
                x index = 3,
                y index = 4,
                ]
		{\1};
                
    \end{axis}
    \end{tikzpicture}}
\hspace{-0.2in}  \scalebox{0.75}{\pgfplotstableread{data/angPOD_LTSVD_drineas_croppedYalefaces.dat}\4
\pgfplotstableread{data/angPOD_CTSVD_drineas_croppedYalefaces.dat}\1

\begin{tikzpicture}[every mark/.append style={line width=2pt, solid}]
\begin{axis}[
    title = {cropped Yale faces}, 
    xlabel = $i$,
    xmax = 11.000000,
     ymax=100,
    ymin = -10,
    xtick={2,4,...,10},
    ytick = {0,20,...,80},
    yticklabels={,,},
    width = 6cm,
    ymajorgrids,
  ]
  \addplot[
		 color=cyan,
		 only marks,
		 mark=triangle,
     mark size = 2pt,
		 x filter/.code={\pgfplotstablegetelem{\coordindex}{[index]0}\of{\4}
\pgfmathtruncatemacro{\temp}{abs(\pgfplotsretval-10)==0? 1 : 0}
\ifnum\temp>0
\pgfplotstablegetelem{\coordindex}{[index]1}\of{\4}
\pgfmathtruncatemacro{\temp}{abs(\pgfplotsretval-0.46)==0? 1 : 0}
\ifnum\temp>0
\pgfplotstablegetelem{\coordindex}{[index]2}\of{\4}
\pgfmathtruncatemacro{\temp}{abs(\pgfplotsretval-0.44)==0? 1 : 0}
\ifnum\temp>0
\relax
			\else
  \def\pgfmathresult{}
\fi
\else
  \def\pgfmathresult{}
\fi
\else
  \def\pgfmathresult{}
\fi
},
		]
                table[ 
                x index = 3,
                y index = 4,
                ]
		{\4};
  \addplot[
		 color=cyan,
		 only marks,
     mark=o,
     mark size = 2pt,
		 x filter/.code={\pgfplotstablegetelem{\coordindex}{[index]0}\of{\1}
\pgfmathtruncatemacro{\temp}{abs(\pgfplotsretval-10)==0? 1 : 0}
\ifnum\temp>0
\pgfplotstablegetelem{\coordindex}{[index]1}\of{\1}
\pgfmathtruncatemacro{\temp}{abs(\pgfplotsretval-0.9)==0? 1 : 0}
\ifnum\temp>0
\pgfplotstablegetelem{\coordindex}{[index]2}\of{\1}
\pgfmathtruncatemacro{\temp}{abs(\pgfplotsretval-0.7)==0? 1 : 0}
\ifnum\temp>0
\relax
			\else
  \def\pgfmathresult{}
\fi
\else
  \def\pgfmathresult{}
\fi
\else
  \def\pgfmathresult{}
\fi
},
		]
                table[ 
                x index = 3,
                y index = 4,
                ]
		{\1};
                \addplot[
		 color=green!50!black,
		 only marks,
		 mark=triangle,
     mark size = 1.5pt,
		 x filter/.code={\pgfplotstablegetelem{\coordindex}{[index]0}\of{\4}
\pgfmathtruncatemacro{\temp}{abs(\pgfplotsretval-5)==0? 1 : 0}
\ifnum\temp>0
\pgfplotstablegetelem{\coordindex}{[index]1}\of{\4}
\pgfmathtruncatemacro{\temp}{abs(\pgfplotsretval-1.0)==0? 1 : 0}
\ifnum\temp>0
\pgfplotstablegetelem{\coordindex}{[index]2}\of{\4}
\pgfmathtruncatemacro{\temp}{abs(\pgfplotsretval-0.1)==0? 1 : 0}
\ifnum\temp>0
\relax
			\else
  \def\pgfmathresult{}
\fi
\else
  \def\pgfmathresult{}
\fi
\else
  \def\pgfmathresult{}
\fi
},
		]
                table[ 
                x index = 3,
                y index = 4,
                ]
		{\4};
                \addplot[
		 color=green!50!black,
		 only marks,
		 mark=o,
		 x filter/.code={\pgfplotstablegetelem{\coordindex}{[index]0}\of{\1}
\pgfmathtruncatemacro{\temp}{abs(\pgfplotsretval-5)==0? 1 : 0}
\ifnum\temp>0
\pgfplotstablegetelem{\coordindex}{[index]1}\of{\1}
\pgfmathtruncatemacro{\temp}{abs(\pgfplotsretval-1.0)==0? 1 : 0}
\ifnum\temp>0
\pgfplotstablegetelem{\coordindex}{[index]2}\of{\1}
\pgfmathtruncatemacro{\temp}{abs(\pgfplotsretval-0.1)==0? 1 : 0}
\ifnum\temp>0
\relax
			\else
  \def\pgfmathresult{}
\fi
\else
  \def\pgfmathresult{}
\fi
\else
  \def\pgfmathresult{}
\fi
},
		]
                table[ 
                x index = 3,
                y index = 4,
                ]
		{\1};
   \addplot[
		 color=red,
		 only marks,
		 mark=triangle,
     mark size = 2pt,
		 x filter/.code={\pgfplotstablegetelem{\coordindex}{[index]0}\of{\4}
\pgfmathtruncatemacro{\temp}{abs(\pgfplotsretval-2)==0? 1 : 0}
\ifnum\temp>0
\pgfplotstablegetelem{\coordindex}{[index]1}\of{\4}
\pgfmathtruncatemacro{\temp}{abs(\pgfplotsretval-1.0)==0? 1 : 0}
\ifnum\temp>0
\pgfplotstablegetelem{\coordindex}{[index]2}\of{\4}
\pgfmathtruncatemacro{\temp}{abs(\pgfplotsretval-0.8)==0? 1 : 0}
\ifnum\temp>0
\relax
			\else
  \def\pgfmathresult{}
\fi
\else
  \def\pgfmathresult{}
\fi
\else
  \def\pgfmathresult{}
\fi
},
		]
                table[ 
                x index = 3,
                y index = 4,
                ]
		{\4};
  \addplot[
		 color=red,
		 only marks,
     mark=o,
     mark size = 2pt,
		 x filter/.code={\pgfplotstablegetelem{\coordindex}{[index]0}\of{\1}
\pgfmathtruncatemacro{\temp}{abs(\pgfplotsretval-2)==0? 1 : 0}
\ifnum\temp>0
\pgfplotstablegetelem{\coordindex}{[index]1}\of{\1}
\pgfmathtruncatemacro{\temp}{abs(\pgfplotsretval-1.0)==0? 1 : 0}
\ifnum\temp>0
\pgfplotstablegetelem{\coordindex}{[index]2}\of{\1}
\pgfmathtruncatemacro{\temp}{abs(\pgfplotsretval-0.8)==0? 1 : 0}
\ifnum\temp>0
\relax
			\else
  \def\pgfmathresult{}
\fi
\else
  \def\pgfmathresult{}
\fi
\else
  \def\pgfmathresult{}
\fi
},
		]
                table[ 
                x index = 3,
                y index = 4,
                ]
		{\1};
                
    \end{axis}
    \end{tikzpicture}}
\hspace{-0.22in}     \scalebox{0.75}{ \pgfplotstableread{data/angPOD_LTSVD_drineas_ext_Yalefaces_wo_Ambient.dat}\4
\pgfplotstableread{data/angPOD_CTSVD_drineas_ext_Yalefaces_wo_Ambient.dat}\1

\begin{tikzpicture}[every mark/.append style={line width=2pt, solid}]
\begin{axis}[
    title = {Yale faces}, 
    xlabel = $i$,
    xmax = 21.000000,
    ymax=100,
    ymin = -10,
    ytick = {0,20,...,80},
    yticklabels={,,},
    width = 6cm,
    ymajorgrids,
  ]
  \addplot[
		 color=cyan,
		 only marks,
		 mark=triangle,
     mark size = 2pt,
		 x filter/.code={\pgfplotstablegetelem{\coordindex}{[index]0}\of{\4}
\pgfmathtruncatemacro{\temp}{abs(\pgfplotsretval-20)==0? 1 : 0}
\ifnum\temp>0
\pgfplotstablegetelem{\coordindex}{[index]1}\of{\4}
\pgfmathtruncatemacro{\temp}{abs(\pgfplotsretval-0.35)==0? 1 : 0}
\ifnum\temp>0
\pgfplotstablegetelem{\coordindex}{[index]2}\of{\4}
\pgfmathtruncatemacro{\temp}{abs(\pgfplotsretval-0.35)==0? 1 : 0}
\ifnum\temp>0
\relax
			\else
  \def\pgfmathresult{}
\fi
\else
  \def\pgfmathresult{}
\fi
\else
  \def\pgfmathresult{}
\fi
},
		]
                table[ 
                x index = 3,
                y index = 4,
                ]
		{\4};
  \addplot[
		 color=cyan,
		 only marks,
     mark=o,
     mark size = 2pt,
		 x filter/.code={\pgfplotstablegetelem{\coordindex}{[index]0}\of{\1}
\pgfmathtruncatemacro{\temp}{abs(\pgfplotsretval-20)==0? 1 : 0}
\ifnum\temp>0
\pgfplotstablegetelem{\coordindex}{[index]1}\of{\1}
\pgfmathtruncatemacro{\temp}{abs(\pgfplotsretval-0.83)==0? 1 : 0}
\ifnum\temp>0
\pgfplotstablegetelem{\coordindex}{[index]2}\of{\1}
\pgfmathtruncatemacro{\temp}{abs(\pgfplotsretval-0.9)==0? 1 : 0}
\ifnum\temp>0
\relax
			\else
  \def\pgfmathresult{}
\fi
\else
  \def\pgfmathresult{}
\fi
\else
  \def\pgfmathresult{}
\fi
},
		]
                table[ 
                x index = 3,
                y index = 4,
                ]
		{\1};
                \addplot[
		 color=green!50!black,
		 only marks,
		 mark=triangle,
     mark size = 1.5pt,
		 x filter/.code={\pgfplotstablegetelem{\coordindex}{[index]0}\of{\4}
\pgfmathtruncatemacro{\temp}{abs(\pgfplotsretval-10)==0? 1 : 0}
\ifnum\temp>0
\pgfplotstablegetelem{\coordindex}{[index]1}\of{\4}
\pgfmathtruncatemacro{\temp}{abs(\pgfplotsretval-0.3)==0? 1 : 0}
\ifnum\temp>0
\pgfplotstablegetelem{\coordindex}{[index]2}\of{\4}
\pgfmathtruncatemacro{\temp}{abs(\pgfplotsretval-0.25)==0? 1 : 0}
\ifnum\temp>0
\relax
			\else
  \def\pgfmathresult{}
\fi
\else
  \def\pgfmathresult{}
\fi
\else
  \def\pgfmathresult{}
\fi
},
		]
                table[ 
                x index = 3,
                y index = 4,
                ]
		{\4};
                \addplot[
		 color=green!50!black,
		 only marks,
		 mark=o,
		 x filter/.code={\pgfplotstablegetelem{\coordindex}{[index]0}\of{\1}
\pgfmathtruncatemacro{\temp}{abs(\pgfplotsretval-10)==0? 1 : 0}
\ifnum\temp>0
\pgfplotstablegetelem{\coordindex}{[index]1}\of{\1}
\pgfmathtruncatemacro{\temp}{abs(\pgfplotsretval-0.62)==0? 1 : 0}
\ifnum\temp>0
\pgfplotstablegetelem{\coordindex}{[index]2}\of{\1}
\pgfmathtruncatemacro{\temp}{abs(\pgfplotsretval-0.9)==0? 1 : 0}
\ifnum\temp>0
\relax
			\else
  \def\pgfmathresult{}
\fi
\else
  \def\pgfmathresult{}
\fi
\else
  \def\pgfmathresult{}
\fi
},
		]
                table[ 
                x index = 3,
                y index = 4,
                ]
		{\1};
   \addplot[
		 color=red,
		 only marks,
		 mark=triangle,
     mark size = 2pt,
		 x filter/.code={\pgfplotstablegetelem{\coordindex}{[index]0}\of{\4}
\pgfmathtruncatemacro{\temp}{abs(\pgfplotsretval-5)==0? 1 : 0}
\ifnum\temp>0
\pgfplotstablegetelem{\coordindex}{[index]1}\of{\4}
\pgfmathtruncatemacro{\temp}{abs(\pgfplotsretval-0.3)==0? 1 : 0}
\ifnum\temp>0
\pgfplotstablegetelem{\coordindex}{[index]2}\of{\4}
\pgfmathtruncatemacro{\temp}{abs(\pgfplotsretval-0.25)==0? 1 : 0}
\ifnum\temp>0
\relax
			\else
  \def\pgfmathresult{}
\fi
\else
  \def\pgfmathresult{}
\fi
\else
  \def\pgfmathresult{}
\fi
},
		]
                table[ 
                x index = 3,
                y index = 4,
                ]
		{\4};
  \addplot[
		 color=red,
		 only marks,
     mark=o,
     mark size = 2pt,
		 x filter/.code={\pgfplotstablegetelem{\coordindex}{[index]0}\of{\1}
\pgfmathtruncatemacro{\temp}{abs(\pgfplotsretval-5)==0? 1 : 0}
\ifnum\temp>0
\pgfplotstablegetelem{\coordindex}{[index]1}\of{\1}
\pgfmathtruncatemacro{\temp}{abs(\pgfplotsretval-0.5)==0? 1 : 0}
\ifnum\temp>0
\pgfplotstablegetelem{\coordindex}{[index]2}\of{\1}
\pgfmathtruncatemacro{\temp}{abs(\pgfplotsretval-0.9)==0? 1 : 0}
\ifnum\temp>0
\relax
			\else
  \def\pgfmathresult{}
\fi
\else
  \def\pgfmathresult{}
\fi
\else
  \def\pgfmathresult{}
\fi
},
		]
                table[ 
                x index = 3,
                y index = 4,
                ]
		{\1};
                
    \end{axis}
    \end{tikzpicture}}\\
  \end{minipage}
  \caption{First $k$ mode angles, $\theta_i$  for different parameter values, as given in Table~\ref{tab:parval}. 
  }
  \label{fig:angPOD_LTSVD_CTSVD}
\end{figure}
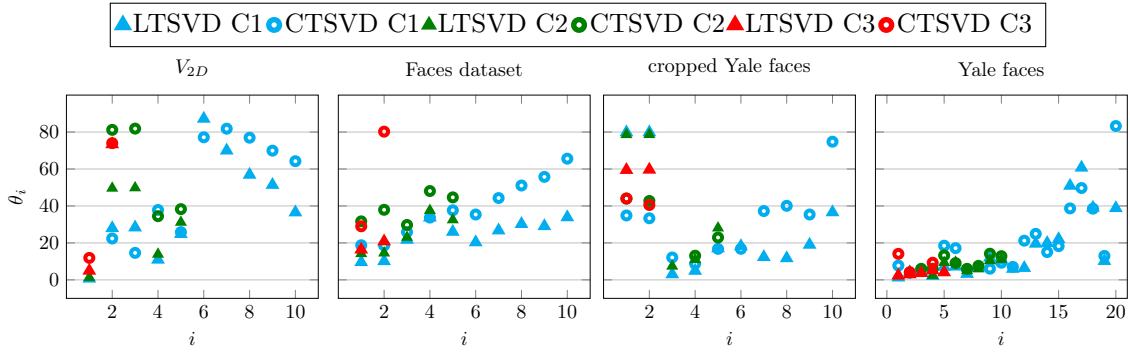

For many applications in fluid flow and in gene expression modelling, it is the POD modes rather than the singular values that are important.
Fig.~\ref{fig:POD_2} contains a few modes from the $V_{2D}$ and FACES datasets. It is seen that that there are visually apparent distortions in the approximated modes,
even though the corresponding singular values match well. 
A measure of the accuracy of the POD modes is the cosine similarity between the approximate mode ($\tilde{\bv{u}}_i$) and the one obtained using a truncated SVD of the entire dataset ($\bv{u}_i$). 
Fig.~\ref{fig:angPOD_LTSVD_CTSVD} shows the angle between the two.
We refer to these angles
as ``mode angles''. It can be seen from Fig.~\ref{fig:angPOD_LTSVD_CTSVD} the modes are approximated very poorly. For reference, the corresponding mode angle is indicated in Fig.\ref{fig:POD_2},
clearly indicating that a large mode angle implies potentially larger distortions. 

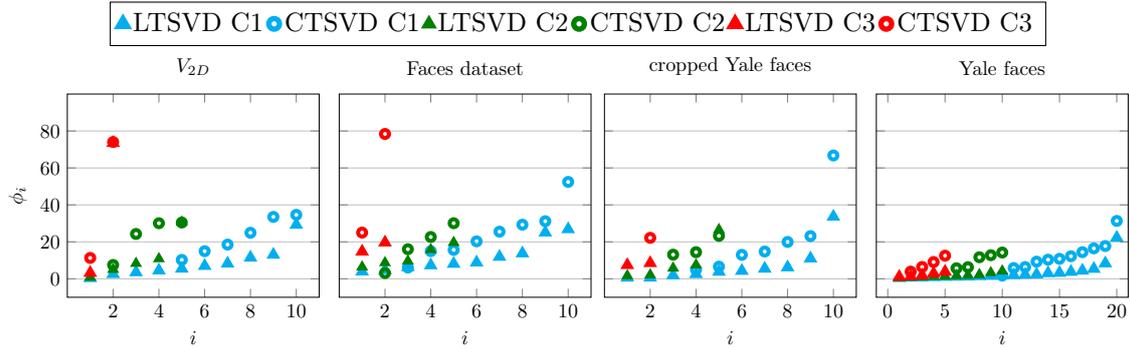
\begin{figure}[h]
  \begin{minipage}{\textwidth}
\centering
  \ref{pc_dri}\\
  \end{minipage}
  \begin{minipage}{\textwidth}
\hspace{-0.4in}    \scalebox{0.75}{  \pgfplotstableread{data/pc_LTSVD_drineas_V_2D.dat}\4
\pgfplotstableread{data/pc_CTSVD_drineas_V_2D.dat}\1

\begin{tikzpicture}[every mark/.append style={line width=2pt, solid}]
\begin{axis}[ 
  title = {$V_{2D}$}, 
  xlabel = $i$,
  ylabel = $\phi_i$,
  legend to name = pc_dri,
  legend style={at={(0.1,1.05)},anchor=south, legend columns = 6, cells={line width=2pt, solid}},
  legend entries = {LTSVD C1, CTSVD C1, LTSVD C2, CTSVD C2, LTSVD C3, CTSVD C3},
  xmax = 11.000000,
  ymax=100,
  ymin = -10,
  xtick={2,4,...,10},
  ytick = {0,20,...,80},
  width = 6cm,
  ymajorgrids,
  ]
  \addplot[
		 color=cyan,
		 only marks,
		 mark=triangle,
     mark size = 2pt,
		 x filter/.code={\pgfplotstablegetelem{\coordindex}{[index]0}\of{\4}
\pgfmathtruncatemacro{\temp}{abs(\pgfplotsretval-10)==0? 1 : 0}
\ifnum\temp>0
\pgfplotstablegetelem{\coordindex}{[index]1}\of{\4}
\pgfmathtruncatemacro{\temp}{abs(\pgfplotsretval-0.7)==0? 1 : 0}
\ifnum\temp>0
\pgfplotstablegetelem{\coordindex}{[index]2}\of{\4}
\pgfmathtruncatemacro{\temp}{abs(\pgfplotsretval-0.45)==0? 1 : 0}
\ifnum\temp>0
\relax
			\else
  \def\pgfmathresult{}
\fi
\else
  \def\pgfmathresult{}
\fi
\else
  \def\pgfmathresult{}
\fi
},
		]
                table[ 
                x index = 3,
                y index = 4,
                ]
		{\4};
  \addplot[
		 color=cyan,
		 only marks,
     mark=o,
     mark size = 2pt,
		 x filter/.code={\pgfplotstablegetelem{\coordindex}{[index]0}\of{\1}
\pgfmathtruncatemacro{\temp}{abs(\pgfplotsretval-10)==0? 1 : 0}
\ifnum\temp>0
\pgfplotstablegetelem{\coordindex}{[index]1}\of{\1}
\pgfmathtruncatemacro{\temp}{abs(\pgfplotsretval-1.04)==0? 1 : 0}
\ifnum\temp>0
\pgfplotstablegetelem{\coordindex}{[index]2}\of{\1}
\pgfmathtruncatemacro{\temp}{abs(\pgfplotsretval-0.94)==0? 1 : 0}
\ifnum\temp>0
\relax
			\else
  \def\pgfmathresult{}
\fi
\else
  \def\pgfmathresult{}
\fi
\else
  \def\pgfmathresult{}
\fi
},
		]
                table[ 
                x index = 3,
                y index = 4,
                ]
		{\1};
                \addplot[
		 color=green!50!black,
		 only marks,
		 mark=triangle,
     mark size = 1.5pt,
		 x filter/.code={\pgfplotstablegetelem{\coordindex}{[index]0}\of{\4}
\pgfmathtruncatemacro{\temp}{abs(\pgfplotsretval-5)==0? 1 : 0}
\ifnum\temp>0
\pgfplotstablegetelem{\coordindex}{[index]1}\of{\4}
\pgfmathtruncatemacro{\temp}{abs(\pgfplotsretval-1.0)==0? 1 : 0}
\ifnum\temp>0
\pgfplotstablegetelem{\coordindex}{[index]2}\of{\4}
\pgfmathtruncatemacro{\temp}{abs(\pgfplotsretval-0.1)==0? 1 : 0}
\ifnum\temp>0
\relax
			\else
  \def\pgfmathresult{}
\fi
\else
  \def\pgfmathresult{}
\fi
\else
  \def\pgfmathresult{}
\fi
},
		]
                table[ 
                x index = 3,
                y index = 4,
                ]
		{\4};
                \addplot[
		 color=green!50!black,
		 only marks,
		 mark=o,
		 x filter/.code={\pgfplotstablegetelem{\coordindex}{[index]0}\of{\1}
\pgfmathtruncatemacro{\temp}{abs(\pgfplotsretval-5)==0? 1 : 0}
\ifnum\temp>0
\pgfplotstablegetelem{\coordindex}{[index]1}\of{\1}
\pgfmathtruncatemacro{\temp}{abs(\pgfplotsretval-1.0)==0? 1 : 0}
\ifnum\temp>0
\pgfplotstablegetelem{\coordindex}{[index]2}\of{\1}
\pgfmathtruncatemacro{\temp}{abs(\pgfplotsretval-0.1)==0? 1 : 0}
\ifnum\temp>0
\relax
			\else
  \def\pgfmathresult{}
\fi
\else
  \def\pgfmathresult{}
\fi
\else
  \def\pgfmathresult{}
\fi
},
		]
                table[ 
                x index = 3,
                y index = 4,
                ]
		{\1};
   \addplot[
		 color=red,
		 only marks,
		 mark=triangle,
     mark size = 2pt,
		 x filter/.code={\pgfplotstablegetelem{\coordindex}{[index]0}\of{\4}
\pgfmathtruncatemacro{\temp}{abs(\pgfplotsretval-2)==0? 1 : 0}
\ifnum\temp>0
\pgfplotstablegetelem{\coordindex}{[index]1}\of{\4}
\pgfmathtruncatemacro{\temp}{abs(\pgfplotsretval-1.0)==0? 1 : 0}
\ifnum\temp>0
\pgfplotstablegetelem{\coordindex}{[index]2}\of{\4}
\pgfmathtruncatemacro{\temp}{abs(\pgfplotsretval-0.8)==0? 1 : 0}
\ifnum\temp>0
\relax
			\else
  \def\pgfmathresult{}
\fi
\else
  \def\pgfmathresult{}
\fi
\else
  \def\pgfmathresult{}
\fi
},
		]
                table[ 
                x index = 3,
                y index = 4,
                ]
		{\4};
  \addplot[
		 color=red,
		 only marks,
     mark=o,
     mark size = 2pt,
		 x filter/.code={\pgfplotstablegetelem{\coordindex}{[index]0}\of{\1}
\pgfmathtruncatemacro{\temp}{abs(\pgfplotsretval-2)==0? 1 : 0}
\ifnum\temp>0
\pgfplotstablegetelem{\coordindex}{[index]1}\of{\1}
\pgfmathtruncatemacro{\temp}{abs(\pgfplotsretval-1.0)==0? 1 : 0}
\ifnum\temp>0
\pgfplotstablegetelem{\coordindex}{[index]2}\of{\1}
\pgfmathtruncatemacro{\temp}{abs(\pgfplotsretval-0.8)==0? 1 : 0}
\ifnum\temp>0
\relax
			\else
  \def\pgfmathresult{}
\fi
\else
  \def\pgfmathresult{}
\fi
\else
  \def\pgfmathresult{}
\fi
},
		]
                table[ 
                x index = 3,
                y index = 4,
                ]
		{\1};
    \end{axis}
    \end{tikzpicture}}
 \hspace{-0.22in}    \scalebox{0.75}{ \pgfplotstableread{data/pc_LTSVD_drineas_faces.dat}\4
\pgfplotstableread{data/pc_CTSVD_drineas_faces.dat}\1

\begin{tikzpicture}[every mark/.append style={line width=2pt, solid}]
\begin{axis}[ 
  title = {Faces dataset}, 
  xlabel = $i$,
  xmax = 11.000000,
  ymax=100,
  ymin = -10,
  xtick={2,4,...,10},
  ytick = {0,20,...,80},
  yticklabels={,,},
  width = 6cm,
  ymajorgrids,
  ]
  \addplot[
		 color=cyan,
		 only marks,
		 mark=triangle,
     mark size = 2pt,
		 x filter/.code={\pgfplotstablegetelem{\coordindex}{[index]0}\of{\4}
\pgfmathtruncatemacro{\temp}{abs(\pgfplotsretval-10)==0? 1 : 0}
\ifnum\temp>0
\pgfplotstablegetelem{\coordindex}{[index]1}\of{\4}
\pgfmathtruncatemacro{\temp}{abs(\pgfplotsretval-0.75)==0? 1 : 0}
\ifnum\temp>0
\pgfplotstablegetelem{\coordindex}{[index]2}\of{\4}
\pgfmathtruncatemacro{\temp}{abs(\pgfplotsretval-0.8)==0? 1 : 0}
\ifnum\temp>0
\relax
			\else
  \def\pgfmathresult{}
\fi
\else
  \def\pgfmathresult{}
\fi
\else
  \def\pgfmathresult{}
\fi
},
		]
                table[ 
                x index = 3,
                y index = 4,
                ]
		{\4};
  \addplot[
		 color=cyan,
		 only marks,
     mark=o,
     mark size = 2pt,
		 x filter/.code={\pgfplotstablegetelem{\coordindex}{[index]0}\of{\1}
\pgfmathtruncatemacro{\temp}{abs(\pgfplotsretval-10)==0? 1 : 0}
\ifnum\temp>0
\pgfplotstablegetelem{\coordindex}{[index]1}\of{\1}
\pgfmathtruncatemacro{\temp}{abs(\pgfplotsretval-1.3)==0? 1 : 0}
\ifnum\temp>0
\pgfplotstablegetelem{\coordindex}{[index]2}\of{\1}
\pgfmathtruncatemacro{\temp}{abs(\pgfplotsretval-1)==0? 1 : 0}
\ifnum\temp>0
\relax
			\else
  \def\pgfmathresult{}
\fi
\else
  \def\pgfmathresult{}
\fi
\else
  \def\pgfmathresult{}
\fi
},
		]
                table[ 
                x index = 3,
                y index = 4,
                ]
		{\1};
                \addplot[
		 color=green!50!black,
		 only marks,
		 mark=triangle,
     mark size = 1.5pt,
		 x filter/.code={\pgfplotstablegetelem{\coordindex}{[index]0}\of{\4}
\pgfmathtruncatemacro{\temp}{abs(\pgfplotsretval-5)==0? 1 : 0}
\ifnum\temp>0
\pgfplotstablegetelem{\coordindex}{[index]1}\of{\4}
\pgfmathtruncatemacro{\temp}{abs(\pgfplotsretval-0.75)==0? 1 : 0}
\ifnum\temp>0
\pgfplotstablegetelem{\coordindex}{[index]2}\of{\4}
\pgfmathtruncatemacro{\temp}{abs(\pgfplotsretval-0.75)==0? 1 : 0}
\ifnum\temp>0
\relax
			\else
  \def\pgfmathresult{}
\fi
\else
  \def\pgfmathresult{}
\fi
\else
  \def\pgfmathresult{}
\fi
},
		]
                table[ 
                x index = 3,
                y index = 4,
                ]
		{\4};
                \addplot[
		 color=green!50!black,
		 only marks,
		 mark=o,
		 x filter/.code={\pgfplotstablegetelem{\coordindex}{[index]0}\of{\1}
\pgfmathtruncatemacro{\temp}{abs(\pgfplotsretval-5)==0? 1 : 0}
\ifnum\temp>0
\pgfplotstablegetelem{\coordindex}{[index]1}\of{\1}
\pgfmathtruncatemacro{\temp}{abs(\pgfplotsretval-1.0)==0? 1 : 0}
\ifnum\temp>0
\pgfplotstablegetelem{\coordindex}{[index]2}\of{\1}
\pgfmathtruncatemacro{\temp}{abs(\pgfplotsretval-1)==0? 1 : 0}
\ifnum\temp>0
\relax
			\else
  \def\pgfmathresult{}
\fi
\else
  \def\pgfmathresult{}
\fi
\else
  \def\pgfmathresult{}
\fi
},
		]
                table[ 
                x index = 3,
                y index = 4,
                ]
		{\1};
   \addplot[
		 color=red,
		 only marks,
		 mark=triangle,
     mark size = 2pt,
		 x filter/.code={\pgfplotstablegetelem{\coordindex}{[index]0}\of{\4}
\pgfmathtruncatemacro{\temp}{abs(\pgfplotsretval-2)==0? 1 : 0}
\ifnum\temp>0
\pgfplotstablegetelem{\coordindex}{[index]1}\of{\4}
\pgfmathtruncatemacro{\temp}{abs(\pgfplotsretval-1.0)==0? 1 : 0}
\ifnum\temp>0
\pgfplotstablegetelem{\coordindex}{[index]2}\of{\4}
\pgfmathtruncatemacro{\temp}{abs(\pgfplotsretval-0.8)==0? 1 : 0}
\ifnum\temp>0
\relax
			\else
  \def\pgfmathresult{}
\fi
\else
  \def\pgfmathresult{}
\fi
\else
  \def\pgfmathresult{}
\fi
},
		]
                table[ 
                x index = 3,
                y index = 4,
                ]
		{\4};
  \addplot[
		 color=red,
		 only marks,
     mark=o,
     mark size = 2pt,
		 x filter/.code={\pgfplotstablegetelem{\coordindex}{[index]0}\of{\1}
\pgfmathtruncatemacro{\temp}{abs(\pgfplotsretval-2)==0? 1 : 0}
\ifnum\temp>0
\pgfplotstablegetelem{\coordindex}{[index]1}\of{\1}
\pgfmathtruncatemacro{\temp}{abs(\pgfplotsretval-1.0)==0? 1 : 0}
\ifnum\temp>0
\pgfplotstablegetelem{\coordindex}{[index]2}\of{\1}
\pgfmathtruncatemacro{\temp}{abs(\pgfplotsretval-0.8)==0? 1 : 0}
\ifnum\temp>0
\relax
			\else
  \def\pgfmathresult{}
\fi
\else
  \def\pgfmathresult{}
\fi
\else
  \def\pgfmathresult{}
\fi
},
		]
                table[ 
                x index = 3,
                y index = 4,
                ]
		{\1};
                
    \end{axis}
    \end{tikzpicture}}
\hspace{-0.2in}  \scalebox{0.75}{\pgfplotstableread{data/pc_LTSVD_drineas_croppedYalefaces.dat}\4
\pgfplotstableread{data/pc_CTSVD_drineas_croppedYalefaces.dat}\1

\begin{tikzpicture}[every mark/.append style={line width=2pt, solid}]
\begin{axis}[
    title = {cropped Yale faces}, 
    xlabel = $i$,
    xmax = 11.000000,
    ymax = 100,
    ymin = -10,
    xtick={2,4,...,10},
    ytick = {0,20,...,80},
    yticklabels={,,},
    width = 6cm,
    ymajorgrids,
  ]
  \addplot[
		 color=cyan,
		 only marks,
		 mark=triangle,
     mark size = 2pt,
		 x filter/.code={\pgfplotstablegetelem{\coordindex}{[index]0}\of{\4}
\pgfmathtruncatemacro{\temp}{abs(\pgfplotsretval-10)==0? 1 : 0}
\ifnum\temp>0
\pgfplotstablegetelem{\coordindex}{[index]1}\of{\4}
\pgfmathtruncatemacro{\temp}{abs(\pgfplotsretval-0.46)==0? 1 : 0}
\ifnum\temp>0
\pgfplotstablegetelem{\coordindex}{[index]2}\of{\4}
\pgfmathtruncatemacro{\temp}{abs(\pgfplotsretval-0.44)==0? 1 : 0}
\ifnum\temp>0
\relax
			\else
  \def\pgfmathresult{}
\fi
\else
  \def\pgfmathresult{}
\fi
\else
  \def\pgfmathresult{}
\fi
},
		]
                table[ 
                x index = 3,
                y index = 4,
                ]
		{\4};
  \addplot[
		 color=cyan,
		 only marks,
     mark=o,
     mark size = 2pt,
		 x filter/.code={\pgfplotstablegetelem{\coordindex}{[index]0}\of{\1}
\pgfmathtruncatemacro{\temp}{abs(\pgfplotsretval-10)==0? 1 : 0}
\ifnum\temp>0
\pgfplotstablegetelem{\coordindex}{[index]1}\of{\1}
\pgfmathtruncatemacro{\temp}{abs(\pgfplotsretval-0.9)==0? 1 : 0}
\ifnum\temp>0
\pgfplotstablegetelem{\coordindex}{[index]2}\of{\1}
\pgfmathtruncatemacro{\temp}{abs(\pgfplotsretval-0.7)==0? 1 : 0}
\ifnum\temp>0
\relax
			\else
  \def\pgfmathresult{}
\fi
\else
  \def\pgfmathresult{}
\fi
\else
  \def\pgfmathresult{}
\fi
},
		]
                table[ 
                x index = 3,
                y index = 4,
                ]
		{\1};
                \addplot[
		 color=green!50!black,
		 only marks,
		 mark=triangle,
     mark size = 1.5pt,
		 x filter/.code={\pgfplotstablegetelem{\coordindex}{[index]0}\of{\4}
\pgfmathtruncatemacro{\temp}{abs(\pgfplotsretval-5)==0? 1 : 0}
\ifnum\temp>0
\pgfplotstablegetelem{\coordindex}{[index]1}\of{\4}
\pgfmathtruncatemacro{\temp}{abs(\pgfplotsretval-1.0)==0? 1 : 0}
\ifnum\temp>0
\pgfplotstablegetelem{\coordindex}{[index]2}\of{\4}
\pgfmathtruncatemacro{\temp}{abs(\pgfplotsretval-0.1)==0? 1 : 0}
\ifnum\temp>0
\relax
			\else
  \def\pgfmathresult{}
\fi
\else
  \def\pgfmathresult{}
\fi
\else
  \def\pgfmathresult{}
\fi
},
		]
                table[ 
                x index = 3,
                y index = 4,
                ]
		{\4};
                \addplot[
		 color=green!50!black,
		 only marks,
		 mark=o,
		 x filter/.code={\pgfplotstablegetelem{\coordindex}{[index]0}\of{\1}
\pgfmathtruncatemacro{\temp}{abs(\pgfplotsretval-5)==0? 1 : 0}
\ifnum\temp>0
\pgfplotstablegetelem{\coordindex}{[index]1}\of{\1}
\pgfmathtruncatemacro{\temp}{abs(\pgfplotsretval-1.0)==0? 1 : 0}
\ifnum\temp>0
\pgfplotstablegetelem{\coordindex}{[index]2}\of{\1}
\pgfmathtruncatemacro{\temp}{abs(\pgfplotsretval-0.1)==0? 1 : 0}
\ifnum\temp>0
\relax
			\else
  \def\pgfmathresult{}
\fi
\else
  \def\pgfmathresult{}
\fi
\else
  \def\pgfmathresult{}
\fi
},
		]
                table[ 
                x index = 3,
                y index = 4,
                ]
		{\1};
   \addplot[
		 color=red,
		 only marks,
		 mark=triangle,
     mark size = 2pt,
		 x filter/.code={\pgfplotstablegetelem{\coordindex}{[index]0}\of{\4}
\pgfmathtruncatemacro{\temp}{abs(\pgfplotsretval-2)==0? 1 : 0}
\ifnum\temp>0
\pgfplotstablegetelem{\coordindex}{[index]1}\of{\4}
\pgfmathtruncatemacro{\temp}{abs(\pgfplotsretval-1.0)==0? 1 : 0}
\ifnum\temp>0
\pgfplotstablegetelem{\coordindex}{[index]2}\of{\4}
\pgfmathtruncatemacro{\temp}{abs(\pgfplotsretval-0.8)==0? 1 : 0}
\ifnum\temp>0
\relax
			\else
  \def\pgfmathresult{}
\fi
\else
  \def\pgfmathresult{}
\fi
\else
  \def\pgfmathresult{}
\fi
},
		]
                table[ 
                x index = 3,
                y index = 4,
                ]
		{\4};
  \addplot[
		 color=red,
		 only marks,
     mark=o,
     mark size = 2pt,
		 x filter/.code={\pgfplotstablegetelem{\coordindex}{[index]0}\of{\1}
\pgfmathtruncatemacro{\temp}{abs(\pgfplotsretval-2)==0? 1 : 0}
\ifnum\temp>0
\pgfplotstablegetelem{\coordindex}{[index]1}\of{\1}
\pgfmathtruncatemacro{\temp}{abs(\pgfplotsretval-1.0)==0? 1 : 0}
\ifnum\temp>0
\pgfplotstablegetelem{\coordindex}{[index]2}\of{\1}
\pgfmathtruncatemacro{\temp}{abs(\pgfplotsretval-0.8)==0? 1 : 0}
\ifnum\temp>0
\relax
			\else
  \def\pgfmathresult{}
\fi
\else
  \def\pgfmathresult{}
\fi
\else
  \def\pgfmathresult{}
\fi
},
		]
                table[ 
                x index = 3,
                y index = 4,
                ]
		{\1};
                
    \end{axis}
    \end{tikzpicture}}
\hspace{-0.22in}     \scalebox{0.75}{ \pgfplotstableread{data/pc_LTSVD_drineas_ext_Yalefaces_wo_Ambient.dat}\4
\pgfplotstableread{data/pc_CTSVD_drineas_ext_Yalefaces_wo_Ambient.dat}\1

\begin{tikzpicture}[every mark/.append style={line width=2pt, solid}]
\begin{axis}[
    title = {Yale faces}, 
    xlabel = $i$,
    ymax=100,
    ymin = -10,
    xmax = 21.000000,
    ytick = {0,20,...,80},
    yticklabels={,,},
    width = 6cm,
    ymajorgrids,
  ]
  \addplot[
		 color=cyan,
		 only marks,
		 mark=triangle,
     mark size = 2pt,
		 x filter/.code={\pgfplotstablegetelem{\coordindex}{[index]0}\of{\4}
\pgfmathtruncatemacro{\temp}{abs(\pgfplotsretval-20)==0? 1 : 0}
\ifnum\temp>0
\pgfplotstablegetelem{\coordindex}{[index]1}\of{\4}
\pgfmathtruncatemacro{\temp}{abs(\pgfplotsretval-0.35)==0? 1 : 0}
\ifnum\temp>0
\pgfplotstablegetelem{\coordindex}{[index]2}\of{\4}
\pgfmathtruncatemacro{\temp}{abs(\pgfplotsretval-0.35)==0? 1 : 0}
\ifnum\temp>0
\relax
			\else
  \def\pgfmathresult{}
\fi
\else
  \def\pgfmathresult{}
\fi
\else
  \def\pgfmathresult{}
\fi
},
		]
                table[ 
                x index = 3,
                y index = 4,
                ]
		{\4};
  \addplot[
		 color=cyan,
		 only marks,
     mark=o,
     mark size = 2pt,
		 x filter/.code={\pgfplotstablegetelem{\coordindex}{[index]0}\of{\1}
\pgfmathtruncatemacro{\temp}{abs(\pgfplotsretval-20)==0? 1 : 0}
\ifnum\temp>0
\pgfplotstablegetelem{\coordindex}{[index]1}\of{\1}
\pgfmathtruncatemacro{\temp}{abs(\pgfplotsretval-0.83)==0? 1 : 0}
\ifnum\temp>0
\pgfplotstablegetelem{\coordindex}{[index]2}\of{\1}
\pgfmathtruncatemacro{\temp}{abs(\pgfplotsretval-0.9)==0? 1 : 0}
\ifnum\temp>0
\relax
			\else
  \def\pgfmathresult{}
\fi
\else
  \def\pgfmathresult{}
\fi
\else
  \def\pgfmathresult{}
\fi
},
		]
                table[ 
                x index = 3,
                y index = 4,
                ]
		{\1};
                \addplot[
		 color=green!50!black,
		 only marks,
		 mark=triangle,
     mark size = 1.5pt,
		 x filter/.code={\pgfplotstablegetelem{\coordindex}{[index]0}\of{\4}
\pgfmathtruncatemacro{\temp}{abs(\pgfplotsretval-10)==0? 1 : 0}
\ifnum\temp>0
\pgfplotstablegetelem{\coordindex}{[index]1}\of{\4}
\pgfmathtruncatemacro{\temp}{abs(\pgfplotsretval-0.3)==0? 1 : 0}
\ifnum\temp>0
\pgfplotstablegetelem{\coordindex}{[index]2}\of{\4}
\pgfmathtruncatemacro{\temp}{abs(\pgfplotsretval-0.25)==0? 1 : 0}
\ifnum\temp>0
\relax
			\else
  \def\pgfmathresult{}
\fi
\else
  \def\pgfmathresult{}
\fi
\else
  \def\pgfmathresult{}
\fi
},
		]
                table[ 
                x index = 3,
                y index = 4,
                ]
		{\4};
                \addplot[
		 color=green!50!black,
		 only marks,
		 mark=o,
		 x filter/.code={\pgfplotstablegetelem{\coordindex}{[index]0}\of{\1}
\pgfmathtruncatemacro{\temp}{abs(\pgfplotsretval-10)==0? 1 : 0}
\ifnum\temp>0
\pgfplotstablegetelem{\coordindex}{[index]1}\of{\1}
\pgfmathtruncatemacro{\temp}{abs(\pgfplotsretval-0.62)==0? 1 : 0}
\ifnum\temp>0
\pgfplotstablegetelem{\coordindex}{[index]2}\of{\1}
\pgfmathtruncatemacro{\temp}{abs(\pgfplotsretval-0.9)==0? 1 : 0}
\ifnum\temp>0
\relax
			\else
  \def\pgfmathresult{}
\fi
\else
  \def\pgfmathresult{}
\fi
\else
  \def\pgfmathresult{}
\fi
},
		]
                table[ 
                x index = 3,
                y index = 4,
                ]
		{\1};
   \addplot[
		 color=red,
		 only marks,
		 mark=triangle,
     mark size = 2pt,
		 x filter/.code={\pgfplotstablegetelem{\coordindex}{[index]0}\of{\4}
\pgfmathtruncatemacro{\temp}{abs(\pgfplotsretval-5)==0? 1 : 0}
\ifnum\temp>0
\pgfplotstablegetelem{\coordindex}{[index]1}\of{\4}
\pgfmathtruncatemacro{\temp}{abs(\pgfplotsretval-0.3)==0? 1 : 0}
\ifnum\temp>0
\pgfplotstablegetelem{\coordindex}{[index]2}\of{\4}
\pgfmathtruncatemacro{\temp}{abs(\pgfplotsretval-0.25)==0? 1 : 0}
\ifnum\temp>0
\relax
			\else
  \def\pgfmathresult{}
\fi
\else
  \def\pgfmathresult{}
\fi
\else
  \def\pgfmathresult{}
\fi
},
		]
                table[ 
                x index = 3,
                y index = 4,
                ]
		{\4};
  \addplot[
		 color=red,
		 only marks,
     mark=o,
     mark size = 2pt,
		 x filter/.code={\pgfplotstablegetelem{\coordindex}{[index]0}\of{\1}
\pgfmathtruncatemacro{\temp}{abs(\pgfplotsretval-5)==0? 1 : 0}
\ifnum\temp>0
\pgfplotstablegetelem{\coordindex}{[index]1}\of{\1}
\pgfmathtruncatemacro{\temp}{abs(\pgfplotsretval-0.5)==0? 1 : 0}
\ifnum\temp>0
\pgfplotstablegetelem{\coordindex}{[index]2}\of{\1}
\pgfmathtruncatemacro{\temp}{abs(\pgfplotsretval-0.9)==0? 1 : 0}
\ifnum\temp>0
\relax
			\else
  \def\pgfmathresult{}
\fi
\else
  \def\pgfmathresult{}
\fi
\else
  \def\pgfmathresult{}
\fi
},
		]
                table[ 
                x index = 3,
                y index = 4,
                ]
		{\1};
                
    \end{axis}
    \end{tikzpicture}}
  \end{minipage}
  \caption{First $k$ principal angles, $\phi_i$, for different datasets for different parameter values, as given in Table~\ref{tab:parval}. 
  }
  \label{fig:pc_LTSVD_CTSVD}
\end{figure}

The mode angles are known to be sensitive to clustering of singular values. This is the reason why the first two mode angles are large for CYF ($\sigma_1\sim \sigma_2$).  
In the limiting case when the singular values are identical, it is only the eigenspace that matters and not the eigenvectors. Also, for applications such as facial recognition, it is the space spanned by the
      the modes rather than the modes themselves that are of interest. 
      For this reason, we also measure the accuracy of the subspaces spanned by $k$ singular vectors rather than the accuracy of each mode. A measure of this accuracy is the principal or canonical angles between the subspaces, whose cosines are the singular values of $\tilde{U}_k^TU_k$ \cite{BjoGol:1973}. 
      Fig.~\ref{fig:pc_LTSVD_CTSVD} shows the $k$ principal angles, $\phi_i$, between the two $k$-dimensional subspaces. Although the
      accuracy of the subspace is better than the modes themselves, in some cases, there is
      significant error even for $c\approx n$.
      The first two principal angles for CYF match well, indicating that the subspace is captured more accurately than the modes themselves (due to almost identical singular values). 


\bibliography{paper}
\bibliographystyle{unsrt}
\end{document}